\documentclass[a4paper]{amsart}
\makeatletter
\RequirePackage{ifthen}
\provideboolean{omarfont}
\setboolean{omarfont}{false}
\provideboolean{usesvjour}
\setboolean{usesvjour}{false}%
\provideboolean{nohyperref}
\setboolean{nohyperref}{false}%
\provideboolean{usecleveref}
\setboolean{usecleveref}{false}%
\provideboolean{landscape}
\setboolean{landscape}{false}%
\makeatletter
      \RequirePackage{ifthen}
      \provideboolean{useutopia}
      \setboolean{useutopia}{false}%
      \provideboolean{usesvjour}
      \setboolean{usesvjour}{false}%
      \provideboolean{usebeamer}%
      \setboolean{usebeamer}{false}%
      \provideboolean{usinghyperref}
      \setboolean{usinghyperref}{false}%
\makeatletter
          \RequirePackage{xparse}
          \RequirePackage{ifthen}
          \RequirePackage{iftex}
          \provideboolean{usemathrsfs}%
          \provideboolean{useutopia}
          \setboolean{useutopia}{false}%
          \provideboolean{usebeamer}%
          \setboolean{usebeamer}{false}%
          \provideboolean{isthesis}%
          \provideboolean{isamsltex}%
          \provideboolean{issiamltex}%
          \ifXeTeX%
            \RequirePackage{xltxtra}%
            \ifthenelse{2=1}{
              \RequirePackage{polyglossia}
              \setdefaultlanguage[variant=british]{english}
              \setotherlanguages{french,vietnamese,russian}%
            }{
              \RequirePackage[american,main=british,french,vietnamese]{babel}
            }
            \DeclareSymbolFont{usualmathcal}{OMS}{cmsy}{m}{n}
            \DeclareSymbolFontAlphabet{\mathcalbf}{usualmathcal}
            \ifthenelse{\boolean{useutopia}}{
              \RequirePackage[utopia,expert,euro]{mathdesign}
              \RequirePackage[OMLmathbf,OMLmathsfit]{isomath}
            }{
              \RequirePackage{amssymb}    %
              \RequirePackage{bbold}    %
            }
            \RequirePackage{mathrsfs} %
            \setboolean{usemathrsfs}{true}%
          \else%
          \RequirePackage[utf8]{inputenc}%
          \ifthenelse{\boolean{useutopia}}{%
            \RequirePackage[utopia,euro]{mathdesign}%
            \RequirePackage[OMLmathrm,OMLmathbf,OMLmathsfit]{isomath}
            \RequirePackage{bbold}    %
            \DeclareSymbolFont{usualmathcal}{OMS}{cmsy}{m}{n}
            \DeclareSymbolFontAlphabet{\mathcalbf}{usualmathcal}
            \RequirePackage{stmaryrd} %
            \providecommand{\diracdelta}[1][]{\ensuremath{\deltaup_{#1}}}
            
            \providecommand{\lap}{\ensuremath{\Deltaup}}
            \providecommand{\pic}{\ensuremath{\mathrm\pi}}
            \providecommand{\measure}[1]{\ensuremath{\mathcalbf{\uppercase{#1}}}}
          }{%
            \RequirePackage{amssymb}    %
            \RequirePackage{mathrsfs} %
            \RequirePackage{bbold}    %
            \RequirePackage{stmaryrd} %
            \setboolean{usemathrsfs}{true}%
            \providecommand{\mathcalbf}{\mathcal}
          }%
          \RequirePackage[american,main=british,vietnamese]{babel}
          \fi%
          \RequirePackage{xstring}%
          \RequirePackage{chngcntr}%
          \RequirePackage{mathtools}%
          \RequirePackage{nicefrac}
          \RequirePackage{stmaryrd} %
          \RequirePackage{amsthm}%
          \ifthenelse{\boolean{usebeamer}}{}{
            \RequirePackage[shortlabels]{enumitem}%
          }
          \RequirePackage{xspace}%
          \RequirePackage{verbatim}%
          \RequirePackage{fvextra}
          \RequirePackage{fancyvrb}
          \RequirePackage{listings}%
          \RequirePackage{xcolor}
          \RequirePackage{epsfig}
          \RequirePackage{graphicx}
          \RequirePackage[space]{grffile}%
          \RequirePackage{booktabs}%
          \RequirePackage{tikz}%
          \usetikzlibrary{calc}%
          \usetikzlibrary{fadings}
      \def\olprovideenvironment{\@star@or@long\provide@environment}
      \def\provide@environment#1{%
              \@ifundefined{#1}%
                      {\def\reserved@a{\newenvironment{#1}}}%
                      {\def\reserved@a{\renewenvironment{dummy@environ}}}%
              \reserved@a
      }
      \def\dummy@environ{}
      \colorlet{a}{magenta}
      \colorlet{b}{green!75!blue}
      \colorlet{c}{yellow!87.5!red}
      \colorlet{d}{cyan}
      \colorlet{e}{red}
      \colorlet{f}{blue}
      \colorlet{g}{white}
      \colorlet{i}{black}
      \colorlet{h}{i!50!g}
      \colorlet{j}{a!75!g}

      \ifthenelse{\boolean{usinghyperref}}{%
        \providecommand{\linkedurl}[1]{\url{1}}%
        \providecommand{\linkedemail}[1]{\href{mailto:#1}{#1}}%

      }{%
        \providecommand{\linkedurl}[1]{\texttt{#1}}%
        \providecommand{\linkedemail}[1]{\texttt{#1}}%
      }
      \providecommand{\email}[1]{{\linkedemail{#1}}}
      \providecommand{\Ignore}[1]{}
      \providecommand{\ignore}[1]{}
      \providecommand{\freeze}[1]{}%
      
      \providecommand{\crossout}[1]{{\color{i!20} #1}}
      \providecommand{\highlightcolor}{a}
      \providecommand{\highlight}[1]{{\color{\highlightcolor}#1}}

      \providecommand{\memo}[1]{%
        \ensuremath{%
          \framebox{\tiny\textbf{\kern-2pt\textsf{#1}}\kern-2pt}%
        }%
        \xspace}

      \RequirePackage{alphalph}
      \provideboolean{shownotes}
      \setboolean{shownotes}{true}%
      \newcounter{margnote}[page]
      \providecommand{\mgcolor}{a}%
      \providecommand{\mgcolorset}[1]{\renewcommand{\mgcolor}{\alphalph{#1}}}
      \providecommand{\mgcolorsetbycounter}[1]{%
        \newcount\@olmodn
        \@olmodn \value{#1}\relax
        \newcount\@olmodd
        \@olmodd 6\relax
        \newcount\@olmodq %
        \newcount\@olmodr %
        \newcount\@olmodc %
        \@olmodc\@olmodd\relax
        \multiply\@olmodc by 2\relax
        \advance\@olmodc-1\relax
        \@olmodq 0\relax
        \@olmodr\@olmodn\relax
        \ifnum\@olmodr>\@olmodc
        \loop
        \advance\@olmodr-\@olmodd\relax
        \advance\@olmodq1\relax
        \ifnum\@olmodr>\@olmodc
        \repeat
        \fi
        \setcounter{tmpcounter}{\the\@olmodr}
        \stepcounter{tmpcounter}
        \mgcolorset{\value{tmpcounter}}
      }
      \providecommand{\mgcolormake}{\mgcolorsetbycounter{margnote}}
      \providecommand{\margnotecolor}{\mgcolormake}
      \providecommand{\margnotemark}{{\colorbox{\mgcolor}{\tiny\color{g}\upshape\texttt{\arabic{page}.\arabic{margnote}}}}\,}
      \providecommand{\margnote}[2][]{%
        \ifthenelse{%
          \boolean{shownotes}%
        }{%
          \stepcounter{margnote}%
          \margnotecolor%
          \margnotemark %
          \marginpar{%
            \color{\mgcolor}%
            \texttt{\bfseries{%
              \begin{minipage}{2cm}%
                \raggedright\tiny%
                \margnotemark%
                #2%
                \\
                {\ifx|#1|{}\else{ - #1}\fi}%
              \end{minipage}%
              }%
            }%
          }%
        }{%
        }%
      }%
      \providecommand{\mathnote}[2][]{%
        \ifthenelse{%
          \boolean{shownotes}%
        }{%
          \stepcounter{margnote}%
          \margnotecolor%
          \text{%
            \colorbox{g}{%
              \color{\mgcolor}%
              \texttt{\bfseries{%
                \tiny%
                    \margnotemark\,%
                    #2%
                    \ifx|#1|{}\else{ - #1}\fi%
                }%
              }%
            }%
          }%
        }{%
        }%
      }%
      \providecommand{\textnote}[2][]{%
        \ifthenelse{%
          \boolean{shownotes}%
        }{%
          \stepcounter{margnote}%
          \margnotecolor%
          \ \\
          \text{%
              \begin{minipage}{\textwidth}
              \color{\mgcolor}%
              \texttt{%
                \margnotemark%
                #2%
                \ifx|#1|{}\else{ - #1}\fi%
              }%
              \end{minipage}
          }%
        }{}%
      }%

      \providecommand{\Todo}[2][To do:]{
        \ifthenelse{\boolean{shownotes}}{
          \begin{tikzpicture}
           \node[fill=a!17]{
             \begin{minipage}{\textwidth}
               \tiny
               \texttt{#1}
               \texttt{\bfseries{#2}}
             \end{minipage}
           };
          \end{tikzpicture}
        }{}}
      \provideboolean{showrevisions}
      \setboolean{showrevisions}{true}
      \provideboolean{emphrevisions}
      \setboolean{emphrevisions}{false}
      \newcommand{\revisionsheader}{\ \clearpage\Warning{the following part is under development/revision}}
      \newcommand{\revisionsfooter}{\ \newline\Warning{end of part under development/revision}\clearpage}

      \providecommand{\HighlightBox}[2][a!6.25]{
        \begin{center}
          \begin{tikzpicture}
            \node[fill=#1]{
              \begin{minipage}{\textwidth}
                #2
              \end{minipage}
            };
          \end{tikzpicture}
        \end{center}
      }
      \providecommand{\Warning}[1]{    
        \HighlightBox[b!25]{%
          \texttt{\bfseries{\small Warning: #1}}
        }
      }

      \provideboolean{showcomments}
      \providecommand{\margincomment}[1]{
      \ifthenelse{\boolean{showcomments}}{\marginpar{\tiny #1}}{}
      }
      \provideboolean{showchanges}
      \setboolean{showchanges}{false}
      \providecommand{\changes}[2][]{%
        \ifthenelse{\boolean{showchanges}}{{%
            \ifx|#1|{}\else\margnote{#1}\fi%
            \highlight{#2}%
        }}{%
          #2}}
      \providecommand{\mathchanges}[2][]{%
        \ifthenelse{\boolean{showchanges}}{{\ifx|#1|{}\else\mathnote{#1}\fi\highlight{#2}}}{#2}}

      \providecommand{\changefromto}[3][replace with]{%
        \ifthenelse{\boolean{showchanges}}{{%
            \crossout{#2}\margnote{#1}%
          }{%
            \highlight{#3}
          }%
        }{%
          #3\xspace%
        }%
      }
      \providecommand{\ChangePar}[3][]{%
        \ifthenelse{\boolean{showchanges}}{
          {\par\textcolor{i!20}{#2}\ifx|#1|\else\margnote{#1}\fi}{\par\textcolor{a}{#3}}
        }{%
          \par #3%
        }%
      }
      \providecommand{\InsertPar}[1]{
        \ifthenelse{\boolean{showchanges}}
        {{\par$\mapsto$ \textcolor{blue}{#1}}}
        {\par #1}
      }
      
      \providecommand{\mathchangefromto}[3][]{\crossout{#2}\ifx|#1|\else\mathnote{#1}\fi\highlight{#3}}

      \let\trueMakeUppercase\MakeUppercase
      \newcommand{\UCmath}[1]{%
        \begingroup
        \ucmathlist\trueMakeUppercase{#1}%
        \endgroup
      }
      \ifthenelse{\(\boolean{useutopia}\)\OR\(\boolean{usesvjour}\)}{
        \newcommand{\ucmathlist}{%
          \def\alpha{A}%
          \def\beta{B}%
          \let\gamma\Gamma
          \let\delta\Delta
          \def\epsilon{E}%
          \def\varepsilon{E}%
          \def\zeta{Z}%
          \def\eta{H}%
          \let\theta\Theta
          \let\vartheta\Theta
          \def\iota{I}%
          \def\kappa{K}%
          \let\lambda\Lambda
          \def\mu{M}%
          \def\nu{N}%
          \let\xi\Xi
          \def\omicron{O}
          \let\pi\Pi
          \let\varpi\Pi
          \def\rho{P}%
          \def\varrho{P}%
          \let\sigma\Sigma
          \def\varsigma{C}
          \def\tau{T}%
          \let\upsilon\Upsilon
          \let\phi\Phi
          \let\varphi\Phi
          \def\chi{X}%
          \let\psi\Psi
          \let\omega\Omega
      }}{
        \newcommand{\ucmathlist}{
          \def\alpha{\mathrm{A}}%
          \def\beta{\mathrm{B}}%
          \let\gamma\Gamma
          \let\delta\Delta
          \def\epsilon{\mathrm{E}}%
          \def\varepsilon{\mathrm{E}}%
          \def\zeta{\mathrm{Z}}%
          \def\eta{\mathrm{H}}%
          \let\theta\Theta
          \let\vartheta\Theta
          \def\iota{\mathrm{I}}%
          \def\kappa{\mathrm{K}}%
          \let\lambda\Lambda
          \def\mu{\mathrm{M}}%
          \def\nu{\mathrm{N}}%
          \let\xi\Xi
          \let\pi\Pi
          \let\varpi\Pi
          \def\rho{\mathrm{P}}%
          \def\varrho{\mathrm{P}}%
          \let\sigma\Sigma
          \def\tau{\mathrm{T}}%
          \let\upsilon\Upsilon
          \let\phi\Phi
          \let\varphi\Phi
          \def\chi{\mathrm{X}}%
          \let\psi\Psi
          \let\omega\Omega
        }
      }
      \providecommand{\mathscript}
      	       {\mathscr}

      \providecommand{\cA}{\ensuremath{\mathscript A}\xspace}

      \providecommand{\cD}{\ensuremath{\mathscript D}\xspace}
      \providecommand{\cE}{\ensuremath{\mathscript E}\xspace}

      \providecommand{\cI}{\ensuremath{\mathscript I}\xspace}
      \providecommand{\cJ}{\ensuremath{\mathscript J}\xspace}

      \providecommand{\cM}{\ensuremath{\mathscript M}\xspace}

      \providecommand{\cR}{\ensuremath{\mathscript R}\xspace}
      \providecommand{\cS}{\ensuremath{\mathscript S}\xspace}
      \providecommand{\cT}{\ensuremath{\mathscript T}\xspace}

      \providecommand{\bbbold}{\mathbb}

      \providecommand{\rN}{\ensuremath{\bbbold N}\xspace}
      
      \providecommand{\rP}{\ensuremath{\bbbold P}\xspace}
      
      \providecommand{\rR}{\ensuremath{\bbbold R}\xspace}
      
      \providecommand{\rT}{\ensuremath{\bbbold T}\xspace}

      \providecommand{\rZ}{\ensuremath{\bbbold Z}\xspace}

      \providecommand{\Ae}[1][]{\ensuremath{\ifx|#1|{\ }\else{\:#1\text{-}}\fi\text{almost everywhere\xspace}}}

      \providecommand{\Aa}[1][]{\ensuremath{\text{ for }\ifx|#1|{}\else{\:#1\text{-}}\fi\text{almost all }}}
      \providecommand{\as}[1][]{\ensuremath{\ifx|#1|{\ }\else{#1\text{-}}\fi\text{almost surely}}\xspace}
      \providecommand{\aposteriori}{aposteriori\xspace}
      \providecommand{\Aposteriori}{Aposteriori\xspace}
      
      \providecommand{\apriori}{apriori\xspace}

       \providecommand{\naturals}{\ensuremath{\rN}}
       \providecommand{\Nat}{\ensuremath{\rZ^+}}
       \providecommand{\NO}[1][]{\ensuremath{\naturals_0\ifx|#1|{}\else^{#1}\fi}}

       \providecommand{\reals}{\rR}

       \providecommand{\fieldmats}[3][F]{\csname#1\endcsname{#2\times#3}}
       
       \providecommand{\fieldtens}[3][F]{\csname#1\endcsname{{#2}_1\times\dotsb\times{#2}_{#3}}}

       \providecommand{\RO}[1][]{{\reals_{0+}\ifx|#1|{}\else^{#1}\fi}}
       \providecommand{\RP}[1][]{{\reals_+\ifx|#1|\else^{#1}\fi}}

       \providecommand{\ring}[1][A]{\csname r#1\endcsname}
       \providecommand{\field}[1][K]{\csname r#1\endcsname}

       \providecommand{\torus}[1]{\rT\ifthenelse{\equal{#1}1}{}{^#1}}
      
       \providecommand{\one}{\ensuremath{\bbbold 1}}
       \providecommand{\zerofun}{\ensuremath{\bbbold 0}}
       \providecommand{\ones}[1][]{\one\ifx|#1|\else_{#1}\fi}
       \providecommand{\zeros}[1][]{\zerofun\ifx|#1|\else_{#1}\fi}

       \providecommand{\diracdelta}[1][]{\ensuremath{{\mathrm{\delta}}\ifx|#1|{}\else_{#1}\fi}}

       \providecommand{\pic}{\ensuremath{\mathrm\pi}}%
       \providecommand{\pifracl}[2][]{\fracl{\ifx|#1|\else#1\fi\pic}{#2}}
       \providecommand{\pifrac}[2][]{\frac{\ifx|#1|\else#1\fi\pic}{#2}}

       \providecommand{\take}{\smallsetminus}
       \providecommand{\takesetof}[1]{\take\setof{#1}}
       \providecommand{\takeset}\takesetof
       \providecommand{\takeel}\oldneg%
       \providecommand{\Tolto}{\take}
       \providecommand{\closure}[2][]{\ifx|#1|\overline{#2}\else\operatorname{clos}_{#1}{#2}\fi}
       \providecommand{\inner}{\cdot}
       \providecommand{\vecprod}{\times}
       \providecommand{\outerp}{\wedge}

       \providecommand{\Ga}{\ensuremath{\varGamma}\xspace}

       \providecommand{\W}{\ensuremath{\varOmega}\xspace}

       \providecommand{\ep}{\ensuremath{\varepsilon}\xspace}
       \providecommand{\epsi}{\ensuremath{\epsilon}\xspace}
       \renewcommand{\epsi}{\ensuremath{\epsilon}\xspace}

       \providecommand{\w}{\ensuremath{\omega}\xspace}

       \providecommand{\qp}[2][]{\ensuremath{\ifx|#1|\left(\else\csname#1\endcsname(\fi{#2}\ifx|#1|\right)\else\csname#1\endcsname)\fi}}

       \providecommand{\qpreg}[1]{\ensuremath{(#1)}}
       \providecommand{\qpbig}[1]{\qp[big]{#1}}%
       \providecommand{\qpBig}[1]{\ensuremath{\Big(#1\Big)}}
       \providecommand{\qpbigg}[1]{\ensuremath{\bigg(\!#1\!\bigg)}}
       \providecommand{\qpBigg}[1]{\ensuremath{\Bigg(\!#1\!\Bigg)}}
       \providecommand{\qb}[2][]{\ifx|#1|\left[\else\csname#1\endcsname[\fi{#2}\ifx|#1|\right]\else\csname#1\endcsname]\fi}
       \providecommand{\qc}[2][]{\ensuremath{\ifx|#1|\left\{\else\csname#1\endcsname\{\fi{#2}\ifx|#1|\right\}\else\csname#1\endcsname\}\fi}}

       \providecommand{\qa}[2][]{\ifx|#1|\left\langle\else\csname#1\endcsname\langle\fi{#2}\ifx|#1|\right\rangle\else\csname#1\endcsname\langle\fi}%
       \providecommand{\qareg}[1]{\ensuremath{\langle#1\rangle}}
       \providecommand{\qabig}[1]{\ensuremath{\big\langle#1\big\rangle}}
       \providecommand{\qaBig}[1]{\ensuremath{\Big\langle#1\Big\rangle}}
       \providecommand{\qabigg}[1]{\ensuremath{\bigg\langle#1\bigg\rangle}}
       \providecommand{\qaBigg}[1]{\ensuremath{\Bigg\langle#1\Bigg\rangle}}
       \providecommand{\opinter}[2]{\ensuremath{\left(#1,#2\right)}\xspace}

       \providecommand{\opclinter}[2]{\ensuremath{\left(#1,#2\right]}\xspace}
       \providecommand{\clopinter}[2]{\ensuremath{\left[#1,#2\right)}\xspace}   
       \providecommand{\opintertopinfty}[1]{\opinter{#1}\infty}
       \providecommand\optyinter\opintertopinfty
       \providecommand{\opinterbotinfty}[1]{\opinter{-\infty}{#1}}
       \providecommand\tyopinter\opinterbotinfty
       \providecommand{\clintertopinfty}[1]{\clopinter{#1}\infty}
       \providecommand{\cltyinter}\clintertopinfty
       
       \providecommand{\clinterbotinfty}[1]{\opclinter{-\infty}{#1}}
       \providecommand{\tyclinter}\clinterbotinfty

       \providecommand{\compowqp}[2]{\ensuremath{\qp{\!#2\!\!}^{\kern -.4em #1}\!}}
       
       \providecommand{\powqpreg}[2]{\ensuremath{%
           \qpreg{#2}^{\kern 0em\lower .1ex\hbox{\scriptsize $#1$}}\kern-.3em}}
       \providecommand{\powqpbig}[2]{\ensuremath{%
           \qpbig{#2}^{\kern -.2em\lower .3ex\hbox{\scriptsize $#1$}}\kern-.3em}}
       \providecommand{\powqpBig}[2]{\ensuremath{%
           \qpBig{#2}^{\kern -.2em\lower .3ex\hbox{\scriptsize $#1$}}\kern-.3em}}
       \providecommand{\powqpbigg}[2]{\ensuremath{%
           \qpbigg{#2}^{\kern -.2em\lower .3ex\hbox{\scriptsize $#1$}}\kern-.3em}}
       \providecommand{\powqpBigg}[2]{\ensuremath{%
           \qpBigg{#2}^{\kern -.2em\lower .3ex\hbox{\scriptsize $#1$}}}}
       \providecommand{\powp}[3][]{{#3}\ifx|#1|^{#2}\else{#1}^{#2}\fi}%
       \providecommand{\pow}[2][]{\ifx|#1|\operatorname{pow}^{#2}\else\powp{#2}{#1}\fi}%

       \providecommand{\norm}[2][]{\ifx|#1|\left|\else\csname#1\endcsname|\fi#2\ifx|#1|\right|\else\csname#1\endcsname|\fi}
       \providecommand{\normon}[3][]{\norm[#1]{#2}_{#3}}

       \providecommand{\abs}[2][]{\ensuremath{\ifx|#1|{\left|#2\right|}\else{\csname#1\endcsname|{#2}\csname#1\endcsname|}\fi}}

       \providecommand{\Norm}[2][]{\ifx|#1|\left\|\else\csname#1\endcsname\|\fi{#2}\ifx|#1|\right\|\else\csname#1\endcsname\|\fi}

       \providecommand{\Normon}[3][]{\Norm[#1]{#2}_{#3}}

       \providecommand{\normonsob}[5][]{\normon[#1]{#2}{\sob{#3}{#4}\if|#5|{}\else(#5)\fi}}
       \providecommand{\Normonsob}[5][]{\Normon[#1]{#2}{\sob{#3}{#4}\if|#5|{}\else(#5)\fi}}
       \providecommand{\normonsobh}[4][]{\normon[#1]{#2}{\sobh{#3}\if|#4|{}\else(#4)\fi}}
       \providecommand{\normonsobhz}[4][]{\normon[#1]{#2}{\sobhz{#3}\if|#4|{}\else(#4)\fi}}
       \providecommand{\Normonsobh}[4][]{\Normon[#1]{#2}{\sobh{#3}\if|#4|{}\else(#4)\fi}}
       \providecommand{\Normonsobhz}[4][]{\Normon[#1]{#2}{\sobhz{#3}\if|#4|{}\else(#4)\fi}}
       \providecommand{\Normonleb}[4][]{\Normon[#1]{#2}{\leb{#3}\if|#4|\else(#4)\fi}}

       \providecommand{\ltwop}[3][]{\ensuremath{\qa{#2,#3}\ifx|#1|\else_{#1}\fi}}
       \providecommand{\ltwopreg}[2]{\ensuremath{\qareg{#1,#2}\ifx|#1|\else_{#1}\fi}}
       \providecommand{\ltwopbig}[2]{\ensuremath{\qabig{#1,#2}\ifx|#1|\else_{#1}\fi}}
       \providecommand{\ltwopBig}[2]{\ensuremath{\qaBig{#1,#2}\ifx|#1|\else_{#1}\fi}}
       \providecommand{\ltwopbigg}[2]{\ensuremath{\qabigg{#1,#2}\ifx|#1|\else_{#1}\fi}}
       \providecommand{\ltwopBigg}[2]{\ensuremath{\qaBigg{#1,#2}\ifx|#1|\else_{#1}\fi}}

       \providecommand{\duality}[3][]{\ensuremath{#1\langle #2\,#1\vert\,#3#1\rangle}}
       
       \providecommand{\average}[2][]{{\qa{#2}\ifx|#1|\else_{#1}\fi}}

       \providecommand{\ensemble}[2]{\ensuremath{\left\{ #1:\;#2 \right\}}}
       \providecommand{\setof}[1]{{\qc{#1}}}

       \providecommand{\conditionalto}[1]{{\left|{#1}\right.}}

      \providecommand{\measure}[1]{\ensuremath{\mathcalbf{\MakeUppercase{#1}}}}
      
      \providecommand{\probmeasure}[2][]{{\measure{#2}}\ifx|#1|\else_{#1}\fi}
      \providecommand{\Prob}{}
      \renewcommand{\Prob}[1][]{\probmeasure[{#1}]{p}}

      \providecommand{\randvars}[1][\Prob]{\operatorname{RV}\ifx|#1|{}\else{(#1)}\fi}
      \providecommand{\discrandvars}[1][\Prob]{\operatorname{DRV}\ifx|#1|{}\else{({#1)}\fi}} 
      \providecommand{\contrandvars}[1][\Prob]{\ensuremath{\operatorname{CDRV}\ifx|#1|{}\else(#1)\fi}} 
       \def\env@matrix{\hskip -\arraycolsep
        \let\@ifnextchar\new@ifnextchar
        \array{*\c@MaxMatrixCols c}}
       \renewcommand*\env@matrix[1][c]{\hskip -\arraycolsep
         \let\@ifnextchar\new@ifnextchar
         \array{*\c@MaxMatrixCols #1}}
       \providecommand{\irow}[2]{#1_{#2}}%
       \providecommand{\icol}[2]{#1^{#2}}%

       \providecommand{\ijrowcol}[3]{\icol{\irow{#1}{#2}}{#3}}
       \providecommand{\entry}[1]{\qb{#1}}
       \providecommand{\vecentry}[2]{\irow{#1}{#2}}
       
       \providecommand{\colvecentry}\vecentry
       \providecommand{\covecentry}[2]{\icol{#1}{#2}}
       \providecommand\rowvecentry\covecentry

       \providecommand{\rowof}[1]{\qb{#1}}
       \providecommand{\disvecof}[1]{\begin{bmatrix}[r]#1\end{bmatrix}}
       \providecommand{\vecof}[1]{\mathchoice{\disvecof{#1}}{\qp{#1}}{\qp{#1}}{\qp{#1}}}

       \providecommand{\getentryi}[2]{\irow{\entry{#1}}{#2}}
       \providecommand{\getcolentry}\getentryi
       
       \providecommand{\getvecentry}[2]{\getentryi{\vec #1}{#2}}

       \providecommand{\discolvecitwo}[1]{\discolvectwo{\vecentry{#1}1}{\vecentry{#1}2}}
       
       \providecommand{\discolvecintwo}\discolvecitwo%

       \providecommand{\dismatof}[2][r]{\begin{bmatrix}[#1]#2\end{bmatrix}}

       \providecommand{\matentry}[3]{\ijrowcol{#1}{#2}{#3}}

       \providecommand{\block}[5]{\ijrowcol{#1}{\ifx#2#3{\rowof{#2}}\else\rowof{{#2}\dotsc{#3}}\fi}{\ifx#4#5{\rowof{#4}}\else\rowof{{#4}\dotsc{#5}}\fi}}
       \providecommand{\colblock}[3]{\getvecentry{#1}{\ifx#2#3{#2}\else\fromto{#2}{#3}\fi}}

       \providecommand{\dismatskeldots}[4]{
         \dismatof[c]{
           #1&\dotsc&#3
           \\
           \vdots & \ddots &\vdots
           \\
           #2&\dotsc&#4
         }
       }
       \providecommand{\dismatcommfromtofromto}[5]{
         \dismatskeldots{#1#2#4}{#1#3#4}{#1#2#5}{#1#3#5}
       }
       \providecommand{\dismatcustfromtofromto}[6][matentry]{
         \dismatcommfromtofromto{\csname#1\endcsname{#2}}#3#4#5#6
       }
       \providecommand{\dismatcustfromtofromto}[6][matentry]{
         \dismatskeldots{%
           \csname#1\endcsname{#2}{#3}{#4}%
         }{%
           \csname#1\endcsname{#2}{#3}{#6}%
         }{%
           \csname#1\endcsname{#2}{#5}{#4}%
         }{%
           \csname#1\endcsname{#2}{#5}{#6}%
         }%
       }%
       \providecommand{\dismatcustfromtofromto}[6][matentry]{
         \dismatof{
           \csname#1\endcsname{#2}{#3}{#4}&\dotsc&\csname#1\endcsname{#2}{#3}{#6}
           \\
           \vdots & \ddots &\vdots
           \\
           \csname#1\endcsname{#2}{#5}{#4}&\dotsc&\csname#1\endcsname{#2}{#5}{#6}
         }
       }

       \providecommand{\dissysaxbdotsnm}[5]{\begin{matrix}[r]%
           \matentry{#1}11\vecentry{#2}1&+\dotsb&+\matentry{#1}1{#5}\vecentry{#2}{#5}
           &
           =
           \ifx|#3|0\else{\vecentry {#3}1}\fi
           \\
           \dotsb
           \\
           \matentry{#1}{#4}1\vecentry{#2}1&+\dotsb&+\matentry{#1}{#4}{#5}\vecentry{#2}{#5}
           &
           =
           \ifx|#3|0\else{\vecentry {#3}{#4}}\fi
       \end{matrix}}
       \providecommand{\seqof}[1]{\qp{#1}}%
       \providecommand{\seqs}[2]{\seqof{#1}_{#2}}
       \providecommand{\sets}[2]{\setof{#1}_{#2}}%
       \providecommand{\seqi}[3][]{\seqs{#2_{#3}}{\ifx|#1|{#3}\else{{#3}\in{#1}}\fi}}%
       \providecommand{\sequ}[3][]{\seqs{#2^{#3}}{\ifx|#1|{#3}\else{{#3}\in{#1}}\fi}}%
       \providecommand{\subseqi}[4][]{\seqs{#2_{{#3}_{#4}}}{\ifx|#1|{#4}\else{{#4}\in{#1}}\fi}}%

       \providecommand{\seti}[3][]{\sets{#2_{#3}}{\ifx|#1|_{#3}\else_{{#3}\in{#1}}\fi}}%
       \providecommand{\setu}[3][]{\sets{#2^{#3}}{\ifx|#1|{#3}\else{{#3}\in{#1}}\fi}}%
       \let\liminf\relax
       \DeclareMathOperator*{\liminf}{liminf}
       \let\limsup\relax
       \DeclareMathOperator*{\limsup}{limsup}
       \providecommand{\limofat}[3][]{\ensuremath{\lim_{\ifx|#1|{}\else{#1\ni}\fi#3}{#2}}}
       \providecommand{\limsupofat}[3][]{\ensuremath{\limsup_{\ifx|#1|{}\else{#1\ni}\fi#3}{#2}}}
       \providecommand{\liminfofat}[3][]{\ensuremath{\liminf_{\ifx|#1|{}\else{#1\ni}\fi#3}{#2}}}

       \providecommand{\stringdotsfrom}[3][]{\ensuremath{#2\ifx|#1|\else#1\fi\,#3\ifx|#1|\else#1\fi\,\dotsc}}
       \providecommand{\listdotsfrom}[3][]{\ensuremath{#2\ifx|#1|\else#1\fi,#3\ifx|#1|\else#1\fi,\dotsc}}
       \providecommand{\stringdotsfromto}[3][]{\ensuremath{#2\ifx|#1|\else#1\fi\,\dotsc\,#3\ifx|#1|\else#1\fi}}
       \providecommand{\listdotsfromto}[3][]{\ensuremath{#2\ifx|#1|\else#1\fi,\dotsc,#3\ifx|#1|\else#1\fi}}
       \providecommand{\listifromto}[5][]{\ensuremath{{#2}_{#3}\ifx|#1|\else#1\fi},\text{ for }\ensuremath{\rangefromto{#3}{#4}{#5}}\xspace}
       \providecommand{\listufromto}[5][]{\ensuremath{{#2}^{#3}\ifx|#1|\else#1\fi},\text{ for }\ensuremath{\rangefromto{#3}{#4}{#5}}\xspace}
       \providecommand{\listitwo}[2][]{\ensuremath{#2_1\ifx|#1|\else#1\fi,#2_2\ifx|#1|\else#1\fi}}
       \providecommand{\listutwo}[2][]{\ensuremath{#2^1\ifx|#1|\else#1\fi,#2^2\ifx|#1|\else#1\fi}}
       \providecommand{\listithree}[2][]{\ensuremath{#2_1\ifx|#1|\else#1\fi,#2_2\ifx|#1|\else#1\fi,#2_3\ifx|#1|\else#1\fi}}
       \providecommand{\listithreez}[2][]{\ensuremath{#2_0\ifx|#1|\else#1\fi,#2_1\ifx|#1|\else#1\fi,#2_2\ifx|#1|\else#1\fi}}
       \providecommand{\listifourz}[2][]{\ensuremath{#2_0\ifx|#1|\else#1\fi,#2_1\ifx|#1|\else#1\fi,#2_2\ifx|#1|\else#1\fi,#2_3\ifx|#1|\else#1\fi}}
       \providecommand{\listuthree}[2][]{\ensuremath{#2^1\ifx|#1|\else#1\fi,#2^2\ifx|#1|\else#1\fi,#2^3\ifx|#1|\else#1\fi}}

       \providecommand{\jump}[2][]{\ensuremath{\left\llbracket #2\right\rrbracket\ifx|#1|{}\else_{#1}\fi}}

       \providecommand{\fromto}[2]{\ensuremath{\setof{#1\dotsc#2}}}%
       \providecommand{\fromstepto}[3]{\ensuremath{\left[#1:#2:#3\right]}}

       \providecommand{\integerbetween}[2]{\ensuremath{={#1},\dotsc,{#2}}}

       \providecommand{\rangefromto}[3]{\ensuremath{#1\integerbetween{#2}{#3}}}

       \providecommand{\d}{}
       \renewcommand{\d}[1][]{\ensuremath{\operatorname{d}\!\ifx|#1|\else{_{#1}}\fi}}
       
       \providecommand{\ds}[1][]{\d{\measure S}}
       \providecommand{\D}[1][]{\ensuremath{\operatorname{D}\!\ifx|#1|\else{_{#1}}\fi}}

      \providecommand{\registered}%
      {\ensuremath{^\text{\textregistered}}}

      \providecommand{\constant}[1]{\ensuremath{C_{#1}}}
      \providecommand{\constext}[2][]{\constant{\textup{#2}{\ifx|#1|{}\else{,\ensuremath{#1}}\fi}}}            %
      \providecommand{\constref}[2][]{\ensuremath{\constant{\textup{\ref{#2}{\ifx|#1|{}\else{,\ensuremath{#1}}\fi}}}}}
      \providecommand{\constdef}[2][]{\label{#2}\ensuremath{\constant{\textup{\ref{#2}{\ifx|#1|{}\else{,\ensuremath{#1}}\fi}}}}}

      \providecommand{\funkref}[3][]{\ensuremath{{#3}_{\textup{\ref{#2}{\ifx|#1|{}\else{,\ensuremath{#1}}\fi}}}}}

      \providecommand{\diam}{\operatorname{diam}}
      
      \providecommand{\curl}{\operatorname{curl}}
      \renewcommand{\curl}[1][]{\nabla\ifx|#1|{}\else\kern-2pt_{#1}\fi\kern-2pt\vecprod}
      \renewcommand{\div}[1][]{\nabla\ifx|#1|{}\else\kern-2pt_{#1}\fi\kern-1pt\inner}
      \providecommand{\divof}[2][]{\div[#1]\ifx|#2|{}\else\qb{#2}\fi}%
      \providecommand{\divideabyb}[2]{\operatorname{div}(a,b)}
      \providecommand{\grad}{}
      \renewcommand{\grad}[1][]{\nabla\ifx|#1|\else_{#1}\fi}

      \providecommand{\rot}[1][]{\nabla\ifx|#1|\else_{#1}\fi\outerp}
      
      \providecommand{\rowdiv}[1][]{\D\ifx|#1|{}\else\kern-1pt_{#1}\kern-2pt\fi\cdot}
      \providecommand{\rowdivof}[2][]{\rowdiv[#1]\ifx|#2|{}\else\qb{#2}\fi}

      \providecommand{\interior}{\operatorname{int}}
      \providecommand{\inv}[1][]{\operatorname{inv}\ifx|#1|\else^{#1}\fi}
      \providecommand{\ivt}[1]{\operatorname{ivt}\ifx|#1|\else^{#1}\fi}
      \providecommand\tensorinvariant\ivt
      
      \providecommand{\inverse}[2][]{\powp[#1]{-1}{#2}}

      \providecommand{\mod}{}
      \renewcommand{\mod}[1][]{\operatorname{mod}\ifx|#1|\else\kern-1pt_{#1}\fi}
      
      \let\oldfrac\frac
      \renewcommand{\frac}[3][]{\ifx|#1|\oldfrac{#2}{#3}\else\begin{array}{#1}{#2}\\\hline{#3}\end{array}\fi}
      \providecommand{\fracl}[3][]{\ifx|#1|\nicefrac{#2}{#3}\else{#2}#1/{#3}\fi}

      \providecommand{\qpfracl}[3][]{\qp{\ifx|#1|\fracl{#2}{#3}\else{#2}#1/{#3}\fi}}
      \providecommand{\qpfrac}[3][]{\qp{\ifx|#1|\frac{#2}{#3}\else{#2}#1/{#3}\fi}}
      
      \providecommand{\absfracl}[3][]{\abs{\ifx|#1|\fracl{#2}{#3}\else{#2}#1/{#3}\fi}}
      \providecommand{\absfrac}[3][]{\abs{\ifx|#1|\frac{#2}{#3}\else{#2}#1/{#3}\fi}}
      \providecommand{\fraclff}[3][]{\ifx|#1|{#2}/{#3}\else{#1}\fracl{#2}{#3}\fi}

      \providecommand{\eye}[1][]{\vec{\mathrm I}\ifx|#1|{}\else_{#1}\fi}%
      \providecommand{\numeye}[1][]{\boldsymbol{\mathsf{I}}\ifx|#1|{}\else_{#1}\fi}%
      \providecommand{\Eye}[1]{
        \begin{bmatrix}
        \ifthenelse{#1>1}{
          \ifthenelse{#1>2}{
            \ifthenelse{#1>3}{
              \ifthenelse{#1>4}{
                1&\zeroentry&\dotso&\zeroentry
                \\
                \zeroentry&1&\dotso&\zeroentry
                \\
                \vdots&\vdots&\ddots&\vdots
                \\
                \zeroentry&\zeroentry&\dotso&1
              }{        
                1&\zeroentry&\zeroentry&\zeroentry
                \\
                \zeroentry&1&\zeroentry&\zeroentry
                \\
                \zeroentry&\zeroentry&1&\zeroentry
                \\
                \zeroentry&\zeroentry&\zeroentry&1
              }
            }{
              1&\zeroentry&\zeroentry
              \\
              \zeroentry&1&\zeroentry
              \\
              \zeroentry&\zeroentry&1
            }
          }{
            1&\zeroentry
            \\
            \zeroentry&1
          }
        }{
          1
        }
        \end{bmatrix}
      }
      \providecommand{\Id}{\operatorname{Id}}                   %

      \providecommand{\lebmeas}[1][]{\measure L^{#1}}     %
      \providecommand{\lebmeasof}[2][]{\ifx|#1|\left|#2\right|\else\lebmeas[#1]\qp{#2}\fi}         %
      \providecommand{\meshsize}[1][]{h\ifx|#1|\else_{#1}\fi}
      
      \providecommand{\maxi}[2]{#1\vee#2}                       %

      \providecommand{\mini}[2]{#1\wedge#2}                     %

      \let\oldneg\neg
      \renewcommand{\neg}[1]{\left[#1\right]_-}
      \providecommand{\dash}[1][']{\ifthenelse{\equal{#1}{'}\OR\equal{#1}{''}}{#1}{^{(#1)}}}
      
      \providecommand{\pdfrac}[2][]{\ensuremath{\frac{\partial\ifx|#1|\phantom{#2}\else{#1}\fi}{\partial{#2}}}} %
      \providecommand{\pdfracpow}[3][]{\ensuremath{\frac{\partial^{#3}\ifx|#1|\phantom{#2}\else{#1}\fi}{\partial{#2}^{#3}}}} %

      \providecommand{\pd}[2][]{\ensuremath{\partial_{#2}}{\ifx|#1|{}\else{\qb{#1}}\fi}} %

      \providecommand{\dd}[2][]{\ensuremath{\ifx|#1|\frac{\d}{\d{#2}}\else\frac[l]{\d{#1}}{\d{#2}}\fi}}    %

      \renewcommand{\Im}{\operatorname{im}}                 %
      \renewcommand{\Re}{\operatorname{re}}                 %
      \providecommand{\imaginpart}[1][]{\Im{\ifx|#1|{}\else\qp{#1}\fi}} %
      \providecommand{\realpart}[1][]{\Re{\ifx|#1|{}\else\qp{#1}\fi}} %
      \providecommand\determinant\det

      \providecommand{\transpose}{\intercal}%

      \providecommand{\transposed}{{}^\transpose}

      \providecommand{\orthogonalto}[1][]{\ensuremath{\perp\ifx|#1|{}\else{\!_{#1}\,}\fi}}
      
      \providecommand{\rowof}[1]{\ensuremath{\mathchoice{\vecof{#1}}{\qb{#1}}{\qb{#1}}{\qb{#1}}}}

      \providecommand{\rowvectwo}[2]{\ensuremath{\vecof{#1,#2}}}

      \providecommand{\colvectwo}[2]{\ensuremath{%
          \mathchoice%
              {\discolvectwo{#1}{#2}}%
              {\rowvectwo{#1}{#2}\transposed}%
              {\rowvectwo{#1,#2}\transposed}%
              {\rowvectwo{#1,#2}\transposed}%
        }
      }
      \providecommand{\coltwovec}\colvectwo

      \providecommand{\discolvec}[2][r]{\ensuremath{\begin{bmatrix}[#1]#2\end{bmatrix}}}

      \providecommand{\discolvectwo}[3][r]{\ensuremath{\discolvec[#1]{#2\\#3}}}

      \providecommand{\discolvecitwo}[1]{\discolvectwo{\vecentry{#1}1}{\vecentry{#1}2}}

      \provideboolean{showzeroentries}%
      \setboolean{showzeroentries}{true}%
      \providecommand{\zeroentry}{\ifthenelse{\boolean{showzeroentries}}{{0}}{\phantom0}}
      \providecommand{\zeroentrywarning}{\ifthenelse{\boolean{showzeroentries}}{}{%
          \ensuremath{\text{($0$ entries omitted)\xspace}}}}

      \providecommand{\smint}{\ensuremath{{\text{\textbf{/}}}\kern-.75em\smallint}}
      \renewcommand{\smint}[1][]{\lower12.3pt\hbox{\begin{tikzpicture}\draw[line width=.75pt] (-3pt,-0.5)--(1pt,-0.5) node[pos=0.6]{$\int$};\path (3pt,-24pt)node {\scriptsize $#1$};\end{tikzpicture}}}

      \providecommand{\lap}{\ensuremath{\mathrm\Delta}}
      \providecommand{\lapin}[1][]{\lap\ifx|#1|\else_{#1}\fi}
      \providecommand{\normalsymbol}{\operatorname{\mathbf{n}}}
      \renewcommand{\normalsymbol}{\vec{\operatorname{n}}}
      \providecommand{\normal}[1][]{\normalsymbol\ifx|#1|\else_{#1}\fi}%
      \providecommand{\norm@l}[1][]{\normalsymbol\ifx|#1|\else_{#1}\fi}%
      \providecommand{\normalto}[2][]{\ensuremath{\norm@l[#2]\ifx|#1|\else\qp{#1}\fi}}
      \providecommand{\normalder}[1][]{\ensuremath{\norm@l\ifx|#1|\else\qp{#1}\fi{\inner\grad}}}

      \providecommand{\tangentialsymbol}{\operatorname{\textbf{t}}}
      \providecommand{\tangentialto}[2][]{\tangentialsymbol\ifx|#1|\else^{#1}\fi\ifx|#2|\else_{#2}\fi}

      \providecommand{\union}[1]{\ensuremath{\bigcup\nolimits_{#1}}}

      \providecommand{\unions}[3][]{\union{#2\in{#3}\ifx|#1|\else:#1\fi}}

      \ifthenelse{\boolean{usesvjour}}{}{
        \let\vec\undefined
        \providecommand{\vec}[1]{\ensuremath{\boldsymbol{#1}}}
        \renewcommand{\vec}[1]{\ensuremath{\boldsymbol{#1}}}
      }

      \providecommand{\hatmat}[1]{\hat{\mat{#1}}}

      \providecommand{\geomat}[1]{\vec{\UCmath{#1}}}

      \providecommand{\mat}[1]{\geomat{#1}} %
      \providecommand{\Prob}[1][]{\ensuremath{\operatorname{Prob}\ifx|#1|{}\else_{#1}\fi}}

      \providecommand{\pdf}[2][]{\ensuremath{\operatorname{pdf}_{#2\ifx|#1|{}\else{\conditionalto{#1}}\fi}}\xspace}

      \providecommand{\expectation}{\ensuremath{\operatorname{E}}}
      \providecommand{\EX}[1][]{\ensuremath{\expectation\ifx|#1|{}\else_{#1}\fi}}

      \providecommand{\gausskernel}[3][x]{%
        \ensuremath{
          \exp\frac{-\if#20{#1}\else(#1-\mu)\fi^2}{%
            2\if#31{}\else\powp2{#3}\fi}%
        }%
      }
      \providecommand{\gaussdistribution}[3][x]{%
        \ensuremath{\frac1{\sqrt{2\pic}\if#31{}\else#3\fi}%
          \gausskernel[#1]{#2}{#3}
        }%
      }%

      \providecommand{\PD}[1]{\operatorname{PD}\qpreg{#1}}
      \providecommand{\pdspace}[1]{\PD{\linspace v}}
      \providecommand{\pdmats}[2][F]{\PD{\csname#1\endcsname{#2}}}
      
      \providecommand{\SPD}{\operatorname{SPD}}
      \providecommand{\spdmats}[2][F]{\SPD(\csname#1\endcsname{#2})}

       \providecommand{\Continuous}{\ensuremath{\operatorname C}\xspace}%
       \providecommand{\Hspace}{\ensuremath{\operatorname H}\xspace}
       \providecommand{\Lebesgue}{\ensuremath{\operatorname L}\xspace}
       \providecommand{\Besovspace}{\ensuremath{\operatorname B}\xspace}

       \providecommand{\Weaklyder}{\ensuremath{\operatorname W}\xspace}
       
       \providecommand{\dual}[1]{\ensuremath{{#1}'}}
       
       \providecommand{\dualspace}[2][]{\dual{\linspace{#2}\ifx|#1|\else{_{#1}}\fi}}
       \providecommand{\bidual}[1]{\ensuremath{{#1}''}}
       \providecommand{\bidualspace}[2][]{\bidual{\linspace{#2}\ifx|#1|\else{_{#1}}\fi}}

       \providecommand{\cont}[1]{\ensuremath{\Continuous^{#1}}}

       \providecommand{\diffable}[2][]{\ensuremath{\cD\ifx|#1|\else^{#1}\fi(#2)}}
       
       \providecommand{\BV}[1]{\ensuremath{\operatorname{BV}}}

       \providecommand{\leb}[1]{\ensuremath{\Lebesgue_{#1}}}

       \providecommand{\lebloc}[1]{\ensuremath{{{\Lebesgue}^{\kern-.20em\lower .1ex\hbox{\tiny\textrm{\textup{loc}}}}_{#1}}}}
       \providecommand{\lebnorm}[3][]{\ensuremath{\Norm{#2}_{\leb{#3}\ifx|#1|{}\else(#1)\fi}}}
       \providecommand{\bes}[3][]{\ensuremath{\Besovspace^{#2}_{#3\ifx|#1|\else,#1\fi}}}

       \providecommand{\sob}[2]{\ensuremath{{\smash\Weaklyder}^{#1}_{#2}}}

       \providecommand{\sobh}[1]{\ensuremath{\Hspace^{#1}}}
       \providecommand{\vecsobh}[1]{\ensuremath{\vec\Hspace^{#1}}}
       \ProvideDocumentCommand{\hdiv}{ O{} O{}}{\vecsobh{\operatorname{div}}\ifx+#1+\else_{0|#1}\fi\ifx|#2|\else(#2)\fi}
       \providecommand{\hcurl}[1][]{\vecsobh{\operatorname{curl}}\ifx|#1|\else(#1)\fi}
       
       \providecommand{\sobhz}[2][]{\sobh{#2}_{0\ifx+#1+\else|#1\fi}}

       \providecommand{\Lip}[1][]{\ensuremath{\operatorname{Lip}}\ifx|#1|{}\else{\qp{#1}}\fi}

       \ProvideDocumentCommand{\polyring}{ O{X} O{A} }{\ring[#2][#1]}
       \ProvideDocumentCommand{\polyfield}{ O{X} O{} }{\field[#2][#1]}
       \providecommand{\polyreals}[1][]{\polyfield[][R]\ifx|#1|\else^{#1}\fi}
       \providecommand{\poly}[2][]{\ensuremath{\rP\ifx#1\else_{#1}\fi^{#2}}}

       \providecommand{\Symmatrices}[2][R]{\ensuremath{\operatorname{Sym}{(\csname#1\endcsname{#2})}}}
       
       \providecommand{\SAmatrices}[2][F]{\ensuremath{\operatorname{SA}{(\csname#1\endcsname{#2})}}}
       \providecommand{\mesh}[2][]{\ensuremath{\mathcalbf{\MakeUppercase{#2}}\ifx|#1|\else_{#1}\fi}}

      \providecommand{\crouzeixraviart}[1][1]{\operatorname{CR}\ifx|#1|{}\else{^{#1}}\fi}

      \providecommand{\linspace}[1]{\mathscript{\MakeUppercase{#1}}}

      \providecommand{\Lin}{\operatorname{Lin}}
      \providecommand{\CL}{\operatorname{CL}}
      \providecommand{\linops}[3][]{\ensuremath{\Lin\ifx|#1|\else^{#1}\fi\qp{{#2}\to{#3}}}}

      \providecommand{\clinops}[3][]{\ensuremath{\CL\ifx|#1|\else^{#1}\fi\qp{{#2}\to{#3}}}}
      
      \providecommand{\fepartition}[2][]{\mathscript{\MakeUppercase{#2}}\ifx|#1|{}\else_{#1}\fi}
      \providecommand{\fespace}[2][]{\mathbb{\MakeUppercase{#2}}\ifx|#1|{}\else_{#1}\fi}
      \providecommand{\hatfespace}[2][]{\widehat{\mathbb{\MakeUppercase{#2}}}\ifx|#1|{}\else_{#1}\fi}

      \providecommand{\vespace}[1][]{\fespace v\ifx|#1|\else_{#1}\fi}
      \providecommand{\hatvespace}[1][]{\hatfespace v\ifx|#1|\else_{#1}\fi}
      \providecommand{\fezerospace}[2][]{\ensuremath{\mathring{\fespace{#2}}\ifx|#1|{}\else_{#1}\fi}}

      \providecommand{\fes}[2]{\ensuremath{\fespace{#1}^{#2}}}

      \providecommand{\fez}[2]{\ensuremath{\fezerospace{#1}^{#2}}}

      \providecommand{\fe}[2][]{\ensuremath{\UCmath{#2}\ifx|#1|\else_{#1}\fi}}%

      \providecommand{\vecfe}[2][]{\ensuremath{\vec{\fe{#2}}\ifx|#1|{}\else{_{#1}}\fi}}%
      
      \providecommand{\matfe}[2][]{\ensuremath{\mat{\fe{#2}}\ifx|#1|{}\else{_{#1}}\fi}}%
      
      \providecommand{\hatmatfe}[2][]{\ensuremath{\hatmat{\UCmath{#2}}\ifx|#1|{}\else{_{#1}}\fi}}%
      
      \providecommand{\EOC}{\ensuremath{\operatorname{EOC}}\xspace}

      \providecommand{\emptyset}{\varnothing}
      \renewcommand{\emptyset}{\varnothing}

      \RequirePackage{amscd}

      \providecommand{\isomorphicto}{\leftrightarrows}
      \providecommand\isomorphic\isomorphicto

      \providecommand{\restriction}[2]{\left.#1\right|_{#2}}
      \renewcommand{\restriction}[2]{\left.#1\right|_{#2}}

      \providecommand{\evalat}[3][]{\qb{#2}_{\ifx|#1|{}\else#1=\fi#3}}
      \providecommand{\evaldiff}[4][]{\qb{#2}^{\ifx|#1|{}\else#1=\fi#3}_{\ifx|#1|{}\else#1=\fi#4}}

      \providecommand{\hhat}[1]{\hat{\hat{#1}}} %
      \providecommand{\ccheck}[1]{\check{\check{#1}}} %
      \providecommand{\bs}{\char '134}   %

      \providecommand{\Program}[1]{\nolinkurl{#1}\xspace}
      \providecommand{\Source}[1]{\nolinkurl{#1}\xspace}

      \providecommand{\matlabplot}[2][]{%
        \begin{center}
          \includegraphics[width=0.9375\linewidth,trim=64 200 64 200,clip]{#2}%
          \ifx|#1|\else\\#1\fi
        \end{center}%
      }

      \providecommand{\texcommand}[1]{\texttt{\bs{\nolinkurl{#1}}}\xspace}
      
      \providecommand{\codename}[1]{\nolinkurl{#1}\xspace}
      \providecommand\colorvar[2][a]{\colorbox{#1!6.25}{#2}}%
      \providecommand{\colorvarname}[2][a]{\colorvar[#1]{\Verb{#2}}}
      
      \providecommand{\codevarname}[1]{\colorvarname[a]{#1}}
      
      \providecommand\olco\codevarname

      \ifthenelse{\boolean{useutopia}}{%
        
      }{%
        
      }

      \providecommand{\matlab}{{\small\Program{MATLAB}}\xspace}%
      \providecommand\MATLAB\matlab

      \ProvideDocumentCommand{\codesnip}{ O{.} O{1.0} m}{%
        \newline
        \begin{minipage}{#2\linewidth}
          \lstinputlisting{#1/#3}
        \end{minipage}
      }
      \providecommand{\codeprint}[2][.]{
        \ \newline
        \begin{minipage}{\linewidth}
          \lstinputlisting{#1/#2}
          \framebox{Contents of file %
            \ifthenelse{\isundefined\pickuppath}{%
             \codename{#2}%
            }{%
              \providecommand{\fullpickuppath}{}%
              \renewcommand{\fullpickuppath}{\pickuppath/\ifx|#1|\else#1/\fi#2}%
              \href{\fullpickuppath}{\codename{#2}}%
          }}
        \end{minipage}
      }
      \providecommand{\codenoprint}[2][.]{
              \providecommand{\fullpickuppath}{}%
              \renewcommand{\fullpickuppath}{\pickuppath/\ifx|#1|\else#1/\fi#2}%
              \href{\fullpickuppath}{\codename{#2}}%
      }

      \providecommand{\indexen}[2][]{{\ifthenelse{\boolean{shownotes}}{\color b}{}#2\ifx|#1|\index{#2}\else\index{#1}\fi}}

      \providecommand{\indexma}[2][]{{\ifthenelse{\boolean{shownotes}}{\color b}{}#2\ifx|#1|\index{\(#2\)}\else\index{<#1@\(#2\)}\fi}}
      
      \providecommand{\secsymbol}{\S}
      
      \providecommand{\secref}[1]{\secsymbol\ref{#1}}

      \providecommand{\ListParameters}{}
      \renewcommand{\ListParameters}%
      {
      	 \setlength{\topsep}{0pt}
      	 \setlength{\leftmargin}{0pt}
               \setlength{\itemsep}{0pt}
      	 \setlength{\parsep}{0pt}
      	 \setlength{\parskip}{0pt}
               \setlength{\labelsep}{0pt}
      	 \setlength{\itemindent}{0pt}
      }
      {%
        \begin{list}%
          {}%
          {\ListParameters%
          
      }}%
      {\end{list}}
      \newcounter{tmpcounter}
      \newcounter{LetterListItem}
      \renewcommand{\theLetterListItem}{(\alph{LetterListItem})}
      \newenvironment{LetterList}[1][0]{%
        \begin{list}{%
            \theLetterListItem\ %
          }{%
            \usecounter{LetterListItem}%
            \ListParameters%
            \ifx|#1|{}\else\setcounter{LetterListItem}{#1}\fi%
          }%
      }{%
        \end{list}%
      }
      \newcounter{CapitalListItem}
      \renewcommand{\theCapitalListItem}{\Alph{CapitalListItem}.}

      \newcounter{NumberListItem}
      \renewcommand{\theNumberListItem}{\arabic{NumberListItem}}
      \newenvironment{NumberList}%
      {
      	\begin{list}%
      	{\theNumberListItem.\ }%
      	{\usecounter{NumberListItem}%
      	 \ListParameters
      	}
      }%
      {\end{list}}
      \newcounter{QuestionListItem}
      \renewcommand{\theQuestionListItem}{\textbf{Question \arabic{QuestionListItem}}}
      {
      	\begin{list}%
      	{\theQuestionListItem.\ }%
      	{\usecounter{QuestionListItem}%
      	 \ListParameters
      	}
      }%
      {\end{list}}
      \newcounter{RomanListItem}
      \renewcommand{\theRomanListItem}{(\roman{RomanListItem})}
      {
      	\begin{list}%
      	{\theRomanListItem\ }%
      	{\usecounter{RomanListItem}
      	 \ListParameters
      	}
      }%
      {\end{list}}
      \newcounter{StepsItem}
      {
      	\begin{list}%
      	{Step \theStepsItem.\ }%
      	{\usecounter{StepsItem}%
      	 \ListParameters
      	}
      }%
      {\end{list}}
      \newcounter{CasesListItem}
      \renewcommand{\theCasesListItem}{\Alph{CasesListItem}}
      {
      	\begin{list}%
      	{\emph{Case \theCasesListItem.}\ }%
      	{\usecounter{CasesListItem}%
      	 \ListParameters
      	}
      }%
      {\end{list}}
      \newcounter{QAListItem}
      \renewcommand{\theQAListItem}{Q\arabic{QAListItem}:}
      {
      	\begin{list}%
      	{\theQAListItem}%
      	{\usecounter{QAListItem}
      	 \ListParameters
      	}
      }%
      {\end{list}}

      \ifthenelse{\boolean{isthesis}}{%
        \setcounter{secnumdepth}{1}%
      }{
        \setcounter{secnumdepth}{2} %
      }
      \providecommand{\ListParameters}{}
      \renewcommand{\ListParameters}
      {
      	 \setlength{\topsep}{0em}
      	 \setlength{\leftmargin}{0em}
               \setlength{\itemsep}{0ex}
      	 \setlength{\parsep}{.5ex}
      	 \setlength{\itemindent}{\labelsep}
      	 \addtolength{\itemindent}{\labelwidth}
      }

        \providecommand{\ObsName}{Remark}%
        \providecommand{\RemName}{Remark}%
        \providecommand{\NotName}{Notation}%
        \providecommand{\BFNName}{Big~Fat~Note}%
        \providecommand{\DefName}{Definition}%
        \providecommand{\ExaName}{Example}%
        \providecommand{\TheName}{Theorem}%
        \providecommand{\LemName}{Lemma}%
        \providecommand{\ProName}{Proposition}%
        \providecommand{\CorName}{Corollary}%
        \providecommand{\PbmName}{Problem}%
        \providecommand{\HypName}{Hypothesis}%
        \providecommand{\AlgName}{Algorithm}%
        \providecommand{\ExeName}{Exercise}%
        \providecommand{\SolName}{Solution}%
        \providecommand{\ClaName}{Claim}%
        \providecommand{\EsyName}{Essay}%
        \providecommand{\Proofname}{Proof}%
        \providecommand{\Derivename}{Derivation}%
      
      \ifthenelse{\boolean{isthesis}}{%
        \providecommand{\Thecounter}{The}
      }{%
        \providecommand{\Thecounter}{subsection}
      }
      \newcommand{\oltikzgetxy}[3]{%
        \tikz@scan@one@point\pgfutil@firstofone#1\relax
        \edef#2{\the\pgf@x}%
        \edef#3{\the\pgf@y}%
      }
      \providecommand{\pdfformat}[1]{
         \provideboolean{pdfoutput}
         \setboolean{pdfoutput}{#1}%
        \ifthenelse{\boolean{pdfoutput}}{
          \typeout{using pdf}
\makeatletter
\usepackage{pdfsync}
          \providecommand{\graphext}{pdf}
          \renewcommand{\graphext}{pdf}
          \providecommand{\graphextex}{pdf_t}
          \renewcommand{\graphextex}{pdf_t}
        }{
          \typeout{using eps}
          \RequirePackage[dvips]{graphicx,xcolor}
          \providecommand{\graphext}{eps}
          \renewcommand{\graphext}{eps}
          \providecommand{\graphextex}{eps_t}
          \renewcommand{\graphextex}{eps_t}
        }
        \RequirePackage{epsfig}
        \RequirePackage{tikz}
        \RequirePackage{rotating}
\makeatletter
        \RequirePackage{graphicx}
        \RequirePackage{xcolor}
        \provideboolean{darkcolortheme}
        \definecolor{SussexFlint}{rgb}{.00,.19,.21}
        \definecolor{SussexGrey}{rgb}{.51,.58,.49}
        \definecolor{SussexOrange}{rgb}{.94,.29,.00}
        \definecolor{SussexYellow}{rgb}{1.00,.73,.00}
        \definecolor{SussexRed}{rgb}{.94,.01,.49}
        \definecolor{SussexPurple}{rgb}{.48,.06,.44}
        \definecolor{SussexGreen}{rgb}{.00,.58,.46}
        \definecolor{OmarGreen}{rgb}{.00,.68,.36}
        \definecolor{SussexBlue}{rgb}{.00,.58,.65}
        \definecolor{OmarBlue}{rgb}{.00,.38,.65}
        \colorlet{a}{OmarBlue}%
        \colorlet{b}{SussexOrange}
        \colorlet{c}{SussexGreen}
        \colorlet{d}{SussexPurple}%
        \colorlet{e}{SussexRed}
        \colorlet{f}{SussexYellow}
        \colorlet{g}{white}%
        \colorlet{h}{SussexGrey}%
        \colorlet{i}{black}%
        \colorlet{j}{SussexFlint}
        \colorlet{colora}{a}
        \colorlet{colorb}{b}
        \colorlet{colorc}{c}
        \colorlet{colord}{d}
        \colorlet{colore}{e}
        \colorlet{colorf}{f}
        \colorlet{colorg}{g}
        \colorlet{colorh}{h}
        \colorlet{colori}{i}
        \colorlet{colorj}{j}
        \newcommand{\mausDarkColorTheme}{
          \colorlet{a}{SussexYellow!50!yellow}
          \colorlet{b}{SussexBlue}%
          \colorlet{c}{SussexRed!50!red}
          \colorlet{d}{SussexOrange!50!yellow}
          \colorlet{e}{SussexGreen!50!green}
          \colorlet{f}{SussexPurple!50!magenta}
          \colorlet{g}{black}%
          \colorlet{h}{SussexFlint!50!black}
          \colorlet{i}{white}%
          \colorlet{j}{SussexGrey}
        }
        \ifthenelse{\boolean{darkcolortheme}}{\mausDarkColorTheme}{}
\makeatletter
      }
      \providecommand{\solution}{\textbf{\SolName.}\xspace}

      \newcounter{phantomedinput}
      \newcounter{phantombox}
      \counterwithin*{phantombox}{phantomedinput}%
      \provideboolean{showphantoms}
      \renewcommand{\thephantombox}{\Alph{phantombox}}%
      \providecommand{\phantombox}[1]{\stepcounter{phantombox}%
        \ensuremath{\boxed{%
            {\ifthenelse{\boolean{showphantoms}}{#1}{\phantom{#1}}}%
            {\texttt{\tiny\ \colorbox{i!50}{\color g\thephantombox}}
            }%
          }%
        }%
      }

      \provideboolean{hidesolution}
      \newcommand{\consolution}[2][]{
        \ifthenelse{\boolean{hidesolution}}{#1\setboolean{showphantoms}{false}}{%
          {\setboolean{showphantoms}{true}\color{i!50}\par \small {\solution}\ #2\par\ \\[5pt]}}
      }
      \provideboolean{showmarks}
      \providecommand{\showmarks}[1]{%
        \ifthenelse{%
          \boolean{showmarks}}{%
          \marginpar{%
            \tiny [$#1$ mark\ifthenelse{\equal{#1}1}{\phantom{s}}s]}%
        }{}}%

      \newcommand{\condibreak}{\ifthenelse{\boolean{hidesolution}}{\clearpage}{}}
      
      \newcommand{\solutibreak}{\ifthenelse{\boolean{hidesolution}}{}{\clearpage}}
      \newcommand{\questionly}[1]{\ifthenelse{\boolean{hidesolution}}{#1}{}}
      \newcommand{\solutionly}[1]{\ifthenelse{\boolean{hidesolution}}{}{#1}}

       \providecommand{\qeyword}[1]{\index{#1}\ifthenelse{\boolean{shownotes}}{{\tiny\color e\colorbox{e!6.25}{#1}}}{}}
       \providecommand{\pathwordbase}{file:///Users/ol20/Sokratex}
       \providecommand{\pathword}[2][]{%
         \label{#2}%
         \ifthenelse{\boolean{shownotes}}{%
           \ \\\index{#2@\tiny\codevarname{#2}}{%
             \tiny{\color f\href{\pathwordbase/#2}{\colorvarname[f]{#2}%
           }}}\\
         }{}%
       }
       \providecommand{\targword}[2][]{%
         \label{#2}%
         \ifthenelse{\boolean{shownotes}}{%
           \index{#2@\tiny\codevarname{#2}}{\ensuremath{\tiny\color f\href{\pathwordbase/#2}{\colorvarname[e]{#2}\ifx|#1|\else\colorvarname[e]{[#1]}\fi}}}\\
         }{}%
       }
       \providecommand{\sourceurl}[2][]{%
         \ifthenelse{\boolean{shownotes}}{{\ \\\tiny\colorbox{d!6.25}{\color d\texttt{source: \ifx|#1|\href{#2}{#1}\else\url{#2}\fi}}}}}
       \providecommand{\sourcecite}[2][]{\ifthenelse{\boolean{shownotes}}{{\ \\\tiny\colorbox{d!6.25}{\color d\texttt{source: \citet[#1]{#2}}}}}{%
       }}
       \providecommand{\conword}[2][]{\ifthenelse{\boolean{shownotes}}{#2}{#1}}
       \providecommand{\solword}[2][]{\ifthenelse{\boolean{hidesolution}}{#1}{#2}}
       \providecommand{\solghost}[1]{\ifthenelse{\boolean{showphantoms}}{#1}{\phantom{#1}}}

      \RequirePackage{lineno}
      \ifthenelse{\boolean{showchanges}}{
        \newcommand{\llabel}[1]{\hypertarget{llineno:#1}{\linelabel{#1}}}
        \newcommand{\lref}[1]{\hyperlink{llineno:#1}{\ref*{#1}}}
      }{
        \newcommand\llabel[1]{}
        \newcommand\lref[1]{}
      }
      \provideboolean{includeresponses}
      \setboolean{includeresponses}{false}
      \providecommand{\mailto}[1]{\href{mailto:#1}{\nolinkurl{#1}}}
      \provideboolean{showoldetails}
      \setboolean{showoldetails}{true}
      \providecommand{\oldetails}[2]{\ifthenelse{\boolean{showoldetails}}{#1}{#2}}
      \RequirePackage{hyphenat}
   \ifthenelse{\boolean{nohyperref}}{
     \PackageWarning{omarstyle}{package hyperref not loaded, please load manually if needed}
   }{
     \RequirePackage[obeyspaces]{url}
     \RequirePackage[bookmarks,hypertexnames=false,debug,pdfpagelabels]{hyperref}
     \RequirePackage{bookmark}
   }

   \newtheoremstyle{plain}%
     {}%
     {}%
     {\mdseries\slshape}%
     {\parindent}%
     {\bfseries}%
     {.}%
     {.5em}%
     {}%
   
   \newtheoremstyle{note}%
     {}%
     {}%
     {}%
     {\parindent}%
     {\bfseries}%
     {.}%
     {.5em}%
     {}%
   
   \newtheoremstyle{claim}%
     {}%
     {}%
     {\mdseries\slshape}%
     {}%
     {\bfseries}%
     {}%
     {.5em}%
     {}%
   
   \newtheoremstyle{exercise}%
     {}%
     {}%
     {}%
     {}%
     {\bfseries}%
     {.}%
     {1em}%
     {}%
   
   \newtheoremstyle{break}%
     {}%
     {}%
     {}%
     {}%
     {\bfseries}%
     {.}%
     {\newline}%
     {}%
   
   \swapnumbers{
     \theoremstyle{plain}
     \ifthenelse
         {\boolean{isthesis}}
         {\newtheorem{The}{\TheName}[section]}%
         {
         }%
   {
      \theoremstyle{plain}
   
      \renewcommand{\Thecounter}{subsection}

      \newtheorem*{The*}{\TheName}
      \newtheorem*{Lem*}{\LemName}
      \newtheorem*{Pro*}{\ProName}
      \newtheorem*{Cor*}{\CorName}
      \newtheorem*{Pbm*}{\PbmName}
      \newtheorem*{Hyp*}{\HypName}
      \newtheorem*{Exe*}{\ExeName}
      \newtheorem*{Txx*}{\ExeName} %
      \newtheorem*{Con*}{Conclusion}
      \newtheorem*{Sum*}{Summary}
    }
    {
      \theoremstyle{claim}

    }
    {
      \theoremstyle{note}

      \newtheorem*{Obs*}{\ObsName}

      \newtheorem*{Def*}{\DefName}
      \newtheorem*{Exa*}{\ExaName}
      \newtheorem*{Alg*}{\AlgName}
    }
   
    {
      \theoremstyle{break}
    }
   }

   \ifthenelse{\boolean{usecleveref}}{
     \usepackage{aliascnt}
     \usepackage{cleveref}
     \newtheorem{OLThe}[subsection]{Theorem}%
     \newenvironment{The}[1][]{%
       \begin{OLThe}[#1]%
         \addcontentsline{toc}{subsection}{\thesubsection. \TheName\ (#1)}
     }{\end{OLThe}}
     \newaliascnt{Lem}{subsection}%
     \newtheorem{OLLem}[Lem]{Lemma}
     \newenvironment{Lem}[1][]{%
       \begin{OLLem}[#1]%
         \addcontentsline{toc}{subsection}{\thesubsection. \LemName\ (#1)}
     }{\end{OLLem}}  
     \newaliascnt{Pro}{subsection}
     \newtheorem{OLPro}[Pro]{Proposition}
     \newenvironment{Pro}[1][]{%
       \begin{OLPro}[#1]%
         \addcontentsline{toc}{subsection}{\thesubsection. \ProName\ (#1)}
     }{\end{OLPro}}  
     \newaliascnt{Def}{subsection}
     \newtheorem{OLDef}[Def]{Definition}
     \newenvironment{Def}[1][]{%
        \begin{OLDef}[#1]%
          \addcontentsline{toc}{subsection}{\thesubsection. \DefName\ of #1}
     }{\end{OLDef}}
     \newaliascnt{Obs}{subsection}
     \newtheorem{OLObs}[Obs]{Remark}
     \newenvironment{Obs}[1][]{%
       \begin{OLObs}[#1]%
         \addcontentsline{toc}{subsection}{\thesubsection. \ObsName\ of #1}
     }{\end{OLObs}}
     \aliascntresetthe{Lem}
     \aliascntresetthe{Pro}
     \aliascntresetthe{Def}
     \aliascntresetthe{Obs}
     \crefname{Lem}{lemma}{lemmata}
     \crefname{Pro}{proposition}{propositions}
     \crefname{Def}{definition}{definitions}
     \crefname{Obs}{remark}{remarks}
   }{
     \newenvironment{The}[1][]{%
       \ifx&#1&%
       \subsection{\TheName\xspace}%
       \else%
       \subsection[\MakeUppercase#1 theorem]{\TheName\ (#1)}%
       \fi%
       \slshape}{%
       \upshape}
     \newenvironment{Pro}[1][]{\subsection{\ProName\xspace{\ifx&#1&{}\else{ (#1)}\fi}}\slshape}{\upshape}
     \newenvironment{Lem}[1][]{\subsection{\LemName\xspace{\ifx&#1&{}\else{ (#1)}\fi}}\slshape}{\upshape}

     \newenvironment{Def}[1][]{\subsection{\DefName\xspace{\ifx&#1&{}\else{ of \indexen{#1}}\fi}}}{}
     \newenvironment{Obs}[1][]{\subsection{\ObsName\xspace{\ifx&#1&{}\else{ (#1)}\fi}}}{}

   }
   \providecommand{\qed}{\vrule height 5pt depth 0pt width 3pt}
   \providecommand{\qqed}{{\raggedright{\ \hfill\qed}}}
   
   \newcounter{passo}

   \newenvironment{Proof}[1][]%
   {\par\noindent{\bf \Proofname\ifx|#1|.\ \else\ #1.\ \fi}\setcounter{passo}{0}}%
   {\qqed\par}
   {\par\noindent{\bf \Derivename\ #1}\setcounter{passo}{0}}%
   {\qqed\par}
   \newenvironment{Proof*}[1][{}]%
   {\subsection{\Proofname\ #1}\setcounter{passo}{0}}
   {\qqed\par}

\makeatletter

\providecommand{\estspace}[1]{\ensuremath{\eta_{#1}}}

\providecommand{\estrec}[1]{\ensuremath{\varepsilon_{#1}}}

\providecommand{\inddata}{\ensuremath{\gamma}}
\providecommand{\indf}{\ensuremath{\beta}}
\providecommand{\indspace}{\ensuremath{\eta}}
\providecommand{\indtime}{\ensuremath{\theta}}
\providecommand{\indrec}{\ensuremath{\varepsilon}}
\providecommand{\honew}{{\sobh1(\W)}} %
\providecommand{\honezw}{{\sobhz1(\W)}} %
\providecommand{\ccr}[2]{\ensuremath{{#1}\vee{#2}}}
\providecommand{\fcc}[2]{\ensuremath{{#1}\wedge{#2}}}

\providecommand{\tm}{{t_m}}
\providecommand{\tms}{{t_m^*}}
\providecommand{\tn}{{t_n}}
\providecommand{\tno}{{t_{n-1}}}

\providecommand{\timestep}{\tau}
\providecommand{\tstep}[1]{\ensuremath{\timestep_{#1}}}
\providecommand{\taun}{{\tstep n}}
\providecommand{\taunm}{\inverse\taun}

\providecommand{\hn}{\ensuremath{h_n}}
\providecommand{\hno}{\ensuremath{h_{n-1}}}
\providecommand{\hathn}{\ensuremath{\hat h_n}}

\providecommand{\barf}{\ensuremath{\overline f}}
\providecommand{\projf}[1]{\ensuremath{{\barf}^{#1}}}%

\providecommand{\U}{}
\renewcommand{\U}{\ensuremath{U}}%

\providecommand{\un}{\U^n}
\providecommand{\uno}{\U^{n-1}}

\providecommand{\rhoms}{\rho^m_*}
\providecommand{\ms}{{m^*}}
\providecommand{\mso}{{m^*-1}}

\providecommand{\elop}{{\cA}}
\providecommand{\pwA}{\elop_{\mathrm{el}}}%
\providecommand{\discelop}[1][A]{\MakeUppercase{#1}}
\providecommand{\disca}[1]{\discelop^{#1}}

\providecommand{\an}{\disca n}%
\providecommand{\ano}{\disca{n-1}}%
\providecommand{\recop}[1]{{\cR^{#1}}}%

\providecommand{\lprojop}{P_0}
\providecommand{\lproj}[1]{\lprojop^{#1}}%

\providecommand{\eproj}[1]{P_1^{#1}}%

\providecommand{\abil}[2]{\ensuremath{a\left(#1,#2\right)}}

\providecommand{\bdisc}{\ensuremath{\partial}}

\providecommand{\fdisc}{\ensuremath{\overset{\smash{\shortrightarrow}}\partial}}
\providecommand{\sdisc}{\ensuremath{{\partial^2}}}
\providecommand{\discn}{\ensuremath{\overline\partial}}

\renewcommand{\Aposteriori}{A~posteriori\xspace}
\renewcommand{\aposteriori}{a~posteriori\xspace}

\renewcommand{\apriori}{a~priori\xspace}

\providecommand{\estspace}[1]{\ensuremath{\eta_{#1}}}

\providecommand{\estrec}[1]{\ensuremath{\varepsilon_{#1}}}

\providecommand{\inddata}{\ensuremath{\gamma}}
\providecommand{\indf}{\ensuremath{\beta}}
\providecommand{\indspace}{\ensuremath{\eta}}
\providecommand{\indtime}{\ensuremath{\theta}}
\providecommand{\indrec}{\ensuremath{\varepsilon}}
\providecommand{\honew}{{\sobh1(\W)}} %
\providecommand{\honezw}{{\sobhz1(\W)}} %
\providecommand{\ccr}[2]{\ensuremath{{#1}\vee{#2}}}
\providecommand{\fcc}[2]{\ensuremath{{#1}\wedge{#2}}}

\providecommand{\tm}{{t_m}}
\providecommand{\tms}{{t_m^*}}
\providecommand{\tn}{{t_n}}
\providecommand{\tno}{{t_{n-1}}}

\providecommand{\taun}{\ensuremath{\tau_n}}
\providecommand{\taunm}{\ensuremath{\tau_n^{-1}}}

\providecommand{\hn}{\ensuremath{h_n}}
\providecommand{\hno}{\ensuremath{h_{n-1}}}
\providecommand{\hathn}{\ensuremath{\hat h_n}}

\providecommand{\barf}{\ensuremath{\overline f}}
\providecommand{\projf}[1]{\ensuremath{{\barf}^{#1}}}%

\providecommand{\U}{\ensuremath{U}}%

\providecommand{\un}{\U^n}
\providecommand{\uno}{\U^{n-1}}

\providecommand{\rhoms}{\rho^m_*}
\providecommand{\ms}{{m^*}}
\providecommand{\mso}{{m^*-1}}

\providecommand{\elop}{{\cA}}
\providecommand{\pwA}{\cA_{\mathrm{el}}}%
\providecommand{\discelop}{A}
\providecommand{\disca}[1]{\discelop^{#1}}
\providecommand{\an}{\disca n}%
\providecommand{\ano}{\disca{n-1}}%
\providecommand{\recop}{{\mathcal R}}%

\providecommand{\lprojop}{P_0}
\providecommand{\lproj}[1]{\lprojop^{#1}}%

\providecommand{\eproj}[1]{P_1^{#1}}%

\providecommand{\abil}[2]{\ensuremath{a\left(#1,#2\right)}}

\providecommand{\bdisc}{\ensuremath{\partial}}

\providecommand{\fdisc}{\ensuremath{\overset{\smash{\shortrightarrow}}\partial}}
\providecommand{\sdisc}{\ensuremath{{\partial^2}}}
\providecommand{\discn}{\ensuremath{\overline\partial}}

\usepackage[authoryear]{natbib}
\usepackage{doi}
\usepackage{hyperref}
\setboolean{isamsltex}{true}
\renewcommand{\duality}[2]{\ltwop{#1}{#2}}
\renewcommand{\Id}{\operatorname{I}}
\renewcommand{\maxi}[2]{\max\left(#1,#2\right)}
\renewcommand{\mini}[2]{\min\left(#1,#2\right)}
\renewcommand{\fes}[1]{\fespace[#1]v}
\renewcommand{\fez}[1]{\fezerospace[#1]v}
\renewcommand{\d}{\operatorname d}

\renewcommand{\changes}[1]{{\blue #1}}
\renewcommand{\changes}[1]{#1}
\setboolean{shownotes}{false}
\makeatother
\begin{document}

\newcommand{\mytitle}{%
  ELLIPTIC RECONSTRUCTION AND A POSTERIORI ERROR
  ESTIMATES FOR FULLY DISCRETE LINEAR PARABOLIC PROBLEMS}
\ifthenelse{\boolean{isamsltex}}{%
  \title[Elliptic reconstruction for fully discrete parabolic
    problems]
  \mytitle
  \author{Omar Lakkis}
  \address{Omar Lakkis
    \newline%
    Department of Mathematics
    \newline%
    University of  Sussex
    \newline%
    Brighton, UK-BN1 9RF United Kingdom}
  \curraddr{}
  \email{lakkis.o.maths@gmail.com}
  \thanks{This work is partially supported
    by the E.U. RTN {\em Hyke} HPRN-CT-2002-00282}
  \author{Charalambos Makridakis}
  \address{
    Charalambos Makridakis
    \newline
    Department of Applied Mathematics
    \newline %
    University of Crete
    \newline
    GR-71409 Heraklion, Greece
    \newline
    and
    \newline
    Institute for Applied and Computational Mathematics
    \newline 
    Foundation for Research and Technology-Hellas
    \newline
    Vasilika Vouton P.O.Box 1527
    \newline
    GR-71110 Heraklion, Greece}
  \curraddr{}
  \email{makr@tem.uoc.gr}
  \date{\today}
  \commby{}
  \dedicatory{}
}{%
  \title{\mytitle} \author{Omar Lakkis and Charalambos Makridakis}
  \date{\today}
}
\begin{abstract}
  We derive \aposteriori error estimates for fully discrete
  approximations to solutions of linear parabolic equations on the
  space-time domain.  The space discretization uses finite element
  spaces, that are allowed to change in time.  Our main tool is an
  appropriate adaptation of the elliptic reconstruction technique,
  introduced by Makridakis and Nochetto (2003).  We derive novel
  optimal order \aposteriori error estimates for the maximum-in-time
  and mean-square-in-space norm and the mean-square in space-time of
  the time-derivative norm.
\end{abstract}
\maketitle
\allowdisplaybreaks{
  \section{Introduction}
  \label{sec:introduction} Adaptive mesh refinement methods for
  variational problems have been the object of intense study in recent
  years.  The main objective of these methods is to reduce the
  computational cost in the numerical approximation of PDE's
  solutions. Their usefulness is especially apparent when the exact
  solution has strong, geometrically localized, variations or exhibits
  singularities.  \Aposteriori estimates have proved to be a
  particularly successful mathematical tool in devising efficient
  adaptive versions of many numerical schemes.  In addition,
  \aposteriori estimates provide a new point of view in the theoretical
  investigation of a scheme's behavior.  This is especially important
  for problems where ``reasonable discretizations" do not perform always
  as expected.

  In the context of finite element methods (FEM), the theory of
  \aposteriori estimates for linear stationary problems is by now
  rather mature %
  \citep[and references therein]{AinsworthOden:00:book:A-posteriori,Verfurth:96:book:A-review}.
  
  The situation for nonlinear and time dependent problems, however, has not been
  as thoroughly explored.
  Even for the linear parabolic equation, in spite of advances made in
  the early 1990s
  \citep[e.g.]{ErikssonJohnson:91:article:Adaptive,ErikssonJohnson:95:article:Adaptive2}
  and subsequent ones
  \citep[e.g.]{%
    Picasso:98:article:Adaptive%
    ,%
    Verfurth:03:article:A-posteriori%
    ,%
    ChenFeng:04:article:An-adaptive%
    ,%
    BergamBernardiMghazli:05:article:A-posteriori%
  },
  many issues have yet to be tackled.  Such issues include, for instance,
  the derivation of optimal order estimates in various norms via the
  energy or other direct methods, the use of nonresidual based
  estimators to control the elliptic part of the error, and estimates
  for various time discretization methods.

  In this paper we address some of these issues in the context of
  fully discrete linear parabolic problems where mesh modification,
  with respect to time, might occur---which is natural to expect in
  adaptive schemes for time-dependent problems.  Our main tool in
  deriving the estimates is an appropriate adaptation to the fully
  discrete case of the {\em elliptic reconstruction} technique
  introduced in \citet{MakridakisNochetto:03:article:Elliptic} for the
  model problem of semidiscrete finite element approximations.  A main
  characteristic of this approach, which contrasts with other direct
  techniques found in the literature, is that we can use any available
  \aposteriori estimates for elliptic equations to control the main part
  of the spatial error.  Thus one can take full advantage of a well
  established theory, instead of trying to adapt the estimates
  case-by-case.  This follows from the fact that, in deriving the
  estimates, instead of comparing directly the exact solution with the
  numerical one, we construct an appropriate auxiliary function that
  fulfills two fundamental properties: (i) we know how to estimate its
  difference to the numerical solution via known \aposteriori results,
  and (ii) it satisfies a variant of the original PDE with right-hand
  side that can be controlled \aposteriori in an optimal way.

  In this paper, we combine the elliptic reconstruction technique with
  {\em \aposteriori energy estimates} for the parabolic equation.
  Although residual based \aposteriori estimates using energy methods
  have been established \citep{BergamBernardiMghazli:05:article:A-posteriori%
    ,ChenFeng:04:article:An-adaptive,Picasso:98:article:Adaptive}, it is
  not immediate to figure how to use other than the residual based
  elliptic estimates in them.  For comparison's sake, in this paper we
  will derive residual based energy estimates; but we emphasize the fact
  that the techniques presented in this paper can be easily turned to
  derive estimators for parabolic problems where the ``elliptic part" of
  the error is controlled by non-residual type estimators, e.g.,
  estimators based on the solution of local subproblems.

  The main new results in this paper are \emph{optimal order} a
  posteriori estimates, via energy techniques, in the following spaces
  (anticipating the notation that we will introduce in
  \secref{sse:setting.notation}):
  \begin{enumerate}[(a)\ ]
  \item
    $\leb\infty(0,T;\leb2(\W))$,
  \item
    \label{itm:higher}
    $\leb\infty(0,T;\honezw)$ and
    $\sobh1(0,T;\leb2(\W))$.
  \end{enumerate}
  (As a by-product we recover also known results, such as optimal order
  estimates for the $\leb2(0,T;\honezw)$ norm.)  We address thus the
  open problem of obtaining optimal a posteriori error bounds with
  energy methods in $\leb\infty(0,T;\leb2(\W))$ for fully discrete
  schemes.

  We stress that we do not require any extraneous conditions on the
  variation of the meshes, and the corresponding finite element spaces,
  between successive time-steps; hence the estimates are applicable to
  schemes that are used in practical computations.  An additional
  important feature of the techniques in the present paper is the
  possibility to obtain estimates in the {\em higher order energy
    norms}, mentioned in \eqref{itm:higher} above.  A direct approach in
  deriving the estimates such as the one employed by \citet{Picasso:98:article:Adaptive}
  and \citet{ChenFeng:04:article:An-adaptive} will not work
  here, unless we impose severe mesh conditions or make some \apriori
  assumptions.  Our approach permits to override this difficulty; we
  refer to \secref{sse:direct.approach.failure} for a detailed
  discussion.

  Note that the elliptic reconstruction approach is not restricted by
  the method used for the stability analysis---in this case it being the
  energy method.  Indeed, in a different paper
  \citep{LakkisMakridakisPryer:14:article:A-comparison} we show how the elliptic reconstruction
  can be used in the context of duality techniques such as those proposed by
  \citet{ErikssonJohnson:91:article:Adaptive}.
  \Aposteriori estimates for fully discrete discontinuous Galerkin
  schemes for linear parabolic problems have been obtained via a
  duality method by Eriksson \& Johnson
  \cite{ErikssonJohnson:91:article:Adaptive}.  The results are of
  optimal order, up to a logarithmic factor, in
  $\leb\infty(0,T;\leb2(\W))$.  Backward Euler is the lowest order
  member in this class of methods.  Energy methods have been used by
  other authors for Backward Euler fully discrete approximations to
  parabolic problems
  \cite{Picasso:98:article:Adaptive%
    ,BergamBernardiMghazli:05:article:A-posteriori%
    ,ChenFeng:04:article:An-adaptive}. While these results are of
  optimal order in $\leb2(0,T;\honew)$, they are not so in
  $\leb\infty(0,T;\leb2(\W))$.  Picasso
  \cite{Picasso:98:article:Adaptive} and Chen \& Jia
  \cite{ChenFeng:04:article:An-adaptive} bound, at each time step, the
  spatial indicators by the error (lower bound); their technique is
  based on the ``bubble functions'' technique introduced by Verf\"urth
  for elliptic problems \cite{Verfurth:96:book:A-review}.  Also
  \citet{BergamBernardiMghazli:05:article:A-posteriori} have
  established lower bounds and their estimators have the additional
  feature of decoupled spatial and temporal error.  
  \citet{ChenFeng:04:article:An-adaptive} provide an adaptive algorithm
  that is rigorously proved to guarantee the control of the error,
  under the assumption of the termination of their algorithm.
  Various other \aposteriori estimates for semidiscrete/fully discrete
  approximations to linear and nonlinear parabolic problems in various
  norms are found in the literature
  \citep{AdjeridFlahertyBabuska:99:article:A-posteriori%
    ,BabuskaFeistauerSolin:01:article:On-one-approach%
    ,BabuskaOhnimus:01:article:A-posteriori%
    ,FrutosNovo:02:article:Postprocessing-the-linear%
    ,ErikssonJohnson:95:article:Adaptive4%
    ,NochettoSavareVerdi:00:article:A-posteriori%
    ,NochettoSchmidtVerdi:00:article:A-posteriori%
    ,Verfurth:98:article:A-posteriori-L%
    ,Verfurth:98:article:A-posteriori-W%
  }.  In particular,
  \citet{BabuskaFeistauerSolin:01:article:On-one-approach} have
  derived estimates in $\leb2(0,T;\leb2(\W))$ \citep[see
    also][]{AdjeridFlahertyBabuska:99:article:A-posteriori%
    ,BabuskaOhnimus:01:article:A-posteriori%
  }.
  \citet{Verfurth:98:article:A-posteriori-L,Verfurth:98:article:A-posteriori-W}
  showed estimates in $\leb r(0,T;\leb\rho(\W)),$ with
  $1<r,\rho<\infty $ for certain fully discrete approximations of
  certain quasilinear parabolic equations.
  \citet{LakkisNochetto:05:article:A-posteriori} used ad-hoc geometric
  energy norms to derive conditional a posteriori estimates for
  quasilinear equations such as the mean curvature flow of graphs.
  \citet{FrutosNovo:02:article:Postprocessing-the-linear} proved a
  posteriori error estimates of the {$p$}-version of space discrete
  schemes for parabolic equations; a similar function to the elliptic
  reconstruction and its improved approximation properties is used
  \citep[see also][]{Garcia-ArchillaTiti:00:article:Postprocessing}.
  Finally, for applications of suitable reconstructions to time
  discretizations of various type, we refer to
  \citet{AkrivisMakridakisNochetto:06:article:A-posteriori%
    ,MakridakisNochetto:06:article:A-posteriori%
  }.
  \subsection{Problem setting and notation}
  \label{sse:setting.notation}
  Let us focus the discourse now and introduce the fully discrete
  scheme, which will be the object of our analysis.  We start with the
  exact problem. Let $\W$ be a bounded domain of the Euclidean space
  $\rR^d$, $d\in\Nat$ and $T\in\rR^+$.
  We assume throughout the paper that \W is a polygonal convex domain,
  noticing that all the results can be extended to certain non-convex
  domains, like domains with reentrant corners in $d=2$.
  Since the difficulties in the analysis bellow in the case of other
  boundaries are mainly coming from the elliptic part of the error, the
  reader interested in \aposteriori error estimates for curved
  boundaries is referred to Dörfler \& Rumpf \cite{DorflerRumpf:98:article:An-adaptive}.

  We will consider the problem of finding a finite element
  approximation of the solution $u\in\leb\infty(0,T;\honezw)$, with
  $\partial_tu\in\leb2(0,T;\leb2(\W))$, to the linear parabolic
  problem
  \begin{equation}
    \begin{split}\label{eqn:continuous.heat}
      &\duality{\partial_tu}{\phi}+\abil u\phi=\ltwop{f}{\phi}
      \quad\forall\phi\in\honezw,\\ &\text{ and }u(0)=g,
    \end{split}
  \end{equation}
  where $f\in\leb2(\W\times(0,T))$ and $g\in\honezw$, and $a$ is a
  bilinear form on $\honezw$ defined by
  \begin{equation}
    \label{eqn:bilinear.form}
    \abil v\psi :=\ltwop{\vec A\grad v}{\grad\psi},\:\forall
    v,\psi\in\honezw
  \end{equation}
  where ``$\grad$'' denotes the spatial gradient and the matrix $\vec
  A\in\leb\infty(\W)^{d\times d}$ is such that
  \begin{gather}
    \label{eqn:bounded.bilinear}
    \abil\psi\phi \leq \beta\norm{\psi}_1\norm{\phi}_1,
    \quad\forall\phi,\psi\in\honezw,\\
    \label{eqn:coercive.bilinear}
    \abil\phi\phi \geq \alpha\norm{\phi}_1^2,
    \quad\forall\phi\in\honezw,
  \end{gather}
  with $\alpha,\beta\in\reals^+$.  Whenever not stated explicitly, we
  assume that the data $f,g,\vec A$ and the solution $u$ of the above
  problem are sufficiently regular for our purposes.
  
  Here and subsequently, for a given Lebesgue measurable set
  $D\subset\rR^d$, we use the common notation
  \begin{align}
    &\ltwop\phi\psi_D:=\int_D\phi(\vec x)\psi(\vec x)\d\mu(\vec x),
    \\ &\Norm{\phi}_D :=\Norm{\phi}_{\leb2(D)} :=\ltwop\phi\phi_D^{1/2},
    \\ &\norm{\phi}_{k,D} :=\Norm{\D^k\phi}_D,\text{ for }k\in\Nat
    \quad(\text{with }\D^1\phi:=\grad\phi,\text{ etc.}),
    \\ &\Norm{\phi}_{k,D}
    :=\bigg(\Norm{\phi}_D^2+\sum_{j=1}^k\norm{\phi}_{j,D}^2\bigg)^{1/2},
    \text{ for }k\in\Nat,
  \end{align}
  where $\d\mu(\vec x)$ is either the Lebesgue measure element $\d\vec
  x$, if $D$ is has positive such measure, or the $(d-1)$-dimensional
  (Hausdorff) measure $\ds(\vec x)$, if $D$ has zero Lebesgue measure.
  
  In many instances, in order to compress notation and when there is
  no danger of engendering confusion, we drop altogether the
  ``differential'' symbol from integrals; this applies also to
  integrals in time.
  
  We use the standard function spaces $\leb2(D)$, $\sobh k(D)$,
  $\sobhz k(D)$ and denote by $\sobh{-1}(D)$ the dual space of
  $\sobhz1(D)$ with the corresponding pairing written as
  $\duality{\cdot}{\cdot}_D$.  We omit the subscript $D$ whenever
  $D=\W$.  We denote the Poincaré constant relative to $\W$ by
  $C_{2,1}$ and take $\norm{\cdot}_1$ to be the norm of $\honezw$.  We
  use the usual duality identification
  \begin{equation}
    \honezw\subset\leb2(\W)\sim\leb2(\W)'\subset\sobh{-1}(\W)
  \end{equation}
  and the dual norm
  \begin{equation}\label{eqn:dual.norm}
    \Norm\psi_{-1}
    :=\sup_{0\neq\phi\in\sobhz1(\W)}\frac{\duality\psi\phi}{\norm{\phi}_1}
    \left(=\sup_{0\neq\phi\in\sobhz1(\W)}\frac{\ltwop\psi\phi}{\norm{\phi}_1},
    \text{ if }\psi\in\leb2(\W)\right).
  \end{equation}
  We also use the {\em energy norm} $\norm{\cdot}_a$ defined as
  \begin{equation}
    \norm{\phi}_a:=\abil\phi\phi^{1/2},\:\forall\phi\in\honezw.
  \end{equation}
  It is equivalent to the norm $\norm{\cdot}_1$ on the space
  $\sobhz1(\W)$, in view of \eqref{eqn:bounded.bilinear} and
  \eqref{eqn:coercive.bilinear}.  In particular, we will often use the
  following inequality
  \begin{equation}
    \label{eqn:one.energy.inequality}
    \norm\phi_1\leq\alpha^{-1/2}\norm\phi_a,\:\forall\phi\in\honezw.
  \end{equation}

  In order to discretize the time variable in
  \eqref{eqn:continuous.heat}, we introduce the partition
  $0=t_0<t_1<\ldots<t_N=T$ of $[0,T]$.  Let
  $I_n:=\clopinter{t_{n-1}}{t_n}$ and we denote by
  $\tau_n:=t_n-t_{n-1}$ the time steps.  We will consistenly use the
  following ``superscript convention'': whenever a function depends on
  time, e.g. $f(\vec x, t)$, and the time is fixed to be
  $t=t_n,\,n\in\fromto0N$ we denote it by $f^n(\vec x)$.  Moreover, we
  often drop the space dependence explictly, e.g, we write $f(t)$ and
  $f^n$ in reference to the previous sentence.

  We use a standard FEM to discretize the space variable.  Let
  $(\cT_n)_{n\in\fromto0N}$ be a family of conforming triangulations
  of the domain $\W$ \cite{BrennerScott:94:book:The-mathematical,Ciarlet:02:book:The-finite}.  These
  triangulations are allowed to change arbitrarily from a timestep to
  the next, as long as they maintain some very mild {\em
    compatibility} requirements.  Our use of the term
  ``compatibility'' is precisely defined in Appendix
  \ref{sec:compatible.triangulations}; it is an extremely mild
  requirement which is easily implemented in practice.

  For each given a triangulation $\cT_n$, we denote by $h_n$ its
  meshsize function defined as
  \begin{equation}
    h_n(\vec x)=\diam(K),\text{ where }{K\in\cT_n}\text{ and }{\vec x\in
      K},
  \end{equation}
  for all $\vec x\in\W$.  We also denote by $\cS_n$ the set of {\em
    internal sides} of $\cT_n$, these are edges in $d=2$---or faces in
  $d=3$---that are contained in the interior of \W; the {\em interior
    mesh of edges} $\Sigma_n$ is then defined as the union of all
  internal sides $\cup_{E\in\cS_n}E$.  We associate with these
  triangulations the {\em finite element spaces}:
  \begin{gather}
    \fez n:=\ensemble{\phi\in\honezw}{\forall
      K\in\cT_n:\restriction\phi K\in\poly\ell}
  \end{gather}
  where $\poly\ell$ is the space of polynomials in $d$ variables of
  degree at most $\ell\in\Nat$.
  Given two successive compatible triangulations $\cT_{n-1}$ and
  $\cT_n$, we define $\hat h_n:=\maxi{h_n}{h_{n-1}}$ (see Appendix
  \ref{sec:compatible.triangulations} and \ref{sec:inequalities}).  We
  will also use the sets $\hat\Sigma_n:=\Sigma_n\cap\Sigma_{n-1}$ and
  $\check\Sigma_n:=\Sigma_n\cup\Sigma_{n-1}$.  Similarly, in
  \secref{sec:lih1}, we will use
  $\hhat\Sigma_n:=\Sigma_n\cap\Sigma_{n-1}\cap\Sigma_{n+1}$,
  $\ccheck\Sigma_n:=\Sigma_n\cup\Sigma_{n-1}\cup\Sigma_{n+1}$, $\hhat
  h_n:=\max_{i\in\fromto{-1}{1}}h_{n+i}$.

  \begin{Def}[fully discrete scheme]
    \label{def:fully.discrete}
    The standard backward Euler-Galerkin method for the discretization
    of problem \eqref{eqn:continuous.heat} associated with the finite
    element spaces $\fez n$, leads to the following recursive
    {\em fully discrete scheme}:
    \begin{equation}
      \label{eqn:fully.discrete.scheme}
      \begin{split}
        &U^0:=I^0u(0),\text{ and } \\ &\taunm\ltwop{\un-\uno}{\phi_n}
        +\abil{\un}{\phi_n}=\ltwop{f^n}{\phi_n},
        \quad\forall\phi_n\in\fez n,\text{for $n\in\fromto1N$}.
      \end{split}
    \end{equation}
    Here the operator $I^0$ is some suitable interpolation or projection
    operator from $\honezw$, or $\leb2(\W)$, into $\fez n$.
  \end{Def}
  In the sequel we shall use a continuous piecewise linear extension in
  time of the sequence $(\un)$ which we denote by $U(t)$ for $t\in[0,T]$
  (see \secref{def:discrete.time} for the precise definition).

  \subsection{A posteriori estimates and reconstruction operators.}
  Suppose we associate with $U$ an auxiliary function
  $\w:[0,T]\rightarrow\honezw$, in such a way that the {\em total error}
  \begin{equation}
    e:=U-u
  \end{equation}
  can be decomposed as follows
  \begin{gather}
    e=\rho-\epsi\\ \epsi:=\w-U,\quad\rho:=\w-u.
  \end{gather}
  The new auxiliary function $\w$ is {\em reconstructed} from the given
  approximation $U$. The success of this splitting for the estimates to
  follow rests in the following properties:
  \begin{NumberList}
  \item The error $\epsi$ is easily controlled by a posteriori
    quantities of optimal order.
  \item The error $\rho$ satisfies a modification of the original PDE
    whose right-hand side depends on $\epsi$ and $U$.  This right-hand
    side can be bounded a posteriori in an optimal way.
  \end{NumberList}
  Therefore in order to successfully apply this idea we must select a
  suitable reconstructed function $\w$.  In our case, this choice is
  dictated by the elliptic operator at hand; the precise definition is
  given in \secref{def:elliptic.reconstruction}.  In addition the effect
  of mesh modification will reflect in the right-hand side of the
  equation for $\rho$.  As a result of our choice for $\w$ we are able
  to derive optimal order estimators for the error in
  $\leb\infty(0,T;\leb2(\W))$, as well as in $\leb\infty(0,T;\honezw)$
  and $\sobh1(0,T;\leb2(\W))$. In addition, our choosing $\w$ as the
  elliptic reconstruction will have the effect of separating the spatial
  approximation error from the time approximation as much as possible.
  We show that the spatial approximation is embodied in $\epsi$ which
  will be referred to as the {\em elliptic reconstruction error} whereas
  the time approximation error information is conveyed by $\rho$, a fact
  that motivates the name {\em main parabolic error} for this term.
  This ``splitting'' of the error is already apparent in the spatially
  discrete case \cite{MakridakisNochetto:03:article:Elliptic}.

  With the above notation, we prove in the sequel that $\rho$ satisfies
  the following variational equation.
  \begin{Lem}[main parabolic error equation]
    \label{lem:main.error.equation}
    For each $n\in\fromto1N$, and for each $\phi\in\sobhz1(\W)$,
    \begin{equation}\label{eqn:main.error.equation}
      \begin{split}
        \duality{\partial_t\rho}{\phi}+\abil{\rho}{\phi}
        =&\ltwop{\partial_t\epsi}{\phi}+\abil{\w-\w^n}{\phi}\\ &+\ltwop{\lproj
          nf^n-f}{\phi} +\taunm\ltwop{\lproj n\uno-\uno}{\phi} \text{ on
          $I_n$.}
      \end{split}
    \end{equation}
    Here $\lproj{n}$ denotes the $L^{2}$-projection into $\fez n .$
  \end{Lem}
  \subsection{How do we derive the estimates? A general overview}
  \label{sse:general.strategy}
  Identity \eqref{eqn:main.error.equation} along with the properties of
  the elliptic reconstruction will allow us to obtain \aposteriori error
  estimates in different norms.  We start from the afore-mentioned
  splitting of the error
  \begin{equation}
    \Norm{e(t)}_{X} \leq\Norm{\epsi(t)}_{X} +\Norm{\rho(t)}_{X}
  \end{equation}
  where $X$ is any suitable space of functions on $\W$.  The choice of
  $X$ might depend on the applications that are in mind, and on the
  ability to bound the terms on the right-hand side.
  Therefore the following observations are fundamental:
  \begin{NumberList}
  \item 
    The first term, $\Norm{\epsi(t)}_X=\Norm{\w(t)-U(t)}_{X}$,
    can be bounded by appropriate a posteriori error estimates for
    elliptic problems.  Indeed, at a time $t_n$, $\w^n$ ($\un$) is the
    exact ($\fez n$-finite element) solution of the elliptic problem
    $\elop v=\an\un$ (with reference to the notation in the
    forthcoming \secref{obs:representation}).  To obtain an a
    posteriori estimate for this term it is enough assume that: (a)
    elliptic \aposteriori error estimates are provided for the
    $X$-norm, through estimator functions that depend on $\un$,
    $\an\un$, the triangulation parameters and the polynomial
    degree---many such estimator functions are available from standard
    a posteriori error analysis for elliptic problems
    \cite[e.g.]{AinsworthOden:00:book:A-posteriori,Braess:01:book:Finite,BrennerScott:94:book:The-mathematical,Verfurth:96:book:A-review};
    and (b) that the ``data'', such as $\an\un$, can be explicitly
    obtained in terms of $\un$---this is possible in our context as we
    observe in \secref{obs:computing.rhs}.  Notice that
    $\Norm{\epsi(t)}_X$ contains exclusively spatial error effects.
  \item
    The second term, $\Norm{\rho(t)}_X=\Norm{\w(t)-u(t)}_{X}$, is
    estimated by an appropriate use of \eqref{eqn:main.error.equation}.
    By the form of its right-hand side we conclude that the resulting
    estimators will include quantities measuring the space error, the
    time error, the variation of $f$ and the effect of mesh changes with
    respect to $n$.
  \end{NumberList}
  \subsection{Comparison with the direct approach}
  The identity \eqref{eqn:main.error.equation} can be appreciated if we
  compare it with the error equation that one obtains from a direct
  comparison of $u$ and $U$.  In the direct approach the error relation
  is given by
  \begin{equation}
    \label{eqn:direct.approach}
    \duality{\partial_t e}{\phi}+\abil e\phi =\ltwop{\partial_t
      U}{\phi}+\abil{U^n}\phi-\ltwop f\phi +\abil{U-U^n}\phi.
  \end{equation}
  Using the fact that $U$ is the solution of the fully discrete scheme
  one sees that
  \begin{equation}
    \label{eqn:direct.approach.orthogonality}
    \begin{split}
      \duality{\partial_te}{\phi}+\abil{e}\phi
      =&
      \ltwop{\partial_tU}{\phi-\phi_n}
      +
      \abil{U^n}{\phi-\phi_n}
      \\
      &
      -\ltwop{f^n}{\phi-\phi_n}
      \\
      &
      -\ltwop{f-f^n}\phi
      +\abil{U-U^n}\phi,\text{ for }\phi_n\in\fez n.
    \end{split}
  \end{equation}

  A comparison between \eqref{eqn:main.error.equation} and
  \eqref{eqn:direct.approach.orthogonality} demonstrates the two main
  differences in the corresponding approaches. Equation
  \eqref{eqn:direct.approach.orthogonality} has all the information from
  the numerical scheme ``built-in'', in particular it satisfies the
  Galerkin orthogonality property; therefore the error (and stability)
  analysis is dictated by the choices of both $\phi $ and $\phi _n$.  In
  the case of the elliptic reconstruction, the Galerkin orthogonality
  property {\em is not used explicitly} in the analysis.  The fact that
  $U$ is solution of the discrete scheme is used only {\em implicitly}
  through the reconstruction $\w$: firstly it is used to estimate
  $\epsi$ and $\partial_t\epsi$ (this estimate comes ``for free'' from
  the elliptic a posteriori theory precisely because of the definition
  of $\w $); secondly, it is used to derive
  \eqref{eqn:main.error.equation}, which allows then to estimate
  $\Norm{\rho}_X$ and $\Norm{\partial_t\rho}_X$, in terms of time and
  data approximation estimators, and spatial estimators dictated only by
  $\epsi$.

  A second important difference, between these two approaches, is the
  presence of the---suboptimal in $\leb2(\W)$---term
  $\abil{U^n}{\phi-\phi_n}$ in
  \eqref{eqn:direct.approach.orthogonality}.  Because of this term,
  the direct approach fails to provide optimal order a posteriori
  estimators in $\leb\infty(0,T;\leb2(\W)) ,$
  \cite{ChenFeng:04:article:An-adaptive,Picasso:98:article:Adaptive};
  the same problem can appear also in a different context
  \cite{BergamBernardiMghazli:05:article:A-posteriori}.  On the other
  hand, it is interesting to note that, the presence of this term is
  also the reason why the direct approach fails to lead satisfactory
  \aposteriori results in the { higher order energy norms}, see
  \secref{sse:direct.approach.failure} for the details.

  Notice, concerning the time discretization error, that the term
  $\abil{U-U^n}\phi$ in \eqref{eqn:direct.approach.orthogonality} is
  very similar to the term $\abil{\w-\w^n}\phi$ in
  \eqref{eqn:main.error.equation}.

  \subsection{Outline}
  The rest of the paper is organized as follows.  We introduce the
  necessary discrete and continuous operators in order to define the
  reconstruction $\w$ in \secref{sec:reconstruction} and state some of
  its basic properties needed in the sequel.

  In \secref{sec:lil2} we provide the \aposteriori analysis in
  $\leb\infty(0,T;\leb2)$.  The estimators in Theorem
  \ref{the:lil2.estimate} are of optimal order and residual type. We
  carry out the analysis with a particular class of ``elliptic"
  estimators in mind.  As mentioned earlier, other choices are possible;
  but in order to use them, the arguments related to the terms involving
  mesh-change effects must be appropriately adapted.

  Next, in \secref{sec:lih1}, we show a posteriori estimates of optimal
  order in the { higher order energy norms} $\leb\infty(0,T;\honezw)$
  and $\sobh1(0,T;\leb2(\W))$.  This case is of particular interest, as
  simplified situation for a class of nonlinear degenerate parabolic
  problems where lower order energy estimates are not available
  \cite{LakkisNochetto:05:article:A-posteriori}.  In this section, we also discuss briefly
  why the direct approach cannot be applied successfully in this case,
  \secref{sse:direct.approach.failure}.

  Finally, in \secref{sec:numerics} we complement the theory with a
  numerical experimentation.  In particular, we show that the estimators
  derived in \secref{sec:lil2} have the expected, optimal, experimental
  order of convergence in the numerical test.

  In the Appendix section we collect some useful facts about the concept
  of compatible triangulations in
  \secref{sec:compatible.triangulations}, elliptic shift---or
  regularity---inequalities in \secref{sse:elliptic.regularity},
  interpolation opertors and inequalities in
  \S\S\ref{sse:scott.zhang.inequalities}--\ref{sse:mesh.change.interpolation.inequalities}
  and our convention on constant's labeling in
  \secref{sse:combining.constants}.

  We refer the reader interested in the practical aspects of our
  estimators to a second paper in which detailed numerical experiments,
  including comparisons with estimators coming from duality, as well as
  numerical investigation of the effect of mesh modification via
  \aposteriori estimators are included \cite {LakkisMakridakisPryer:14:article:A-comparison}.

  \subsection*{Acknowledgment} We would like to thank Ricardo~Nochetto,
  for many useful discussions about this paper's results, and
  Andreas~Veeser for his interesting remarks.
  \section{The elliptic reconstruction: definition and preliminaries}
  \label{sec:reconstruction}
  We introduce basic tools and the elliptic reconstruction.  Although
  the definitions in the first part of this section are independent of
  the time discretization and could be applied to any finite element
  space, we still use the space $\fez n$ defined in the introduction.
  \begin{Def}[representation of the elliptic operator, 
      discrete elliptic operator, projections]
    \label{obs:representation}
    Suppose a function $v\in\fez n$, the bilinear form can be then
    represented as
    \begin{equation}\label{eqn:representation.long.notation}
      \abil v\phi = \sum_{K\in\cT_n} \ltwop{\elop v}{\phi}_K +
      \sum_{E\in\cS_n}\ltwop{J[v]}\phi_E, \:\forall\phi\in\honezw,
    \end{equation}
    where $J[v]$ is the {\em spatial jump of the field $\vec A\grad v$
      across an element side $E\in\cS_n$} defined as
    \begin{equation}
      \restriction{J[v]}{E}(\vec x)=\jump{\vec A\grad v}_E(\vec x)
      :=\lim_{\ep\rightarrow0}\left(\vec A\grad v(\vec x+\ep\vec\nu_E)-
      \vec A\grad v(\vec x-\ep\vec\nu_E)\right)\inner\vec\nu_E
    \end{equation}
    where $\vec\nu_E$ is a choice, which does not influence this
    definition, between the two possible normal vectors to $E$ at the
    point $\vec x$.
    
    Since we use the representation
    \eqref{eqn:representation.long.notation} quite often, we introduce
    now a practical notation that makes it shorter and thus easier to
    manipulate in convoluted computations.  For a finite element
    function, $v\in\fez n$ (or more generally for any Lipschitz
    continuous function $v$ that is $\cont2(\interior(K))$, for each
    $K\in\cT_n$), denote by $\pwA v$ the {\em regular part} of the
    distribution $\elop v$, which is defined as a piecewise continuous
    function such that
    \begin{equation}
      \ltwop{\pwA v}\phi = \sum_{K\in\cT_n}\ltwop{\elop v}\phi,
      \quad\forall\phi\in\honezw.
    \end{equation}
    The operator $\pwA $ is sometime refered to as the \emph{elementwise
    elliptic operator}, as it is the result of the application of
    $\elop{\cdot}$ only on the interior of each element $K\in\cT_n$.
    This observation justifies our subscript in the notation.  We shall
    write the representation \eqref{eqn:representation.long.notation} in
    the shorter form
    \begin{equation}\label{eqn:representation.short.notation}
      \abil v\phi = \ltwop {\pwA v}\phi + \ltwop{J[v]}\phi_{\Sigma_n},
      \quad\forall\phi\in\honezw.
    \end{equation}
  \end{Def}
  \begin{Def}[discrete elliptic operator]  
    Let us now recall some more basic definitions that we will be using.
    The {\em discrete elliptic operator} associated with the bilinear form
    $a$ and the finite element space $\fez n$ is the operator
    $\an:\honezw\rightarrow\fez n$ defined by
    \begin{equation}
      \label{eqn:discrete.elliptic.operator}
      \ltwop{\an v}{\phi_n}=\abil{v}{\phi_n},\quad\forall\phi_n\in\fez
      n,
    \end{equation}
    for $v\in\honezw$.
  \end{Def}
  The {\em $\leb2$-projection operator} is defined as the operator
  $\lproj n: \leb2(\W)\rightarrow\fez n$ such that
  \begin{equation}
    \ltwop{\lproj
      nv}{\phi_n}=\ltwop{v}{\phi_n},\quad\forall\phi_n\in\fez n
  \end{equation} 
  for $v\in\leb2(\W)$; and the {\em elliptic projection operator}
  $\eproj n: \honezw\rightarrow\fez n$ is defined by
  \begin{equation} \abil{\eproj
      nv}{\phi_n}=\abil{v}{\phi_n},\quad\forall\phi_n\in\fez n.
  \end{equation} 

  The elliptic reconstruction, which we define next, is a partial right
  inverse of the elliptic projection
  \cite{MakridakisNochetto:03:article:Elliptic}. (Notice that a similar operator has
  been introduced by Garc\'\i a-Archilla \& Titi, albeit with different
  applications in mind than ours \cite{Garcia-ArchillaTiti:00:article:Postprocessing}.)
  \begin{Def}[elliptic reconstruction]
    \label{def:elliptic.reconstruction}
    We define the {\em elliptic reconstruction operator} associated with
    the bilinear form $a$ and the finite element space $\fez n$ to be
    the unique operator $\cR^n:\honezw\rightarrow\honezw$ such that
    \begin{equation}
      \label{eqn:reconstruction}
      \abil{\cR^nv}{\phi}=\ltwop{\an v}{\phi},\quad\forall\phi\in\honezw,
    \end{equation} for a given $v\in\honezw$. The function $\recop n v$ is
    referred to as the {\em elliptic reconstruction} of $v$.  A crucial
    property of $\recop n$ is that $v-\recop nv$ is orthogonal, with respect
    to $a$, to \fez n:
    \begin{equation}\label{eqn:reconstruction.orthogonality}
      \abil{v-\recop n v}{\phi_n}=0,\quad\forall\phi_n\in\fez n.
    \end{equation}
    From this property and by now standard techniques in \aposteriori
    error estimates for elliptic problems
    \cite{AinsworthOden:00:book:A-posteriori,Braess:01:book:Finite,Verfurth:96:book:A-review}, it is possible
    to obtain the following result.
  \end{Def}
  \changes{
    \begin{Lem}[elliptic reconstruction error estimates]
      \label{lem:elliptic.aposteriori.estimates}
      For any $v\in\fez n$ the following estimates hold true
      \begin{align}
        \label{eqn:elliptic.estimate.h1}
        \norm{\recop nv-v}_1 &\leq\frac{C_{3,1}}\alpha \Norm{(\pwA v-\an
          v)\hn}+\frac{C_{5,1}}\alpha \Norm{J[v]\hn^{1/2}}_{\Sigma_n}
        ,\\
        \label{eqn:elliptic.estimate.l2}
        \Norm{\recop nv-v} &\leq C_{6,2}\Norm{(\pwA v-\an
          v)\hn^2}+C_{10,2}\Norm{J[v]\hn^{3/2}}_{\Sigma_n},
      \end{align}
      where the constants $C_{k,j}$ are defined in Appendix
      \ref{sec:inequalities}.
    \end{Lem}
  }%
  \begin{Def}[discrete time extensions and derivatives]
    \label{def:discrete.time}
    Given any discrete function of time---that is, a sequence of values
    associated with each time node $t_n$---e.g., $(\un)$, we associate to
    it the continuous function of time defined by the Lipschitz continuous
    piecewise linear interpolation, e.g.,
    \begin{equation}
      \label{eqn:pwl.extension}
      U(t):=l_{n-1}(t)\uno+l_n(t)\un, \text{ for $t\in I_n$ and
        $n\in\fromto1N$};
    \end{equation}
    where the functions $l$ are the hat functions defined by
    \begin{equation}
      l_n(t):=\frac{t-t_{n-1}}{\tau_n}\one_{I_n}(t)-\frac{t-t_{n+1}}{\tau_{n+1}}\one_{I_{n+1}}(t),
      \text{ for $t\in[0,T]$ and $n\in\fromto0N$},
    \end{equation}
    $\one_X$ denoting the characteristic function of the set $X$.
    The {\em time-dependent elliptic reconstruction} of $U$ is the
    function
    \begin{equation}
      \w(t):=l_{n-1}(t)\recop {n-1}\uno+l_n(t)\recop n\un, \text{ for $t\in
        I_n$ and $n\in\fromto1N$}.
    \end{equation}
    Notice that $\w$ is a Lipschitz continuous function of time.

    We introduce the following definitions whose purpose is to make
    notation more compact:
    \begin{LetterList}
    \item 
      {\em Discrete (backward) time derivative}
      \begin{equation}
        \bdisc\un:=\frac{\un-\uno}{\taun}
      \end{equation}
      Notice that $\bdisc\un=\partial_tU(t)$, for all $t\in I_n$, hence we
      can think of $\bdisc\un$ as being the value of a discrete function at
      $t_n$.  We thus define $\bdisc U$ as the piecewise linear extension of
      $(\bdisc\un)_n$, as we did with $U$.
    \item {\em Discrete (centered) second time derivative}
      \begin{equation}
        \sdisc\un:=\frac{\bdisc U^{n+1}-\bdisc\un}{\taun}.
      \end{equation}
    \item  
      {\em Averaged ($\leb2$-projected) discrete time derivative}
      \begin{equation}
        \label{eqn:averaged.discrete.time.derivative}
        \discn\un := \lproj n\bdisc\un = \frac{\un-\lproj n\uno}{\taun},
        \quad\forall n\in\fromto 1N.
      \end{equation}
      The reason we introduce this notation for is that $\bdisc\un$, in
      general, does not belong to the current finite element space, $\fez
      n$, whereas $\discn\un$ does.
    \item {The \em $\leb2$-projection of $f^n$}
      \begin{equation}
        \projf n := \lproj n f^n.
      \end{equation}
      Since this is a discrete function of time, consistently with notation
      \eqref{eqn:pwl.extension}, we denote by ${\barf}$ the piecewise linear
      interpolation of $(\projf n)_n$.
    \end{LetterList}
  \end{Def}

  \begin{Obs}[pointwise form]
    The discrete elliptic operators $\an$ can be employed to write the
    fully discrete scheme \eqref{eqn:fully.discrete.scheme} in the
    following {\em pointwise form}
    \begin{equation}
      \label{eqn:pointwise.form}
      \discn\un(\vec x) + \an\un(\vec x) = \projf n(\vec
      x),\:\forall\vec x\in\W.
    \end{equation}
    Indeed, in view of $\discn\un+\an\un-\projf n\in\fez n$,
    \eqref{eqn:fully.discrete.scheme}, and
    \eqref{eqn:discrete.elliptic.operator}, we have
    \begin{equation}
      \begin{split}
        \ltwop{\an\un+\discn\un-\projf n}{\phi}
        &=\ltwop{\an\un+\discn\un-\projf n}{\lproj
          n\phi}\\ &=\abil{\un}{\lproj
          n\phi}+\ltwop{\taunm(\un-\uno)-f^n}{\lproj n\phi} =0,
      \end{split}
    \end{equation}
    for any $\phi\in\honezw$.  Therefore the function
    $\discn\un+\an\un-\projf n$ vanishes.
  \end{Obs}

  \subsection{Proof of Lemma \ref{lem:main.error.equation}}
  The definitions in \ref{def:discrete.time} and
  \eqref{eqn:pointwise.form} yield
  \begin{equation}
    \begin{split}
      \ltwop{\bdisc\un+\an\un-\lproj nf^n}{\phi} &-\taunm\ltwop{\lproj
        n\uno-\uno}{\phi}\\ &=\ltwop{\discn\un+\an\un-\projf
        n)}{\phi}=0,
    \end{split}
  \end{equation}
  for each $\phi\in\honezw$ and $n\in\fromto1N$.  In view of the
  elliptic reconstruction definition we obtain
  \begin{equation}
    \label{eqn:fully.discrete.global}
    0=\ltwop{\bdisc\un}{\phi}+\abil{\w^n}{\phi} -\ltwop{\lproj
      nf^n}{\phi}-\taunm\ltwop{\lproj n\uno-\uno}{\phi}.
  \end{equation}
  On the other hand \eqref{eqn:continuous.heat} implies
  \begin{equation}
    \begin{split}
      \duality{\partial_t\rho}\phi+\abil\rho\phi
      =\ltwop{\partial_t\w}\phi+\abil\w\phi-\ltwop f\phi,
    \end{split}
  \end{equation}
  from which we subtract equation \eqref{eqn:fully.discrete.global} and
  obtain \eqref{eqn:main.error.equation}.\hfill\qed

  \begin{Obs}[How to compute $\an\un$?]
    \label{obs:computing.rhs}
    The computer evaluation of the estimators appearing in
    \eqref{eqn:elliptic.estimate.h1}--\eqref{eqn:elliptic.estimate.l2}
    involves an exact differentiation, for $\pwA \un$, and a jump
    evaluation, for $J[\un]$, which are straightforward to implement for
    polynomial finite elements.  Notice that the computation also involves
    the term $\an\un$.  In the context of the elliptic equation $\elop
    u=g$ this term is $\lproj ng$ (the $\leb2$-projection of g), whereas
    in the context of the heat equation~\eqref{eqn:continuous.heat} this
    term can be easily obtained from the discrete equation in pointwise
    form~\eqref{eqn:pointwise.form}.  This observation motivates the next
    definition.
  \end{Obs}

  \begin{Def}[residuals]
    The residuals constitute the building blocks of the \aposteriori
    estimators used in this paper.  We associate with equations
    \eqref{eqn:continuous.heat} and \eqref{eqn:pointwise.form} two
    residual functions: the {\em inner residual} is defined as
    \begin{equation}
      \label{eqn:inner.residual.function}
      \begin{split}
        R^0 &:= \pwA U^0-A^0U^0,\\ R^n &:= \pwA \un-\an\un = \pwA \un-\projf
        n+\discn\un,\text{ for }n\in\fromto1N,
      \end{split}
    \end{equation}
    and the {\em jump residual} which is defined as
    \begin{equation}
      \label{eqn:jump.residual.function}
      J^n := J[\un].
    \end{equation}
    We note that, with definition \secref{obs:representation} in mind, the
    inner residual terms can be written in the following, more familiar
    but also more cumbersome, fashion
    \begin{equation}
      \ltwop{R^n}\phi =\sum_{K\in\cT_n} \ltwop{\elop v-\lproj nf(t_n)
        +\frac{U^n-\lproj nU^{n-1}}{\taun}}\phi _K.
    \end{equation}
  \end{Def}

  \section{\Aposteriori error estimates in the lower order norms}
  \label{sec:lil2}
  We start by introducing the following error estimators that are local
  in time.  The full estimators, that will appear in Theorem
  \ref{the:lil2.estimate}, are accumulations in time of these local
  estimators.  The accumulations, which can be of $\leb1$, $\leb2$ or
  $\leb\infty$ type, are anticipated by the first subscript in the
  estimators.  \changes{
    \begin{Def}[$\leb\infty(\leb2)$ and $\leb2(\sobh1)$ error estimators]
      We introduce, for $n\in\fromto0N$, the {\em elliptic reconstruction
        error estimators}
      \begin{gather}
        \label{eqn:lil2.indrec}
        \indrec_{\infty,n} :=C_{6,2}\Norm{\hn^2 R^n
        }+C_{10,2}\Norm{\hn^{3/2} J^n }_{\Sigma_n} ,\\
        \indrec_{2,n}:=\frac{C_{3,1}}\alpha\Norm{\hn R^n}
        +\frac{C_{5,1}}\alpha\Norm{\hn^{1/2} J^n}_{\Sigma_n},
        \intertext{and, for $n\in\fromto1N$, the {\em space error
            estimator}}
        \label{eqn:lil2.indspace.1}
        \indspace_{1,n} :=C_{6,2}\Norm{\hathn^2\bdisc R^n}
        +C_{10,2}\Norm{\hathn^{3/2}\bdisc J^n}_ {\hat\Sigma_n}
        +C_{14,2}\Norm{\hathn^{3/2}\bdisc J^n}_
        {\check\Sigma_n\Tolto\hat\Sigma_n} , \intertext{the {\em data
            approximation error estimators} for {\em space} and {\em time}
          respectively}
        \label{eqn:lil2.inddata.2}
        \inddata_{2,n} :=\frac{C_{3,1}}{\sqrt\alpha} \Norm{\hn(\lproj n
          -\Id)\left(f^n+\frac\uno\taun\right)} ,\, \indf_{1,n}
        :=\frac1\taun\int_\tno^\tn\Norm{f^n-f(t)}\d t \intertext{and the
          {\em time error estimator}} \indtime_{1,n} := \begin{cases}
          \frac12\Norm{\bdisc (\projf n-\discn \un)}\taun &\text{ for
            $n\in\fromto2N$,}\\ \frac12\Norm{\projf 1-\discn U^1-A^0U^0}
          &\text{ for $n=1$}.
        \end{cases}
      \end{gather}
      We refer to Appendix \ref{sec:inequalities} for an explanation of
      the constants $C_{k,j}$ involved here.
    \end{Def}
    \begin{The}[$\leb\infty(\leb2)$ and $\leb2(\sobh1)$ 
        \aposteriori error estimates]
      \label{the:lil2.estimate}
      For each $m\in\fromto1N$ the following error estimates hold
      \begin{gather}
        \label{eqn:lil2.estimate}
        \max_{t\in[0,\tm]}\Norm{u(t)-U(t)} \leq \Norm{\recop 0U^0-u(0)}
        +\max_{n\in\fromto0m}\indrec_{\infty,n} +4
        \big(\cE_{1,m}^2+\cE_{2,m}^2\big)^{1/2} \\
        \label{eqn:l2h1.estimate}
        \begin{split}
          \left(\int_0^\tm\norm{u(t)-U(t)}_1^2\right)^{1/2} \leq&
          \Norm{\recop 0U^0-u(0)} +\left( \sum_{n=1}^m
          \left(\indrec_{2,n}^2+\indrec_{2,n-1}^2\right) \taun
          \right)^{1/2}\\ &+4\big(\cE_{1,m}^2+\cE_{2,m}^2\big)^{1/2}
        \end{split}
      \end{gather}
      \label{def:lil2.estimators}
      where
      \begin{gather}
        \cE_{1,m}:= \sum_{n=1}^m\left(\indtime_{1,n}
        +\indf_{1,n}+\indspace_{1,n}\right)\taun,
        \\ \cE_{2,m}^2:=\sum_{n=1}^m \inddata_{2,n}^2 \taun.
      \end{gather}
    \end{The}
  }%
  \begin{Proof} Following the general strategy of
    \secref{sse:general.strategy} the error is decomposed as follows
    \begin{equation}
      \label{eqn:lil2.error.decomposition}
      \Norm{U(t)-u(t)}=\Norm{e(t)}\leq\Norm{\epsi(t)}+\Norm{\rho(t)}.
    \end{equation}
    To bound the first term, which is the elliptic reconstruction error,
    we apply Lemma \ref{lem:elliptic.aposteriori.estimates}. For the
    estimate in \eqref{eqn:lil2.estimate} we use
    \eqref{eqn:elliptic.estimate.l2} as follows
    \begin{equation*}
      \begin{split}
        \Norm{\epsi(t)} & =\Norm{l_{n-1}(t)\epsi^{n-1}+l_n(t)\epsi^n}
        \leq\maxi{\Norm{\epsi^n}}{\Norm{\epsi^{n-1}}} \\ &
        \leq\max_{n\in\fromto0m}{\Norm{\epsi^n}}
        \leq\max_{n\in\fromto0m}\indrec_{\infty,n},
      \end{split}
    \end{equation*}
    for all $t\in I_n$ and all $n\in\fromto1m$.  In analogous way, we use
    \eqref{eqn:elliptic.estimate.h1} to obtain the estimate of the
    elliptic reconstruction error for \eqref{eqn:l2h1.estimate}.

    The second term on the right-hand side of
    \eqref{eqn:lil2.error.decomposition}, which is the main parabolic
    error, will be estimated via Lemma \ref{lem:lil2.main.error.estimate}
    which we establish next.
  \end{Proof}
  \begin{Lem}[$\leb\infty(\leb2)$
      \aposteriori estimate for the main parabolic error]
    \label{lem:lil2.main.error.estimate}
    For each $m\in\fromto1N$, the following estimate holds
    \begin{equation}
      \label{eqn:lil2.main.error.estimate}
      \begin{split}
        \Big(\max_{[0,t_m]}\Norm{\rho(t)}^2
        +2\int_0^\tm\norm{\rho(t)}_a^2\d t\Big)^{1/2} \leq&\Norm{\rho^0}
        + 4\big(\cE_{1,m}^2+\cE_{2,m}^2\big)^{1/2}.
      \end{split}
    \end{equation}
  \end{Lem}
  We divide the proof of this result in several steps which constitute
  the paragraphs
  \S\S~\ref{sse:lil2.basic.estimate}--\ref{sse:lil2.lemma.proof}.
  \subsection{The basic estimate}
  \label{sse:lil2.basic.estimate}

  To obtain $\leb\infty(\leb2)$ and $\leb2(\sobh1)$ estimates we employ
  standard energy techniques.  We replace $\phi$ in
  \eqref{eqn:main.error.equation} by the main parabolic error
  $\rho=\w-u$ and we integrate in time; thus we have
  \begin{equation}
    \label{eqn:lil2.estimate.base}
    \begin{split}
      \frac12&\Norm{\rho^{m}}^2 -\frac12\Norm{\rho^0}^2
      +\int_0^{t_m}\norm{\rho(t)}_a^2\d t\\ &\leq
      \sum_{n=1}^m\int_{t_{n-1}}^{t_n}
      \norm{\ltwop{\partial_t\epsi(t)}{\rho(t)}} +
      \norm{\abil{\w(t)-\w^n}{\rho(t)}} \\ &+ \norm{\ltwop{ {\lproj
            nf^n-f^n} +\taunm\left({\lproj n\uno-\uno}\right)}{\rho(t)
      }} +\norm{\ltwop{f^n-f(t)}{\rho(t)}} \d t\\ &=:\sum_{n=1}^m
      \left(\cI^1_n+\cI^2_n+\cI^3_n+\cI^4_n\right) =:\cI_m.
    \end{split}
  \end{equation}
  If we denote by $t^m_*\in[0,t_m]$ the time for which
  \begin{equation}
    \max_{t\in[0,t_m]}\Norm{\rho(t)}=\Norm{\rho(t_m^*)}=:\Norm{\rho^m_*}
  \end{equation}
  we deduce that
  \begin{equation}
    \frac12\Norm{\rho^m_*}^2-\frac12\Norm{\rho_0}^2
    +\int_0^{t_m^*}\norm{\rho}_a^2\leq\cI_m.
  \end{equation}
  Consequently we have
  \begin{equation}
    \label{eqn:lil2.estimate.base.star}
    \frac12\Norm{\rho^m_*}^2 +\int_0^{t_m}\norm{\rho}_a^2\leq
    \frac12\Norm{\rho_0}^2+2\cI_m.
  \end{equation}
  We proceed by estimating each of the summands $\cI^i_n$ appearing on
  the right-hand side of \eqref{eqn:lil2.estimate.base}.

  \subsection{Time error estimate}
  \label{sse:lil2.time.estimate}
  In order to bound $\cI_n^2$ in \eqref{eqn:lil2.estimate.base}, which
  accounts for the time discretization error, we use directly the
  elliptic reconstruction definition \eqref{eqn:reconstruction} as
  follows
  \begin{equation*}
    \begin{split}
      \cI_n^2 &= \int_\tno^\tn \norm{\abil{\w(t)-\w^n}{\rho(t)}}\d t
      \\ &=\int_\tno^\tn \norm{\abil{l_{n-1}(t)\recop {n-1}\uno
          +l_n(t)\recop n\un-\recop n\un}{\rho(t)}}\d t \\ &=\int_\tno^\tn
      l_{n-1}(t) \norm{\abil{\recop {n-1}\uno-\recop n\un}{\rho(t)}}\d t
      \\ &=\int_\tno^\tn l_{n-1}(t)
      \norm{\ltwop{\ano\uno-\an\un}{\rho(t)}}\d t .
    \end{split}
  \end{equation*}

  Therefore
  \begin{equation}
    \begin{split}
      \cI_n^2 &\leq \int_\tno^\tn
      l_{n-1}(t)\Norm{\ano\uno-\an\un}\Norm{\rho(t)}\d t,
    \end{split}
  \end{equation}
  which leads to
  \begin{equation}
    \label{eqn:lil2.time.estimate.l1}
    \sum_{n=1}^m\cI_n^2 \leq\Norm{\rho^m_*}
    \sum_{n=1}^m\indtime_{1,n}\taun.
  \end{equation}
  \subsection{Spatial error estimate}
  \label{sse:lil2.space.estimate}
  To estimate the term $\cI_n^1$ on the right-hand side of
  \eqref{eqn:lil2.estimate.base}, which measures the space error and
  mesh change, we will exploit the orthogonality property of the
  elliptic reconstruction (\ref{eqn:reconstruction.orthogonality}).
  Observe that for each $n\in\fromto1N$ we have
  \begin{equation}
    \label{eqn:space.est.on.In}
    \begin{split}
      \cI_n^1&=\int_\tno^\tn
      \norm{\ltwop{\partial_t\epsi(t)}{\rho(t)}}\d
      t\\ &=\taunm\int_\tno^\tn
      \norm{\ltwop{\recop n\un-\recop {n-1}\uno-\un+\uno}{\rho(t)}}\d t.
    \end{split}
  \end{equation}
  Since $\recop n\un-\un$ is orthogonal to $\fez n$ with respect to
  $\abil\cdot\cdot$, the first term inside the brackets is orthogonal to
  $\fez n\cap\fez{n-1}$.  We can therefore use standard residual based
  \aposteriori estimation techniques.  Let
  $\psi:[0,T]\rightarrow\sobhz1(\W)$ be such that
  \begin{equation}
    \abil\chi{\psi(t)}=
    \ltwop{\rho(t)}\chi,\quad\forall\chi\in\honezw,\,t\in[0,T].
  \end{equation}

  \changes{ By \eqref{eqn:reconstruction.orthogonality}, Definition
    \ref{obs:representation} and the use of the interpolation operator
    $\hat\Pi^n$ defined in
    \secref{sse:mesh.change.interpolation.inequalities} }%
  follows that
  \begin{equation}
    \begin{split}
      \langle\recop n\un&-\recop {n-1}\uno-\un+\uno,{\rho(t)}\rangle\\ =&\abil{\recop n\un-\recop {n-1}\uno-\un+\uno}{\psi(t)}\\ =&\abil{\recop n\un-\recop {n-1}\uno-\un+\uno}{{\psi(t)}-\hat\Pi^n{\psi(t)}}\\ =&\ltwop{\an\un-\ano\uno-\pwA
        (\un-\uno)}{{\psi(t)}-\hat\Pi^n{\psi(t)}}\\ &+\ltwop{J^n-J^{n-1}}{{\psi(t)}-\hat\Pi^n{\psi(t)}}_{\check\Sigma_n},
    \end{split}
  \end{equation}
  for each $t\in I_n$.  Using the pointwise form of the fully discrete
  scheme \eqref{eqn:pointwise.form} we can rewrite these terms in a more
  compact form
  \begin{equation}
    \begin{split}
      \an\un-\ano\uno-\pwA (\un-\uno) &=\Big((\discn+\pwA
      )(U^n-U^{n-1})-\projf n+\projf{n-1}\Big)
      \\ \Ignore{&=\taun\taun^{-1}\Big((\discn+\pwA )
        (U^n-U^{n-1})-\projf n+\projf{n-1}\Big)\\}
      &=\taun\left(\bdisc\left((\discn+\pwA )\un-\projf n\right)\right)
      =\taun\bdisc R^n,
    \end{split}
  \end{equation}
  on each interval $I_n$.  Here $R^n:=(\discn+\pwA )U^n-\projf n$ is the
  {\em internal residual} function at time $t_n$.  Likewise we have
  \begin{equation}
    J^n-J^{n-1}=\taun\bdisc J^n.
  \end{equation}
  Whence, with $j\in\Nat$ being at our disposal, in view of the
  interpolation inequalities in \secref{sse:scott.zhang.inequalities} we
  may conclude that
  \begin{equation}
    \label{eqn:lil2.space.estimate.base}
    \begin{split}
      \cI_n^1 \leq &
      \bigg(\int_\tno^\tn\norm{\psi(t)}_j\bigg)\times\\ &\bigg(
      C_{3,j}\Norm{\hathn^j \bdisc R^n} +C_{5 ,j}\Norm{\hathn^{j-1/2}
        \bdisc J^n}_ {\hat\Sigma_n} +C_{7 ,j}\Norm{\hathn^{j-1/2} \bdisc
        J^n}_ {\check\Sigma_n\Tolto\hat\Sigma_n} \bigg).
    \end{split}
  \end{equation}
  \changes{ Since $\ell\geq1$, we may take $j=2$ in
    \eqref{eqn:lil2.space.estimate.base} and use the elliptic regularity
    \eqref{eqn:elliptic.regularity} with $k=0$, to get
    \begin{equation}
      \cI_n^1\leq\max_{t\in I_n}\Norm{\rho(t)}\indspace_{1,n}\taun
    \end{equation}
    where $\indspace_{1,n}$ is given by \eqref{eqn:lil2.indspace.1}.
    Hence
    \begin{equation}
      \label{eqn:lil2.space.estimate.l1}
      \begin{split}
        \sum_{n=1}^m\cI_n^1
        \leq\Norm{\rho^m_*}\sum_{n=1}^m\indspace_{1,n}\taun.
      \end{split}
    \end{equation}
  }%
  \subsection{Data approximation and mesh change estimates}
  \label{sse:lil2.data.estimate}
  We now bound the term $\cI_n^3$ in \eqref{eqn:lil2.estimate.base}.
  Here we exploit the orthogonality of the $\leb2$-projection.  Since
  $\fez n\subset\ker(\lproj n - \Id)$ we have
  \begin{equation}
    \begin{split}
      \cI_n^3 &=\int_\tno^\tn \norm{ \ltwop{(\lproj n
          -\Id)(f^n+\taunm\uno)} {\rho(t)-\Pi^n\rho(t)}}\d
      t\\ &\leq\int_\tno^\tn \norm{\ltwop{\hn(\lproj n
          -\Id)(f^n+\taunm\uno)} {\hn^{-1}(\rho(t)-\Pi^n\rho(t))}}\d
      t\\ &\leq\alpha^{-1/2}C_{3,1} \taun^{1/2}\Norm{\hn(\lproj n
        -\Id)(f^n+\taunm\uno)} \left(\int_\tno^\tn\norm{\rho(t)}_a\d
      t\right)^{1/2}.
    \end{split}
  \end{equation}
  We can thus conclude that
  \begin{equation}
    \label{eqn:lil2.data.estimate.l2}
    \sum_{n=1}^m \cI_n^3 =\sum_{n=1}^m
    \left(\int_\tno^\tn\norm{\rho}_a^2\right)^{1/2}\inddata_{2,n}\taun^{1/2}.
  \end{equation}

  \label{sse:lil2.f.estimate}
  We conclude this paragraph by estimating the fourth term in the
  right-hand side of \eqref{eqn:lil2.estimate.base} in a simple way as
  follows
  \begin{equation}
    \cI_n^4 \leq \left(\max_{t\in I_n}
    \Norm{\rho(t)}\right)\int_\tno^\tn\Norm{f^n-f(t)}\d t
  \end{equation}
  thus
  \begin{equation}
    \label{eqn:lil2.f.estimate.l1}
    \sum_{n=1}^m\cI_4\leq\Norm{\rho^m_*}\sum_{n=1}^m\indf_{1,n}\taun.
  \end{equation}
  \subsection{\changes{Proof of Lemma 
      \ref{lem:lil2.main.error.estimate}: last step}}
  \label{sse:lil2.lemma.proof}
  What remains to do in order to conclude the proof is to combine
  appropriately the results from
  \S\S\ref{sse:lil2.time.estimate}--\ref{sse:lil2.f.estimate} with
  inequalities \eqref{eqn:lil2.estimate.base} and
  \eqref{eqn:lil2.estimate.base.star}.  We can write
  \begin{equation}
    \label{eqn:lil2.prelemma}
    \begin{split}
      \frac12\Norm{\rho^m_*}^2 +\int_0^{t_m}\norm{\rho}_a^2\leq
      &\frac12\Norm{\rho_0}^2 +2\Norm{\rho^m_*}\sum_{n=1}^m
      \left(\indtime_{1,n}+\indf_{1,n}+\indspace_{1,n}\right)\tau_n\\ &+2\sum_{n=1}^m\left(\int_\tno^\tn\norm{\rho}_a^2\right)^{1/2}
      \inddata_{2,n}\taun^{1/2}.
    \end{split}
  \end{equation}
  We can apply now the elementary fact that, for $\vec
  a=(a_0,\ldots,a_m),\,\vec b=(b_0,\ldots,b_m)\in\rR^{m+1}$ and
  $c\in\rR$, if
  \begin{equation}
    \label{eqn:swallowed.sphere.hyp}
    \norm{\vec a}^2\leq c^2 + \vec a\inner\vec b
  \end{equation}
  then
  \begin{equation}
    \norm{\vec a}\leq \norm{c}+\norm{\vec b}.
  \end{equation}
  In particular, in reference to \eqref{eqn:lil2.prelemma}, we take
  \begin{gather}
    a_0=\frac1{\sqrt 2}\Norm{\rho^m_*},\quad
    a_n=\left(\int_\tno^\tn\norm{\rho}_a^2\right)^{1/2},\quad
    c=\frac1{\sqrt 2}\Norm{\rho^0},\\ b_0=2\sqrt
    2\sum_{n=1}^m\left(\indtime_{1,n}+\indf_{1,n}+\indspace_{1,n}\right)\tau_n,\quad
    b_n=2\inddata_{2,n}\taun^{1/2},
  \end{gather}
  for $n\in\fromto1m$, obtain \eqref{eqn:lil2.main.error.estimate},
  which concludes the proof of the lemma.\hfill\qed

  \begin{Obs}[\changes{relation to the semidiscrete case}]
    The spatial error estimators, containing $\estspace{}$'s, should be
    compared with the ones corresponding to the (space) semi-discrete
    scheme given by
    \begin{equation}
      \left(\Norm{h^3\partial_t((\partial_t-\pwA )u_h-f)}^2
      +\Norm{h^{5/2}\partial_tJ[u_h]}_{\Sigma_h}^2 \right)^{1/2}
    \end{equation}
    where one triangulation $\cT_h$ is given for all (continuous) time
    $t\in[0,T]$ \citep[eqn. (4.4)]{MakridakisNochetto:03:article:Elliptic}.
  \end{Obs}

  \begin{Obs}[mesh change]
    We interpret the presence of the coarsest meshsize $\hathn$ in the
    estimator as a reflection of the discrepancy between the finite
    element spaces $\fez{n-1}$ and $\fez n$, which might be different in
    general.  Mesh change, a delicate issue in evolution problems, can
    lead to non-convergent schemes despite the global meshsize going to
    zero.  This happens in an example by \citet[\secsymbol
      4]{Dupont:82:article:Mesh} for which
    \begin{equation}
      \max_n{\sup_\W h_n}\rightarrow0\quad \left(\text{ but }
      \max_n{\sup_\W\hat h_n}\not\rightarrow0 \right)
    \end{equation}
    yet the discrete solution does not converge to the exact solution.

  \end{Obs}

  \section{\Aposteriori error estimates in higher order norms}
  \label{sec:lih1}
  In this section we give estimates in the seminorms corresponding to
  the spaces $\sobh1(0,T;\leb2(\W))$ and $\leb\infty(0,T;\sobh1(\W))$.
  The exposition of this section parallels that of \secref{sec:lil2}.
  We start by introducing \aposteriori error estimators that are local
  in time and which will be used in the subsequent main result.  Notice
  that the error estimators in the next definition, compared to those of
  \secref{sec:lil2}, are redefined.
  \begin{Def}[error estimators for the $\leb\infty(\sobh1)$%
      and $\sobh1(\leb2)$ semi-norms] \changes{ We introduce the {\em
        elliptic reconstruction error estimator}
      \begin{align}
        \label{eqn:lih1.indrec}
        \indrec_{\infty,n} &:=\frac{C_{3,1}}\alpha \Norm{\hn R^n}
        +\frac{C_{5,1}}\alpha \Norm{\hn^{1/2} J^n} , \\ \indrec_{2,n}
        &:=C_{6,2}\Norm{\hathn^2\bdisc R^n}
        +C_{10,2}\Norm{\hathn^{3/2}\bdisc J^n}_{\hat\Sigma_n}
        +C_{14,2}\Norm{\hathn^{3/2}\bdisc
          J^n}_{\check\Sigma_n\Tolto\hat\Sigma_n}
        \intertext{the {\em space error estimator}}
        \label{eqn:lih1.indspace.2}
        \indspace_{2,n} &:= C_{6,2}\Norm{\hathn^2\bdisc R^n}
        +C_{10,2}\Norm{\hathn^{3/2}\bdisc J^n}_{\hat\Sigma_n}
        +C_{14,2}\Norm{\hathn^{3/2}\bdisc
          J^n}_{\check\Sigma_n\Tolto\hat\Sigma_n} ,
        \intertext{the {\em data space approximation error estimators}}
        \label{eqn:lih1.inddata.1}
        \inddata_{1,n} & :=\alpha^{-1/2}C_{3,1}
        \Norm{\hathn\bdisc\left((\lproj n-\Id)(f^n-\taun\uno)\right)} ,
        \\
        \label{eqn:lih1.inddata.infty}
        \inddata_{\infty,n} & :=\alpha^{-1/2}C_{3,1} \Norm{\hn(\lproj
          n-\Id)(f^n-\taun\uno)} ,
        \intertext{the {\em data time approximation error estimator} and
          the {\em time error estimator}}
        \label{eqn:lih1.indf.2}
        \indf_{2,n} &
        :=\left(\frac1\taun\int_\tno^\tn\Norm{f^n-f(t)}^2\d
        t\right)^{1/2},
        \quad \indtime_{2,n}
        =\frac1{\sqrt3}\Norm{\partial_t(\projf{}-\discn U)}\taun .
      \end{align}
    }%
    definition of constants $C_{k,j}$ involved above.
  \end{Def}
  \changes{
    \begin{The}[$\leb\infty(\sobh1)\cap\sobh1(\leb2)$ 
        \aposteriori estimates]
      \label{the:lih1.estimate}
      Suppose the exact solution $u$ of \eqref{eqn:continuous.heat}
      satisfies
      \begin{gather}
        \label{eqn:lih1.regularity.1}
        \partial_tu\in\leb2(0,T;\leb2(\W)),\\
        \label{eqn:lih1.regularity.2}
        \partial_tu(t)\in\honezw\text{ for a.e. }t\in[0,T],\\
        \label{eqn:lih1.regularity.3}
        \grad u\in\leb2(0,T;\honew).
      \end{gather}
      Then the following \aposteriori error estimates hold
      \begin{gather}
        \begin{split}
          \left(\int_0^\tm \Norm{\partial_t\left(U(t)-u(t)\right)}^2\d t
          \right)^{1/2} \leq& \norm{\recop 0U^0-u(0)}_a\\ &+4\big(
          \cE_{1,m}^2+\cE_{2,m}^2 \big)^{1/2} +\estrec m'
        \end{split}
        \\ \max_{t\in[0,\tm]}\norm{U(t)-u(t)}_a \leq \norm{\recop 0U^0-u(0)}_a
        +4\big( \cE_{1,m}^2+\cE_{2,m}^2 \big)^{1/2} +\estrec m,
      \end{gather}
      where
      \begin{gather*}
        \cE_{1,m} :=2\max_{n\in\fromto1m}
        \inddata_{\infty,m}+\sum_{n=2}^m\inddata_{1,n}\taun , \\ \cE_{2,m}^2
        := \sum_{n=1}^m\left(
        \indtime_{2,n}^2+\indf_{2,n}^2+\indspace_{2,n}^2 \right)\taun,
        \\ \estrec m := \max_{n\in\fromto0m}\indrec_{\infty,n} \quad\text{
          and }\quad (\estrec m')^2 := \sum_{n=1}^m\indrec_{2,n}^2.
      \end{gather*}
    \end{The}
  }%
  \begin{Proof} 
    Following the general strategy of \secref{sse:general.strategy}
    the error is decomposed as follows
    \begin{equation}
      \label{eqn:lih1.error.decomposition}
      \norm{U(t)-u(t)}_Y=\norm{e(t)}_Y\leq\norm{\epsi(t)}_Y+\norm{\rho(t)}_Y,
    \end{equation}
    where $\norm{\cdot}_Y$ is the appropriate seminorm\footnote{%
    this section we deliberately use $\norm{\cdot}_a$ instead of
    $\norm{\cdot}_1$ as the norm for $\honezw$ in order to keep the
    exposition clear. The changes to replace $\norm{\cdot}_a$ by
    $\norm{\cdot}_1$ are straightforward.  } with $Y$ being either one
    of $\sobh1(0,t_m;\leb2(\W))$ or $\leb\infty(0,t_m;\honezw)$.  The
    first term on the left-hand side of
    \eqref{eqn:lih1.error.decomposition}, which is the elliptic
    reconstruction error, can be estimated in the residual based
    context via Lemma \ref{lem:elliptic.aposteriori.estimates}, as to
    obtain
    \begin{equation}
      \norm{\epsi(t)}_Y\leq
      \begin{cases}
        \estrec m',&\text{ if $Y=\sobh1(0,t_m;\leb2(\W))$},\\ \estrec
        m,&\text{ if $Y=\leb\infty(0,t_m;\honezw)$}.
      \end{cases}
    \end{equation}
    The second term on the left-hand side of
    \eqref{eqn:lih1.error.decomposition} is estimated with the help of
    Lemma \ref{lem:lih1.main.error.estimate} which we state and prove
    next.
  \end{Proof}
  \begin{Lem}[$\leb\infty(\sobh1)\cap\sobh1(\leb2)$ estimates for the 
      main parabolic error]
    \label{lem:lih1.main.error.estimate}
    For each $m\in\fromto1N$ the following \aposteriori estimate is
    valid
    \begin{equation}
      \label{eqn:lih1.main.error.estimate}
      \left(\max_{t\in[0,\tm]}\norm{\rho(t)}_a^2
      +2\int_0^\tm\Norm{\partial_t\rho}^2\right)^{1/2} \leq
      \norm{\rho^0}_a+4\left(\cE_{1,m}^2+\cE_{2,m}^2\right)^{1/2},
    \end{equation}
    with reference to the notation of Theorem \ref{the:lih1.estimate}.
  \end{Lem}
  As in \secref{sec:lil2}, the proof of this lemma is subdivided into
  several steps that constitute paragraphs
  \S\S~\ref{sse:lih1.basic.estimate}--\ref{sse:lih1.lemma.proof}.  Since
  the arguments are very similar to those of \secref{sec:lil2}, the we
  condense the discourse and stress only the main differences.  To
  motivate the proof, we discuss first the semidiscrete case.

  \subsection{The spatially discrete case}
  \label{sse:lih1.semidiscrete}
  The estimates of Theorem \ref{the:lih1.estimate} are based on the
  energy estimate in higher order norms for problem
  \eqref{eqn:continuous.heat} which reads
  \begin{equation}
    \left(\norm{u(t)}_a^2+2\int_0^t\Norm{\partial_t u}^2\right)^{1/2}
    \leq\norm{u(0)}_a+\sqrt2\Norm{f}_{\leb2(0,t;\leb2(\W))}.
  \end{equation}
  This estimate can be obtained by testing the PDE with $\partial_tu$
  and integrating in time.  For this more regularity of $u$ is required
  in this case than in the lower order norms case of \secref{sec:lil2}.
  Sufficient regularity requirements on $u$ are those of conditions
  \eqref{eqn:lih1.regularity.1}--\eqref{eqn:lih1.regularity.3}.

  We stress that in some particular situations these stronger energy
  norms can play an important role.  For instance, an estimate that is
  based on nonlinear quantities similar to these higher order norms has
  been derived for the error related to the mean curvature flow of
  graphs \cite{LakkisNochetto:05:article:A-posteriori}; in that situation there is no
  reasonable way to obtain estimates by testing with the solution. The
  only approach that works is by testing with the time derivative of the
  solution.

  Let us turn now our discussion toward the use of this energy estimate
  in the semidiscrete case; namely only spatially discrete, with
  $\fes{h}$ as a finite element space.  The semidiscrete case is
  simpler, and motivated by, the more involved fully discrete case which
  we will deal with in the next paragraphs.  The semidiscrete case has
  been extensively studied by Makridakis \& Nochetto for the usual
  (lower order) norms \cite{MakridakisNochetto:03:article:Elliptic}.  We further
  simplify our discussion by assuming also that $f\in\fez{h}$.

  The starting point of the error estimate is, like in
  \secref{sec:lil2}, the semidiscrete analog of
  \eqref{eqn:main.error.equation} which is given by \changes{
    \begin{equation}
      \duality{\partial_t\rho}\phi+\abil\rho\phi
      =\ltwop{\partial_t\epsi}\phi.
    \end{equation}
  }%
  \cite[Eqn. (3.2)]{MakridakisNochetto:03:article:Elliptic}.)  Taking this identity with
  $\phi=\partial_t\rho$%
  extra regularity properties
  \eqref{eqn:lih1.regularity.1}--\eqref{eqn:lih1.regularity.3}---}%
and integrating in time we obtain
\begin{equation}
  \begin{split}
    \int_0^T\Norm{\partial_t\rho}^2
    &+\frac12\norm{\rho(T)}_a^2-\frac12\norm{\rho(0)}_a^2
    =\int_0^T\ltwop{\partial_t\epsi}{\partial_t\rho}.
  \end{split}
\end{equation}
We have now the choice to control the right-hand side in two different
ways.
\begin{LetterList}
\item
  We use a straightforward \changes{$\leb2(0,T;\leb2(\W))$} estimate
  that leads to
  \begin{equation}
    \begin{split}
      \int_0^T\ltwop{\partial_t\epsi}{\partial_t\rho}
      &=\left(\int_0^T\Norm{\partial_t\epsi}^2\right)^{1/2}
      \left(\int_0^T\Norm{\partial_t\rho}^2\right)^{1/2}
      \\ &=\left(\int_0^T\cE[\fes{h},A_h\partial_tu_h;\leb2(\W)]^2\right)^{1/2}
      \left(\int_0^T\Norm{\partial_t\rho}^2\right)^{1/2}
    \end{split}
  \end{equation}
  where $\cE[\cdot]$ is an {\em elliptic error estimator function}
  \cite{MakridakisNochetto:03:article:Elliptic}.  This could be for instance, but
  not necessarily so, the residual based estimators of Lemma
  \ref{lem:elliptic.aposteriori.estimates}.
\item
  An alternative estimate that is often useful when quadratic or
  higher finite elements are employed or when only the energy norm
  $\norm{\rho}_a$ can be controlled---in a nonlinear setting, for
  instance---involves an integration by parts in time as follows:
  \begin{equation}
    \label{eqn:semidiscrete.d2}
    \begin{split}
      &\int_0^T\!\!\!\!\!\!\ltwop{\partial_t\epsi}{\partial_t\rho}
      -\ltwop{\partial_t\epsi(0)}{\rho(0)}
      =\ltwop{\partial_t\epsi(T)}{\rho(T)}
      -\int_0^T\!\!\!\!\!\!\ltwop{\partial_t^2\epsi}{\rho}\\ &\leq
      C\max_{[0,T]}\norm{\rho}_a \!\left(
      \max_{[0,T]}\cE[\fes{h},A_h\partial_tu_h;\sobh{-1}(\W)]
      +\int_0^T\!\!\!\!\!\!\cE[\fes{h},A_h\partial_t^2u_h;\sobh{-1}(\W)]
      \right).
    \end{split}
  \end{equation}
  The extra time differentiation will not affect the order of
  convergence of the right-hand side because the elliptic
  reconstruction error $\ep$ is purely elliptic in nature.  We
  note, however, that in order to make estimate
  \eqref{eqn:semidiscrete.d2} rigorous, it is necessary to impose
  extra time-regularity assumptions on the approximate solution
  and to have an $\sobh{-1}(\W)$-norm elliptic error estimator
  function $\cE[\cdot,\cdot;\sobh{-1}(\W)]$ available.  Such
  estimators can be obtained with optimal order, using the duality
  techinque, under the assumption that the domain is smooth.  This
  is an issue that is beyond the scope of this paper so we will
  limit our analysis to the first alternative.
\end{LetterList}

\subsection{The basic estimate for the fully discrete case}
\label{sse:lih1.basic.estimate}
We proceed now with the proof of Lemma
\ref{lem:lih1.main.error.estimate}. The first step consists of taking
$\phi=\partial_t\rho$ in identity \eqref{eqn:main.error.equation} and
integrate by parts in time as follows
\begin{equation}
  \label{eqn:lih1.estimate.base}
  \begin{split}
    \int_0^\tm&\Norm{\partial_t\rho}^2 +\frac12\norm{\rho^m}_a^2
    =\frac12\norm{\rho^0}_a^2\\ &+\sum_{n=1}^m \int_\tno^\tn
    \ltwop{\partial_t\epsi}{\partial_t\rho(t)}
    +\abil{\w(t)-\w^n}{\partial_t\rho(t)}\\ &+\ltwop{\lproj n
      f^n-f^n+\taunm(\lproj n\uno-\uno)}{\partial_t\rho(t)}
    \\ &+\ltwop{f^n-f(t)}{\partial_t\rho(t)}\d
    t\\ &=:\frac12\norm{\rho(0)}_a^2+\sum_{n=1}^m(\cI_n^1+\cI_n^2+\cI_n^3+\cI_n^4)
    =:\cI_m.
  \end{split}
\end{equation}
Introduce $\tms\in[0,\tm]$ such that
\begin{equation}
  \norm{\rhoms}_a:=\norm{\rho(\tms)}_a=\max_{t\in[0,\tm]}\norm{\rho(t)}_a.
\end{equation}
Let $\ms\in\fromto1m$ be the index for which $\tms\in I_\ms$. We can
write
\begin{equation}
  \label{eqn:lih1.estimate.base.star}
  \begin{split}
    \int_0^\tms\Norm{\partial_t\rho}^2 +\frac12\norm{\rhoms}_a^2
    =&\frac12\norm{\rho(0)}_a^2
    +\sum_{n=1}^\ms\left(\cJ_n^1+\cJ_n^2+\cJ_n^3+\cJ_n^4\right)
  \end{split}
\end{equation}
where
\begin{equation}
  \cJ_n^i =
  \begin{cases}
    \cI_n^i &\text{ for }n\in\fromto1{\ms-1} , \\ \int_{t_\mso}^\tms
    I_n^i(t)\d t &\text{ for }n=\ms .
  \end{cases}
\end{equation}
with $I_n^i$ being the same integrand as that of $\cI_n^i$.  To prove
the lemma we must bound $\sum_{n=1}^m\cI^k_n$ and
$\sum_{n=1}^\ms\cJ^k_n$, for each $k\in\fromto14$.

\subsection{Time error estimate}
\label{sse:lih1.time.estimate}
We estimate the term due to time discretization.  By the definition of
elliptic reconstruction \eqref{eqn:reconstruction} we have
\begin{equation*}
  \begin{split}
    \cI_n^2 &\leq \int_\tno^\tn
    \norm{\abil{\w(t)-\w^n}{\partial_t\rho(t)}}\d t \\ &=\int_\tno^\tn
    l_{n-1}(t) \norm{\ltwop{\ano\uno-\an\un}{\partial_t\rho(t)}}\d
    t\\ &\leq\left(\int_\tno^\tn\Norm{\partial_t\rho(t)}^2\d
    t\right)^{1/2}\sqrt {1/3}\Norm{\ano\uno-\an\un}\taun^{1/2}\\
  \end{split}
\end{equation*}
It follows that
\begin{equation}
  \sum_{n=1}^m\cI_n^2\leq
  \sum_{n=1}^m\left(\int_\tno^\tn\Norm{\partial_t\rho}^2\right)^{1/2}
  \indtime_{2,n}\taun^{1/2}.
\end{equation}
The same bound applies to $\sum_{n=1}^\ms\cJ_n^2$.

\subsection{Spatial error estimate} 
\label{sse:lih1.space.estimate}
The spatial error estimator term can be bounded in a similar way to
the one in \secref{sse:lil2.space.estimate}.  Introduce first the
auxiliary function \changes{$\psi:[0,T]\rightarrow\honezw$} such that
\changes{
  \begin{equation}
    \abil\chi{\psi(t)}=\ltwop{\partial_t\rho(t)}\chi,
    \:\forall\chi\in\honezw.
  \end{equation}
}%

\changes{ Noting that $\partial_t\epsi$ is a piecewise constant
  function of time, and in view of
  \eqref{eqn:reconstruction.orthogonality},
  \eqref{eqn:scott.zhang.inner}--\eqref{eqn:scott.zhang.coral} (with
  $j=2$) and \eqref{eqn:elliptic.regularity} (with $k=0$), we can
  write
  \begin{equation}
    \begin{split}
      \cI_n^1 &=\int_\tno^\tn\abil{\partial_t\epsi^n}{\psi(t)}\d t
      =\int_\tno^\tn\abil{\partial_t\epsi^n}{\psi(t)-\hat\Pi^n\psi(t)}
      \d t \\ &=\int_\tno^\tn\ltwop{\bdisc
        R^n}{\psi(t)-\hat\Pi^n\psi(t)} +\ltwop{\bdisc
        J^n}{\psi(t)-\hat\Pi^n\psi(t)}_{\check\Sigma_n} \d t\\ &\leq
      \left(\int_\tno^\tn\Norm{\partial_t\rho(t)}^2\right)^{1/2}
      \indspace_{2,n}\taun^{1/2}
    \end{split}
  \end{equation}
}%
\eqref{eqn:lih1.indspace.2}.  Thus, upon summing in time, we conclude
that
\begin{equation}
  \sum_{n=1}^m\cI_n^1\leq\sum_{n=1}^m
  \left(\int_\tno^\tn\Norm{\partial_t\rho(t)}^2\right)^{1/2}
  \indspace_{2,n}\taun^{1/2}.
\end{equation}
The same estimate holds for $\sum_{n=1}^\ms\cJ_n^1$.

\subsection{Data approximation and mesh change estimates}
\label{sse:lih1.data.estimate}
We conclude the estimates in this section by bounding the last two
terms in \eqref{eqn:lih1.estimate.base} regarding data approximation
and mesh changes.

The data space approximation error can be bounded as follows
\begin{equation}
  \begin{split}
    \sum_{n=1}^m\cI_n^3 =&\sum_{n=1}^m\int_\tno^\tn\ltwop{(\lproj
      n-\Id)(f^n-\taun\uno)}{\partial_t\rho}\\ =&\sum_{n=1}^m\ltwop{(\lproj
      n-\Id)(f^n-\taun\uno)}{\rho^n-\rho^{n-1}}\\ =&\sum_{n=1}^{m-1}
    \ltwop{(\lproj n-\Id)(f^n-\taun\uno) -(\lproj
      {n+1}-\Id)(f^{n+1}-\tau_{n+1}\un)} {\rho^n}\\ &+\ltwop{(\lproj
      m-\Id)(f^m-\tau_m U^{m-1})}{\rho^m} -\ltwop{(\lproj
      1-\Id)(f^1-\tau_1 U^0)}{\rho^0} \\ =&\sum_{n=2}^{m}
    \ltwop{\bdisc\left((\lproj
      n-\Id)(f^n-\taun\uno)\right)}{\rho^{n-1}-\hat\Pi^n\rho^{n-1}}\taun
    \\ &+\ltwop{(\lproj m-\Id)(f^m-\tau_m
      U^{m-1})}{\rho^m-\Pi^m\rho^m}\\ &-\ltwop{(\lproj
      1-\Id)(f^1-\tau_1 U^0)}{\rho^0-\Pi^0\rho^0}.
  \end{split}
\end{equation}
Owing to \eqref{eqn:scott.zhang.inner} and
\eqref{eqn:one.energy.inequality} we may thus conclude that
\begin{equation}
  \sum_{n=1}^m\cI_n^3\leq\norm{\rho^m_*}_a
  \left(\inddata_{\infty,m}+\sum_{n=2}^m\inddata_{1,n}\taun+\inddata_{\infty,1}\right)
\end{equation}
where the estimators $\gamma_{i,n}$ are defined in
\eqref{eqn:lih1.inddata.1} and \eqref{eqn:lih1.inddata.infty} for
$i=1$ and $\infty$ respectively.  Likewise we obtain the following
bound
\begin{equation}
  \sum_{n=1}^\ms\cJ_n^3\leq\norm{\rho^m_*}_a
  \left(\inddata_{\infty,\ms}+\sum_{n=2}^\ms\inddata_{1,n}\taun
  +\inddata_{\infty,1}\right).
\end{equation}
The last term in \eqref{eqn:lih1.estimate.base} is handled in a
straightforward way, as in \secref{sse:lil2.data.estimate} and the
following bound is readily derived
\begin{equation}
  \sum_{n=1}^m\cI_n^4 \leq\sum_{n=1}^m
  \left(\int_\tno^\tn\Norm{\partial_t\rho}^2\right)^{1/2}
  \indf_{2,n}\taun^{1/2}.
\end{equation}
Also here the same bound applies to $\sum_{n=1}^\ms\cJ_n^4$.

\subsection{Proof of Lemma \ref{lem:lih1.main.error.estimate}: last step}
\label{sse:lih1.lemma.proof}
As in \ref{sse:lil2.lemma.proof}, we collect appropriately the results
from the preceding paragraphs and we use
\eqref{eqn:lih1.estimate.base} and
\eqref{eqn:lih1.estimate.base.star}.  We thus have
\begin{equation}
  \begin{split}
    \frac12\norm{\rhoms}_a^2 +\int_0^{t_m}\Norm{\partial_t\rho}^2\leq
    \frac12\norm{\rho_0}_a^2 +2\norm{\rhoms}_a \cE_{1,m}
    \\ +2\sum_{n=1}^m\left(\int_\tno^\tn\Norm{\partial_t\rho}^2\right)^{1/2}
    \left(\indtime_{2,n}+\indf_{2,n}+\indspace_{2,n}\right)\taun^{1/2}
    .
  \end{split}
\end{equation}
We can proceed now by using the same elementary fact used in
\secref{sse:lil2.lemma.proof}, with
\begin{gather}
  a_0=\frac1{\sqrt 2}\norm{\rhoms}_a, \quad
  a_n=\left(\int_\tno^\tn\Norm{\partial_t\rho}^2\right)^{1/2}, \quad
  c=\frac1{\sqrt 2}\norm{\rho^0}_a,\\ b_0=2\sqrt 2\cE_{1,m},\quad
  b_n=2\left(\indtime_{2,n}+\indf_{2,n}+\indspace_{2,n}\right)\taun^{1/2},
\end{gather}
for $n\in\fromto1m$.  Some simple manipulations yield estimate
\eqref{eqn:lih1.main.error.estimate}.\hfill\qed
\subsection{The difficulty with the direct approach}
\label{sse:direct.approach.failure}
We conclude this section by exhibiting the main problem with the
direct approach to derive energy estimates in higher norms.

\changes{To see this we go back to equation
  \eqref{eqn:direct.approach.orthogonality} which we now take with
  $\phi=\partial_t e$ and integrate in time.  We are thus required to
  estimate the term
  \begin{equation}
    \sum_{n=1}^m\int_\tno^\tn\abil \un{\partial_t e-\Pi_n\partial_t e}.
  \end{equation}
  Since the only practical way to proceed seems to be by decreasing the
  number of derivatives acting upon the error term---integration by
  parts in space being of no help here---we perform a ``summation by
  parts'' in time as follows:
  \begin{equation*}
    \begin{split}
      \sum_{n=1}^m&\left( \abil\un{e^n-\Pi^ne^n}
      -\abil\un{e^{n-1}-\Pi^ne^{n-1}}\right)\\ &=\abil\un{e^m-\Pi^me^m}-\abil{U^0}{e^0-\Pi^1e^0}
      +\sum_{n=1}^{m-1}\abil{\un-U^{n+1}}{e^n-\Pi^ne^n}\\ &-\sum_{n=1}^{m-1}\abil{\un}{(\Pi^n-\Pi^{n+1})e^n}.
    \end{split}
  \end{equation*}
  The difficulty, which should be apparent now, is how to control the
  last term.  There seems to be no practical way to do this without
  imposing strong assumptions on $(\Pi^n-\Pi^{n+1})e^n$.  Notice that
  this term vanishes if there is no mesh change.  }%
\changes{
  \section{Numerical results}
  \label{sec:numerics}
  We present the results of a series of numerical experiments to
  exemplify some of the practical aspects of the \aposteriori estimates
  of Theorem \ref{the:lil2.estimate}.  The main goal here is to
  approximate the asymptotic behavior of the various estimators and
  compare this behavior with that of the norms.
  
  \subsection{Benchmark solutions}
  We perform the numerical experiment by approximating in each case,
  either one of the following two exact solutions:
  \begin{equation}
    \label{eqn:benchmark}
    u(\vec x,t) =
    \begin{cases}
      \sin(\pi t)\exp(-10\norm{\vec x}^2),&\text{(slow)}
      \\ \frac1{10}\sin(20\pi t)\exp(-10\norm{\vec x}^2),&\text{(fast)}.
    \end{cases}
  \end{equation}
  These solutions are used as benchmarks for the problem with $\vec A =
  \left[\begin{smallmatrix}1&0\\0&1\end{smallmatrix}\right]$. The
  right-hand side $f$ of the problem, is thus easily calculated by
  applying the parabolic operator to each $u$.  Therefore the exact
  errors are computable and we can compare them with the error
  estimators.  The domain on which we compute this solution is the
  square $[-1,1]^2$ and the time interval is $[0,1]$.  Notice that the
  boundary conditions are not exactly zero but of the order of $10^{-6}$
  so that special care has to be taken with very small numbers.  The
  initial conditions are exactly zero in both cases and there is no
  initial error to be computed.
  
  \subsection{Choice of parameters}
  Since we are interested in understanding the asymptotic behavior of
  the estimators, we conduct tests on uniform meshes with uniform
  timestep.  For each numerical experiment we choose a sequence of
  meshsizes $(h(i):i\in\fromto1I)$, to which we couple a sequence
  stepsizes $(\tau(i):i\in\fromto1I)$, $\tau(i)=c_0h(i)^k$, with $k$
  equal either $1$ or $2$ and $I$, the number of runs, either $5$ or $6$
  in each of the 4 cases that we run.  The following table summarizes
  the choice of the various parameters for each for each of the 4
  simulations:
  \begin{center}
    \begin{tabular}{cccccccc}
      \hline simulation&problem&$\ell$ &$k$ &$h_1$ & $\tau_1$ & $I$
      (runs) & figure \\\hline\hline 1 & slow & 1 & 2 & 0.5 & 0.04 & 6 &
      \ref{fig:P1-tauh2} \\\hline 2 & fast & 1 & 1 & 0.25 & 0.01 & 5 &
      \ref{fig:P1-fast-tauh1} \\\hline 3 & slow & 1 & 3 & 0.125 & 0.08 &
      4 & \ref{fig:P2-tauh3} \\\hline 4 & fast & 2 & 2 & 0.125 & 0.02 &
      4 & \ref{fig:P2-fast-tauh2} \\\hline
    \end{tabular}
  \end{center}

  \subsection{Computed quantities}
  Notice that we report numerical results for the estimates of
  \secref{sec:lil2} only. For each simulation and for each run
  $i\in\fromto1I$, we calculate the following quantities
  \begin{itemize}
  \item  
    the error norms $\Norm{e}_{\leb\infty(0,t_m;\leb2(\W))}$ and
    $\norm{e}_{\leb2(0,t_m;\sobh1(\W))}$,
  \item the reconstruction error estimators $\max_0^m\indrec_{\infty,n}$
    and $(\tau(i)\sum_1^m\indrec_{2,n}^2)^{1/2}$,
  \item the space estimator $\tau(i)\sum_1^m\indspace_{1,n}$,
  \item and the time estimator $\tau(i)\sum_1^m\indtime_{1,n}$,
  \end{itemize}
  for each time $t_m\in\fromstepto{0=t_0}{\tau(i)}{t_N=1}$.  Of course,
  all the errors and estimators depend on the run $i$, but for sake of
  conciseness, we do not add this index.  We deliberately ignore the
  estimators for data approximation (and mesh change), $\indf_{1,n}$ and
  $\inddata_{2,n}$, as our examples are designed so to make these
  estimators either negligible or comparable with respect one of those
  we study.  Indeed, we take the mesh and step sizes small enough as to
  resolve the data and we keep the mesh unchanged accross timesteps; so
  $\inddata_{2,n}$ can be shown to be of order $h^2$ and $\indf_{1,n}$
  to be of order $\tau$ since the function $f$ is smooth enough
  \cite{LakkisMakridakisPryer:14:article:A-comparison}.

  For each computed norm or estimator we look at its \emph{experimental
  order of convergence (EOC)}. The EOC is defined as follows: for a
  given finite sequence of uniform triangulations
  $\{\cT_{h(i)}\}_{i=1,\ldots,I}$ of meshsize $h(i)$, the EOC of a
  corresponding sequence of some triangulation-dependent quantity $E(i)$
  (like an error or an estimator), is itself a sequence defined by
  \begin{equation}
    \EOC E(i)=\frac{\log(E(i+1)/E(i))}{\log(h(i+1)/h(i))}
  \end{equation}
  Since the timesteps $\tau(i)$ are coupled to $h(i)$ this is well
  defined.

  Finally we look at the {\em effectivity index}, for each
  error-estimator pair, defined by
  \begin{align}
    \label{eqn:effectivity.index.lil2}
    \frac{\max_{n\in\fromto0m}\Norm{e(t_m)}}{
      \max_{n\in\fromto0m}\indrec_{\infty,n}
      +\tau(i)\sum_{n=1}^m\left(\indspace_{1,n}+\indtime_{1,n}\right) }
    \\
    \label{eqn:effectivity.index.l2h1}
    \frac{\left(\tau\sum\norm{e(t_m)}_1^2\right)^{1/2}}{
      \left(\tau\sum_{n=0}^m\indrec_{2,n}^2\right)^{1/2}
      +\tau(i)\sum_{n=1}^m\left(\indspace_{1,n}+\indtime_{1,n}\right) }
    .
  \end{align}
  The intial estimator is zero in this case and the data approximation
  and the mesh change estimators are dropped from this study.

  \begin{Obs}[effectivity index]
    We take all the constants involved in the estimators, including the
    interpolation constants, to be equal to $1$.  This, of course, is
    not true and a fine tuning of constants should be performed, but
    since our purpose here is mainly to check that the asymptotic
    behavior of the error and the estimator is the optimal one, the
    efficiency index is to be taken only qualitatively .  The fine
    tuning of constants is not in the scope of this paper.  Notice also
    that the time estimator, involving $\theta_{1,n}$, is in the
    denominators in \eqref{eqn:effectivity.index.lil2} and
    \eqref{eqn:effectivity.index.l2h1}.  This has the effect of making
    our effectivity index quite different than those used in the
    numerics for elliptic problems.
  \end{Obs}

  \subsection{Conclusions}
  The main conclusion of our numerical tests is that the estimators have
  the optimal rate of convergence which matches that of the error's
  norm.  It is important to notice that in order to observe the
  optimality of the estimators for different norms, different coupling
  of the meshsize $h$ and the stepsize $\tau$ must be chosen: for
  $\poly1$ elements, it is necessary to take $\tau\approx h^2$ to see
  that the $\leb\infty(\leb2)$ error norm has EOC $2$ while the
  $\leb2(\sobh1)$ norm has EOC $1$ (Figure \ref{fig:P1-tauh2}).  If the
  coupling $\tau\approx h$ is taken, for a problem where the time
  discretization errror dominates, such as (\ref{eqn:benchmark} fast),
  then both errors have EOC $1$ (Figure \ref{fig:P1-fast-tauh1}).

  The same observations are valid for tests with $\poly2$ elements,
  albeit the couplings are $\tau\approx h^3$ and $\tau\approx h^2$
  respectively in this case (see Figures \ref{fig:P2-tauh3} and
  \ref{fig:P2-fast-tauh2}).

}%

In a different article
\cite{LakkisMakridakisPryer:14:article:A-comparison} we conduct a more
thorough numerical experimentation, where mesh changes and data
approximation effects are included.

\newcommand{\figwidth}{.90\textwidth}
\begin{figure}[t]
  \begin{center}
    \includegraphics[width=\figwidth,trim= 0 0 0 0,clip]{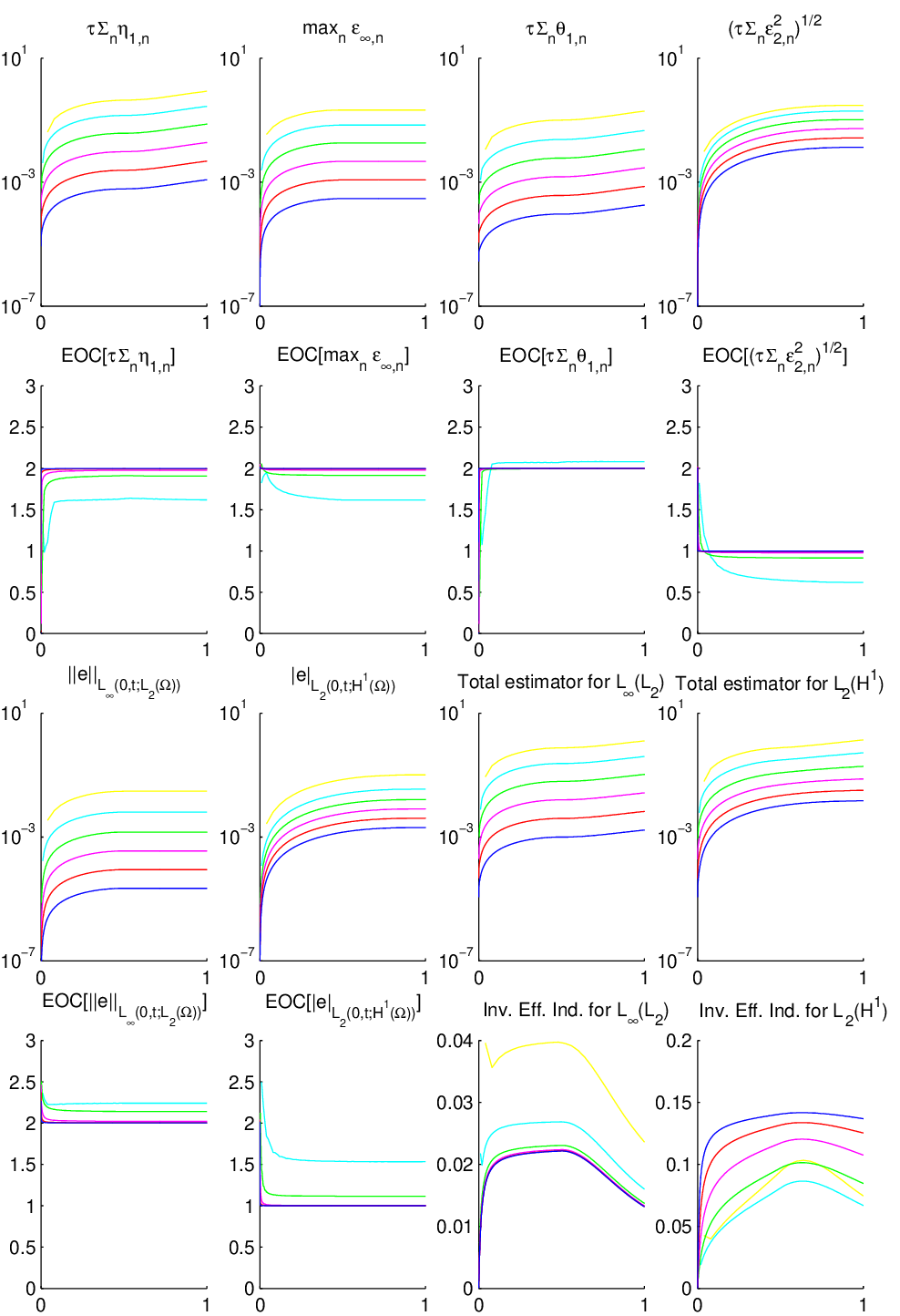}
  \end{center}
  \caption
      [Numerical results for slow problem, $\poly1$ elements and
        $\tau\approx h^2$] {
        \label{fig:P1-tauh2} 
        \changes{Numerical results for problem with exact solution
          (\ref{eqn:benchmark} slow) with $\poly1$ elements and $\tau\approx
          h^2$.  The abscissa represents time which ranges in $[0,1]$.  In the
          top-most row we plot the various estimators and in the second row we
          show the corresonding EOCs.  Notice that $\max\indrec_{\infty,n}$
          has EOC $2$ whereas $\left(\tau\sum\indrec_{2,n}^2\right)^{1/2}$ has
          EOC $1$. These are the leading terms in the total estimators (the 3
          and 4 plots in the 3rd row) and match in EOC, respectively, the
          $\leb\infty(\leb2)$ error and of the $\leb2(sobh1)$ error, as shown
          in the first 2 plots of the 3rd row and the 4th row.  Thus
          \eqref{eqn:lil2.estimate} and \eqref{eqn:l2h1.estimate} are seen to
          be sharp and optimal. The last two plots in the 4th row are the
          effectivity indexes for each norm.}
      }
\end{figure}
\begin{figure}
  \begin{center}
    \includegraphics[width=\figwidth,trim= 0 0 0 0,clip]{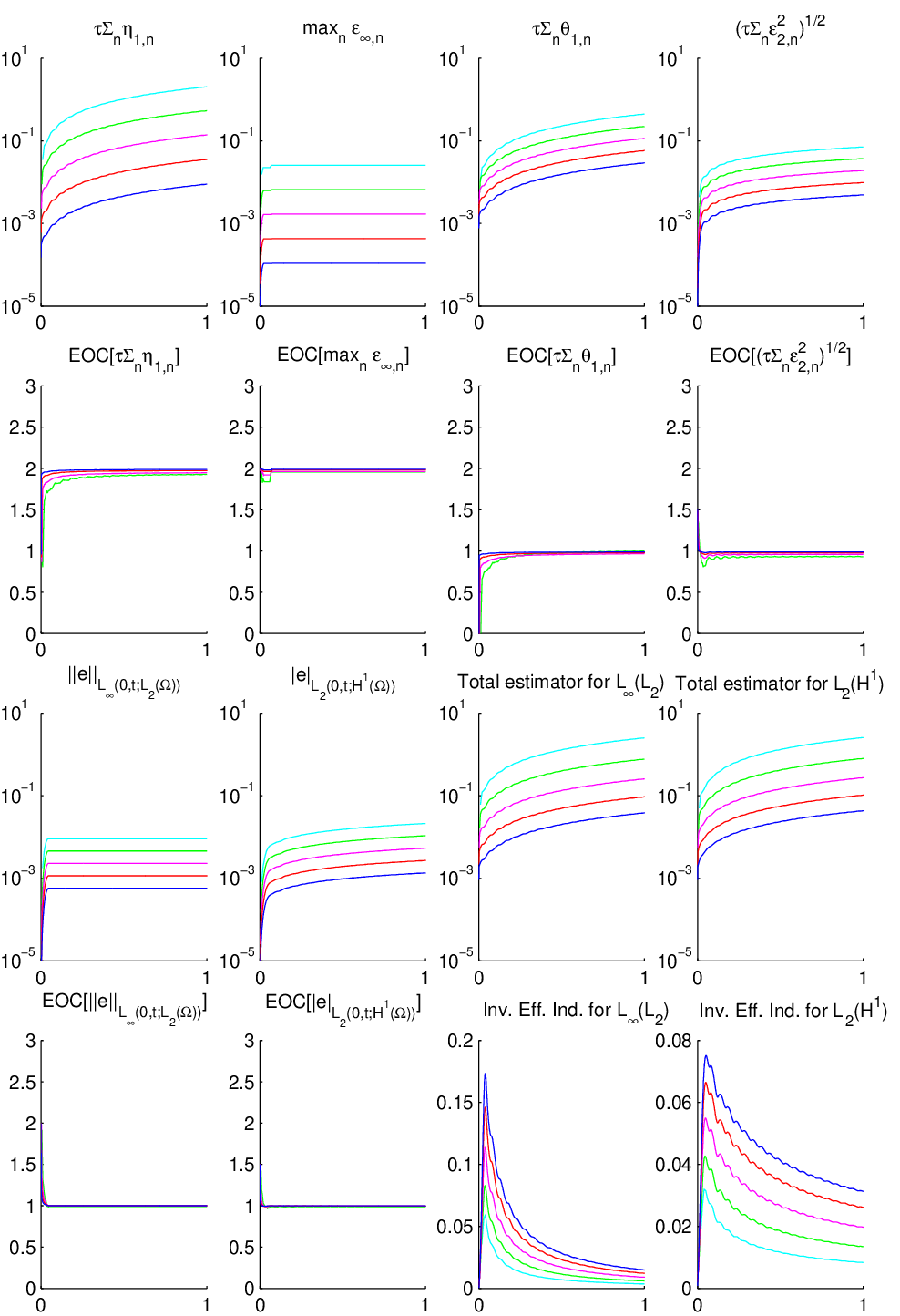}
  \end{center}
  \caption
      [Numerical results for fast problem, $\poly1$ elements and
        $\tau\approx h$]{
        \label{fig:P1-fast-tauh1}
        Simulation with $\poly1$ elements and $\tau\approx h$.
        Dominant time discretization error is created by taking the problem
        with fast time-oscillating exact solution (\ref{eqn:benchmark}
        fast).  The abscissa represents time which ranges in $[0,1]$.  In
        the top-most row we plot the various estimators and in the second
        row we show the corresonding EOCs.  Notice that
        $\tau\sum\indtime_{\infty,n}$ has EOC $1$, as opposed to $2$ in the
        previous example. This is now a leading term in both total
        estimators (plots 3 and 4 in the 3rd row) which have both EOC $1$.
        This EOC matches that of the $\leb\infty(\leb2)$ error and of the
        $\leb2(sobh1)$ error, as shown in plots 1 and 2 of the 3rd row and
        the 4th row.  Thus \eqref{eqn:lil2.estimate} and
        \eqref{eqn:l2h1.estimate} are both sharp, but only the second one is
        optimal due to the coupling. The last two plots in the 4th row are
        the effectivity indexes for each norm.
      }
\end{figure}
\begin{figure}
  \begin{center}
    \includegraphics[width=\figwidth,trim= 0 0 0 0,clip]{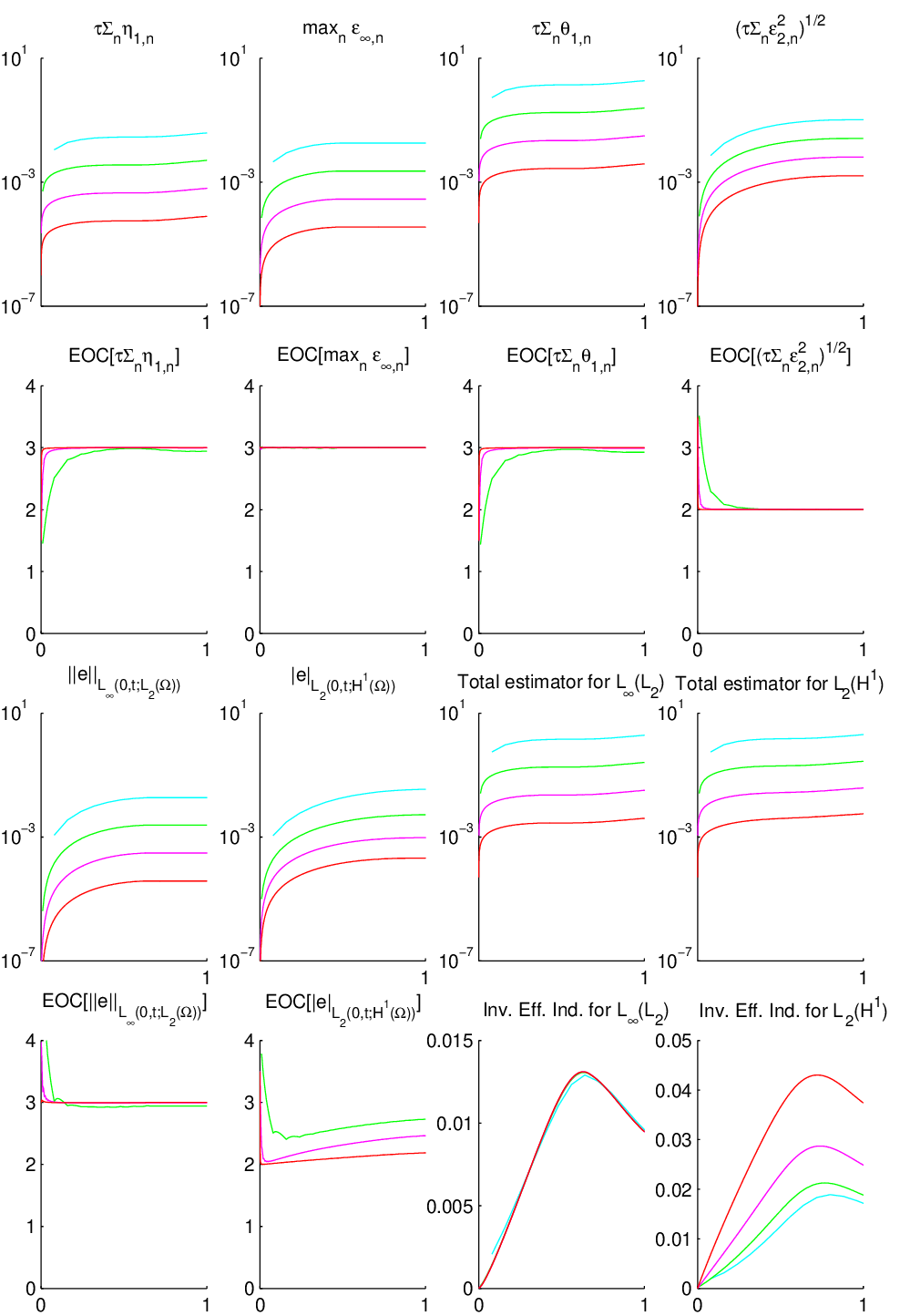}
  \end{center}
  \caption[Numerical results for slow problem, $\poly2$ elements and
    $\tau\approx h^3$] { 
    \label{fig:P2-tauh3}
    Numerical results for problem with exact solution (\ref{eqn:benchmark}
    slow) with $\poly1$ elements and $\tau\approx h^3$.  The abscissa
    represents time which ranges in $[0,1]$.  In the top-most row we plot
    the various estimators and in the second row we show the corresonding
    EOCs.  Notice that $\max\indrec_{\infty,n}$ has EOC $3$ whereas
    $\left(\tau\sum\indrec_{2,n}^2\right)^{1/2}$ has EOC $2$. These are
    the leading terms in the total estimators (the 3 and 4 plots in the
    3rd row) and match in EOC, respectively, the $\leb\infty(\leb2)$ error
    and of the $\leb2(sobh1)$ error, as shown in plots 1 and 2 of the 3rd
    and 4th rows.  Here the estimates \eqref{eqn:lil2.estimate} and
    \eqref{eqn:l2h1.estimate} are both sharp and optimal. The last two
    plots in the 4th row are the effectivity indexes for each norm.
  }
\end{figure}
\begin{figure}
  \begin{center}
    \includegraphics[width=\figwidth, trim= 0 0 0 0,clip]{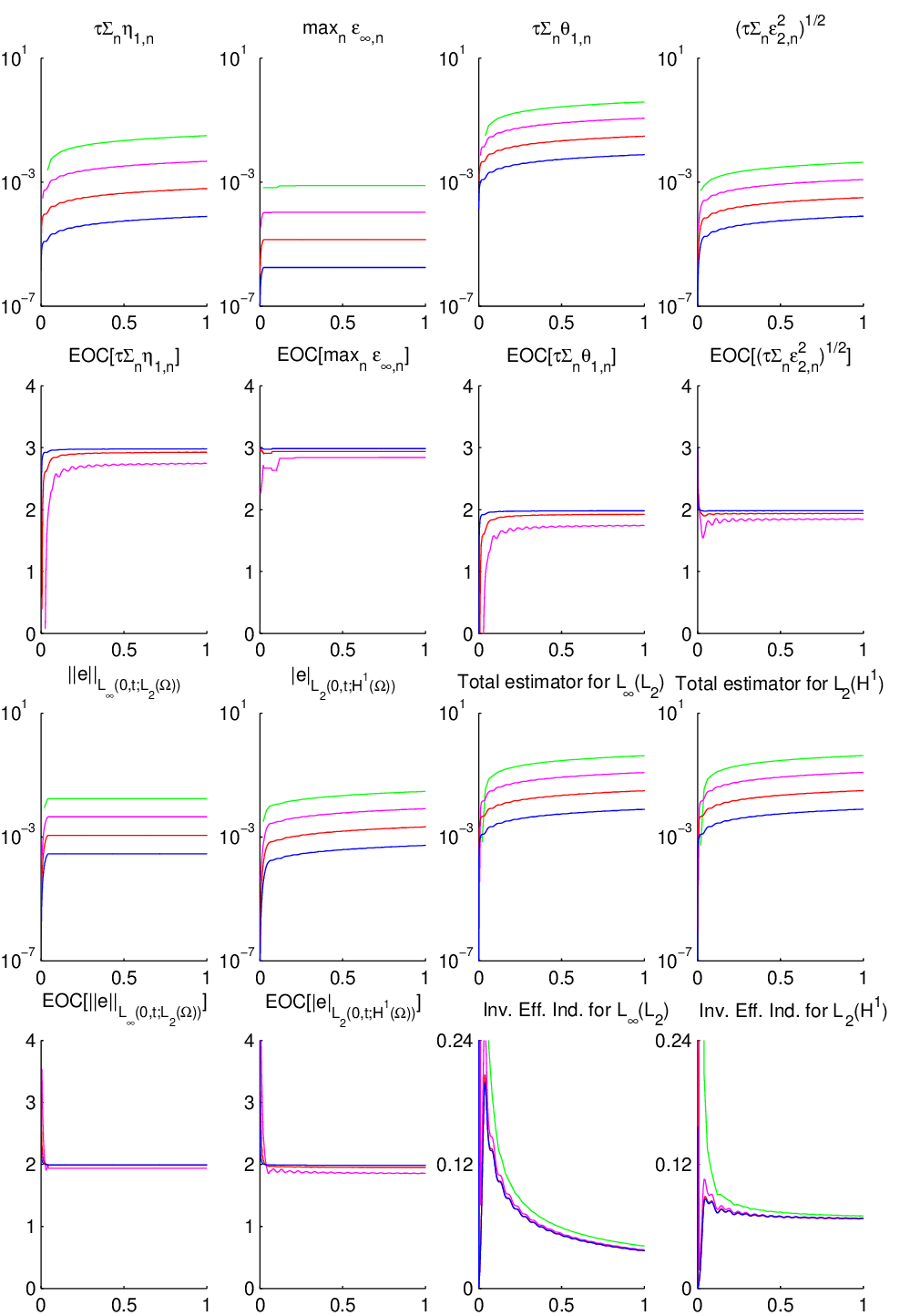}
  \end{center}
  \caption[Behavior of errors in time and effectivity indexes] {
    \label{fig:P2-fast-tauh2}
    Simulation with $\poly2$ elements and coupling $\tau\approx
    h^2$.  The time discretization error is dominant since exact is
    solution (\ref{eqn:benchmark} fast).  The abscissa represents time
    which ranges in $[0,1]$.  In the top-most row we plot the various
    estimators and in the second row we show the corresonding EOCs.
    Notice that $\tau\sum\indtime_{\infty,n}$ has EOC $2$, as opposed to
    $3$ in the previous example, because of the different coupling of
    the mesh and step sizes. The time estimator is now a leading term in
    both total estimators (plots 3 and 4 in the 3rd row) which have both
    EOC $2$.  This EOC matches that of the $\leb\infty(\leb2)$ error and
    of the $\leb2(sobh1)$ error, as shown in plots 1 and 2 of the 3rd
    row and the 4th row.  Thus \eqref{eqn:lil2.estimate} and
    \eqref{eqn:l2h1.estimate} are both sharp, but only the second one is
    optimal. The last two plots in the 4th row are the effectivity
    indexes for each norm.  }
\end{figure}
\appendix
\section{Compatible triangulations}
\label{sec:compatible.triangulations}
Each triangulation $\cT_n$, for $n\in\fromto1N$, is a refinement of a
{\em macro-triangulation} $\cM$ which is a triangulation of the domain
$\W$ that satisfies the same conformity and shape-regularity
\cite{BrennerScott:94:book:The-mathematical} assumptions made on its refinements in
\secref{sec:introduction}.  A refinement procedure is {\em admissible}
if it satisfies the following criteria:
\begin{NumberList}
\item the refined triangulation is conforming;
\item the shape-regularity of an arbitrarily deep refinement depends
  only on the shape-regularity of the macro-triangulation $\cM$;
\item if $\cT$ and $\cT'$ are both refinements, then for any two
  elements $K\in\cT$ and $K'\in\cT'$, then
  \begin{equation}
    K\cap K'=\emptyset\quad\text{or}\quad K\subset K'\text{or} \quad
    K'\subset K.
  \end{equation}
\end{NumberList}
Refinement procedures that satisfy these criteria exist.  For example,
the refinement by bisection described in the ALBERT manual
\cite{SchmidtSiebert:05:book:Design}, which is known to work for the space dimensions
$d=1,2,3,4$, is admissible for simplex triangulations.  All the
refinements by bisection of the macro-triangulation $\cM$ can be
stored in a single binary tree whose nodes represent a simplex.

We say that two triangulations are {\em compatible} if they are
refinements of the same macro-triangulation.  A set of compatible
triangulations can be endowed with a partial order relation: namely,
given two compatible triangulations $\cT$ and $\cT'$ we write
$\cT\leq\cT'$ if $\cT'$ is a refinement of $\cT$.  This partial
ordering permits to define in the natural way the {\em coarsest common
  refinement} of $\cT$ and $\cT'$, which we denote by
$\ccr{\cT}{\cT'}$, and the {\em finest common coarsening}, which we
denote by $\fcc{\cT}{\cT'}$.  An immediate property of these
definitions is
\begin{equation}
  \hat h=\maxi h{h'},\quad\text{and}\quad \check h=\mini h{h'}
\end{equation}
where $h,h',\hat h,\check h$ denote the meshsize of $\cT,\cT',
\ccr{\cT}{\cT'},\fcc{\cT}{\cT'}$ respectively.

\section{Inequalities}
\label{sec:inequalities}
\subsection{Elliptic regularity}
\label{sse:elliptic.regularity}
The estimates based on the space $\sobh{-1}(\W)$ or $\leb2(\W)$ norms
are based on \eqref{eqn:dual.norm}, duality arguments and the elliptic
regularity.  \changes{Whenever those estimates are used, we are
  tacitly assuming enough regularity on $a$ and $\W$, in such a way
  that for $k\in\NO$ there exists a constant $C_{2,k+2}=C_2[\W,k]$
  such that if $\phi\in\sobh k(\W)$ and $\psi\in\honezw$ are functions
  related by the (dual) elliptic problem}%
\begin{equation}
  \label{eqn:elliptic.duality}
  \abil\chi\psi=\ltwop\phi\chi,\quad\forall\chi\in\honezw,
\end{equation}
then
\begin{equation}
  \label{eqn:elliptic.regularity}
  \psi\in\sobh{k+2}(\W)\quad \text{and}\quad \norm{\psi}_{k+2}\leq
  C_{2,k+2}\norm{\phi}_k.
\end{equation}
\changes{For instance, this is true with $k=0$ if $\W$ is a convex
  polygonal domain and the coefficient matrix $\vec A$ in
  \eqref{eqn:bilinear.form} is a smooth space-function.  Whenever we
  consider $k=1$, or higher, it means that the domain $\W$ is smooth.
  We refer to \citet{Grisvard:85:book:Elliptic} for details about regularity on
  non-smooth domains.}%
\subsection{Interpolation inequalities}
\label{sse:scott.zhang.inequalities}
We will use the Clément-type interpolation operator
$\Pi^n:\honezw\rightarrow\fez n$ introduced by Scott \& Zhang
\cite{ScottZhang:90:article:Finite} which, under the needed regularity assumptions of
$\psi$ and finite element polynomial degree $\ell$, satisfies the
following interpolation inequalities for $j\leq\ell+1$:
\begin{gather}
  \label{eqn:scott.zhang.inner}
  \Norm{\hn^{-j}(\psi-\Pi^n\psi)}\leq C_{3,j} \norm{\psi}_j\\
  \label{eqn:scott.zhang.boundary}
  \Norm{\hn^{1/2-j}(\psi-\Pi^n\psi)}_{\Sigma_n}\leq C_{5,j}
  \norm{\psi}_j
\end{gather}
where the constants $C_{3,j}$ and $C_{5,j}$ depend only on the
shape-regularity, of the family of triangulations.
\subsection{Interpolation and mesh change}
\label{sse:mesh.change.interpolation.inequalities}
Since the triangulations $\cT_n$ and $\cT_{n-1}$ can be different when
adaptive mesh refinement strategies are employed, we introduce
$\hat\Pi^n$, the Clément--Scott--Zhang interpolator relative to the {\em
  finest common coarsening} of $\cT_n$ and $\cT_{n-1}$,
$\hat\cT_n:=\fcc{\cT_n}{\cT_{n-1}}$, whose meshsize is given by
$\hathn:=\maxi\hn\hno$.  Then the following inequality holds:
\begin{equation}
  \label{eqn:scott.zhang.coral} 
  \Norm{\hathn^{1/2-j}(\psi-\hat\Pi^n\psi)}%
  _{\Sigma_n\cup\Sigma_{n-1}\Tolto\Sigma_n\cap\Sigma_{n-1}} \leq
  C_{7,j} \norm{\psi}_j
\end{equation}
where the constant $C_{7,j}$ depends on the shape-regularity of the
triangulations and on the number of refinement steps (bisections)
necessary to pass from $\cT_n$ to $\cT_{n-1}$.
\subsection{Combining constants}
\label{sse:combining.constants}
\changes{ We often use combination of the constants introduced in this
  appendix or throughout the paper.  Since many constants appearing in
  theorems are products of basic constants our convention is that
  whenever a constant $C_{k,j}$ appears with the index $k$ being a
  non-prime integer then $C_{k,j}=C_{i,j}C_{l,j}$ where $k=il$.  If
  the index $k$ is a prime then the constant is a ``basic'' one and is
  defined in the text.  We recall that the constant $C_{2,1}$ is the
  Poincaré constant for the domain $\W$.  It is our hope that this
  convention simplifies the tracking of constants and does not
  complicate the reading.  }%
\bibliographystyle{plainnat}%

\end{document}